\documentclass[final,oneeqnum,onefignum,onetabnum]{siamltex}
\usepackage{amssymb}
\usepackage{amsmath}
\usepackage{amscd}
\usepackage{graphics}
\usepackage{epsfig}
\usepackage{url}
\usepackage{mdwlist}
\usepackage{tabls}
\usepackage{array}
\usepackage{algorithmic}
\usepackage{algorithm}
\usepackage{slashbox}
\usepackage{color}
\usepackage{multirow}

\usepackage{pgf}
\usepackage{tikz}
\usepgflibrary{arrows}
\usepackage{verbatim} % this is just for the multi line comments

\newtheorem{preremark}[theorem]{Remark}
{\begin{preremark} \upshape}{\end{preremark}}

\newcommand{\bm}[1]{\mbox{\boldmath${#1}$}}

\newcommand{\domain}{\Omega}
\newcommand{\cdomain}{\widebar{\Omega}}
\newcommand{\boundary}{\partial \domain}

\newcommand{\x}{\bm{x}}

\newcommand{\y}{\bm{y}}

\newcommand{\widebar}{\overline}
\newcommand{\drop}[1]{}

\newif\ifnever\neverfalse

\newif\iffullversion\fullversiontrue
\newif \ifearlyterm \earlytermtrue
\newcommand{\marginfix}{
\setlength{\parskip}{0.01cm}
\setlength{\textwidth}{6.0in}
\setlength{\oddsidemargin}{-0.0 in}
\setlength{\evensidemargin}{0.0 in}
\setlength{\topmargin}{-0.5in}
\setlength{\textheight}{9.0 in}
}

\marginfix

  \renewenvironment{thebibliography}[1]{%
    \begin{oldthebibliography}{#1}%
      \setlength{\parskip}{.3ex}%
      \setlength{\itemsep}{.3ex}%
  }%
  {%
    \end{oldthebibliography}%
  }

\begin{document}

%COMMENT REMOVED
%COMMENT REMOVED
%COMMENT REMOVED
%COMMENT REMOVED
%COMMENT REMOVED
%COMMENT REMOVED

%COMMENT REMOVED
%COMMENT REMOVED
%COMMENT REMOVED
%COMMENT REMOVED
%COMMENT REMOVED
%COMMENT REMOVED

\iffullversion
{\Large{\bf
\centerline{A parallel Heap-Cell Method for Eikonal equations.}
}}
\else
{\Large{\bf
\centerline{A parallel two-scale method for Eikonal equations.}
}}
\fi

\vspace*{.1in}
{\Large
\centerline{Adam Chacon and Alexander Vladimirsky\footnote{
\mbox{This research is supported in part by
the National Science Foundation grant DMS-1016150.}
\mbox{The first author's research was also supported in part
by Alfred P. Sloan Foundation Graduate Fellowship.}}
}

\vspace*{.1in}
\centerline{Center for Applied Mathematics and Department of Mathematics }
\centerline{Cornell University, Ithaca, NY 14853}
}

\vspace*{.1in}

\vspace*{.1in}

\begin{abstract}
\noindent
Numerous applications of Eikonal equations prompted the development of many efficient numerical algorithms.
The Heap-Cell Method (HCM) is a recent serial two-scale technique that has been shown to have advantages over other serial state-of-the-art solvers for a wide range of problems \cite{ChacVlad}.
This paper presents a parallelization of HCM for a shared memory architecture.
The numerical experiments in $R^3$ show that the parallel HCM exhibits good algorithmic behavior and scales well, resulting in a very fast and practical solver.
\iffullversion

We further explore the influence on performance and scaling of data precision, early termination criteria, and the hardware architecture.  A shorter version of this manuscript (omitting these more detailed tests) has been submitted to SIAM Journal on Scientific Computing \cite{ChacVlad_pHCM_short} in 2012.
\fi
\end{abstract}

\section{Introduction}
\label{s:intro}

The Eikonal equation is a nonlinear first order static PDE used in a range of applications, including robotic navigation, wavefront propagation, seismic imaging, optimal control, and shape-from-shading calculations.
The computational efficiency on a fixed grid is an important practical consideration in many of these applications.
%COMMENT REMOVED
Several competing ``fast'' serial algorithms have been introduced to %COMMENT REMOVED
solve the grid-discretized Eikonal equation in the last two decades.
The two most widely used among them are the Fast Marching Method (FMM) and the Fast Sweeping Method (FSM).  The Heap-Cell Method (HCM), introduced in the authors' previous work \cite{ChacVlad},  is a two-scale technique based on combining the ideas of FMM and FSM.  The current paper focuses on the parallelization of HCM for a shared memory architecture. We will start by briefly
describing the relevant discretization of the Eikonal PDE (section \ref{s:Eik_upwind}) and the prior algorithms for solving it (sections \ref{s:prior_serial_work} and \ref{s:prior_parallel}).  HCM is reviewed in section \ref{s:HCM}, and the new parallel HCM (pHCM) is described in detail in section \ref{s:pHCM}.  The numerical experiments in section \ref{s:experiments} demonstrate that pHCM is
%COMMENT REMOVED
efficient and achieves good parallel scalability for a wide range of grid resolutions and domain decomposition parameters.
%COMMENT REMOVED
\iffullversion
\else
Additional experimental results are included in an extended version of this manuscript \cite{pHCM_expanded}.
\fi
%COMMENT REMOVED
We conclude by discussing the limitations of our approach and the directions of future work in section \ref{s:conclusions}.

\subsection{Eikonal PDE and its upwind discretization}
\label{s:Eik_upwind}

%COMMENT REMOVED
%COMMENT REMOVED
An important subclass of Hamilton-Jacobi equations is formed by static Eikonal PDEs:
\begin{eqnarray}
\nonumber
| \nabla u (\x) | F(\x) &=& 1, \text{ on } \domain \subset R^d;\\
u(\x) &=& q(\x),  \text{ on } \boundary.
\label{eq:Eikonal}
\end{eqnarray}
Classical (smooth) solutions of (\ref{eq:Eikonal}) generally do not exist, and weak solutions are not unique \cite{BardiDolcetta}.
However, existence and uniqueness can be shown for the \emph{viscosity solution} \cite{CranLion}.  Moreover, the viscosity solution has an important natural interpretation as the value function of an isotropic time-optimal control problem:
$F$ can be viewed as a speed of motion, $q$ -- as an exit time-penalty on the boundary,
and $u(\x)$ -- as the minimum time-to-exit $\cdomain$ starting from $\x \in \domain.$
The gradient lines of $u$ can be interpreted both as the characteristics of the Eikonal PDE and as the optimal trajectories for the corresponding optimal control problem.

In this paper we will also consider slightly more general problems, where
exiting is only allowed through a closed nonempty ``exit set'' $Q \subset \boundary$,
with a prohibitively large exit time-penalty (e.g., $q = +\infty$) on $\boundary \backslash Q$.
This corresponds to a time-optimal control problem ``state-constrained'' to motion inside $\cdomain \backslash Q$,
with $u$ interpreted as a {\em constrained viscosity solution}  on $\cdomain$.
The boundary conditions on $Q$ are satisfied as usual (with $u=q$), while $\boundary \backslash Q$ is treated
as a non-inflow boundary, where the boundary conditions are ``satisfied in a viscosity sense'';
see \cite{BardiDolcetta}.
%COMMENT REMOVED
%COMMENT REMOVED
%COMMENT REMOVED
%COMMENT REMOVED
%COMMENT REMOVED
%COMMENT REMOVED
%COMMENT REMOVED
%COMMENT REMOVED
%COMMENT REMOVED
%COMMENT REMOVED

Several discretizations have been developed for equation \eqref{eq:Eikonal}, but here we focus on simple first-order upwind finite differences (similar to those presented in \cite{RouyTour}) on a 3D uniform Cartesian grid with stepsize $h$.
A typical gridpoint in $R^3$ will be denoted  $\x_{ijk} = (x_i, y_j, z_k) = (ih, jh, kh)$,
where $0 \leq i,j,k \leq n$ with $nh =1$ and $M = (n+1)^3$ is the total number of gridpoints in $\cdomain = [0,1]^3.$
We will use $U_{ijk}$ as a numerical approximation of $u(x_i, y_j, z_k)$,
with $U$ reserved to denote the entire grid solution vector.
For simplicity of exposition, we will assume that the exit set $Q$ is well-approximated on this grid,
and that all gridpoint values outside this computational cube are equal to $+\infty$.
Using the above notation, the upwind discretization can be written as
\begin{equation}
\begin{split}
\label{eq:Eik_discr}
\left( \max \left( D^{-x}_{ijk}U, \, -D^{+x}_{ijk}U, \, 0 \right)\right)^{2}
\; + \;
\left( \max \left( D^{-y}_{ijk}U, \, -D^{+y}_{ijk}U, \, 0 \right)\right)^{2} \\
\; + \;
\left( \max \left( D^{-z}_{ijk}U, \, -D^{+z}_{ijk}U, \, 0 \right)\right)^{2}
\; = \; \frac{1}{F^2_{ijk}},
\end{split}
\end{equation}

$$
\text{ where } \qquad
u_x(x_i, y_j, z_k) \approx D^{\pm x}_{ijk} U =
\frac{ U_{i \pm 1, j,k} - U_{i,j,k} }{ \pm h}, \qquad
\text{ etc.}
$$

Equation (\ref{eq:Eik_discr}) must hold at each gridpoint $(x_i, y_j, z_k) \in \cdomain \backslash Q,$ yielding a system of  coupled, nonlinear equations.  Since the set $Q \subset \boundary$ is
usually
lower dimensional, the total number of these equations is still $O(M)$.
If the %COMMENT REMOVED
values of $U$ at the neighboring gridpoints were available, the equation \eqref{eq:Eik_discr} could be
solved directly for $U_{ijk}$.
Since those neighboring $U$ values are not a priori known,
the resulting system can be solved iteratively (e.g., using Gauss-Seidel iterations), with $V_{ijk}$ denoting temporary values (or ``temporary labels'').
When these temporary gridpoint values stop changing, the iterative process terminates and $V \equiv U$.

However, the upwind nature of the discretization guarantees that not all neighboring values are relevant; i.e., only those neighboring values {\em smaller} than $U_{ijk}$ are actually needed to compute it from \eqref{eq:Eik_discr}. This is also a reflection of the Eikonal equation's \emph{causality property}, which is often exploited in the construction of fast algorithms.

\section{Prior Serial Methods}
\label{s:prior_serial_work}

The literature on serial methods for the Eikonal is vast; see \cite{ChacVlad} for a recent review.  Here we describe only those methods directly relevant to HCM and its parallelization.  For simplicity we describe these methods for Eikonal equations on Cartesian grids, but we note that some of them have been developed in much greater generality (\cite{BoueDupuis, KaoOsherQian,Li_1} and \cite{KimmSethTria,SethVlad1,SethVlad2}, for example).

The FSM \cite{TsaiChengOsherZhao, Zhao} solves system (\ref{eq:Eik_discr}) by Gauss-Seidel iterations with an alternating ordering of the gridpoints at each iteration.  These orderings, or ``sweep directions," are given by the standard loop orderings for Cartesian grids. E.g., in 2D these are

\begin{alignat*}{4}
i = 0,...,n &\qquad \qquad i = n,...,0 &\qquad \qquad i = 0,...,n &\qquad \qquad i = n,...,0\\
j = 0,...,n &\qquad  \qquad j = 0,...,n &\qquad \qquad j = n,...,0 &\qquad  \qquad j = n,...,0
\end{alignat*}

There are $2^d$ loop orderings in $d$ dimensions.
The efficiency of FSM is due to the fact that each characteristic of the solution can be divided into a finite number of contiguous portions where the characteristic directions in each portion are within only one quadrant.  Every $2^d$ sweeps all gridpoints along at least one of these portions will obtain their final values $U$.  The number of sweeps to convergence is related to $\rho$, the maximum number of times a characteristic changes direction from quadrant to quadrant.
Also, this number of sweeps is typically largely independent of the grid size as $h \rightarrow 0$ \cite{Zhao}, resulting in $O(\rho M)$ algorithmic complexity.  Unfortunately, even for fixed functions $F$ and $q$,
%COMMENT REMOVED
this $\rho$ is a priori unknown; moreover, $\rho$ is dependent on the grid orientation.

An additional speedup technique called the Locking Sweeping Method (LSM) \cite{Renzi} uses boolean ``active flags'' at each gridpoint to determine whether or not it will be updated in the next sweep.  (A value at a given gridpoint definitely will not change if none of its neighboring values have changed in the previous sweep.)  Initially only the gripdoints adjacent to the exit set are marked as ``active."
%COMMENT REMOVED
This technique reduces the total number of gridpoint computations per sweep but does not reduce the number of sweeps to convergence.

The FMM \cite{SethFastMarcLeveSet, SethSIAM} is a noniterative
method that uses the Eikonal equation's causality property to dynamically determine
an order to process the gridpoints.  A set $L$ of ``considered" gridpoints is maintained.  At each step of the algorithm, the considered gridpoint with the smallest temporary value is permanently ``accepted," and its not-yet-accepted neighbors are updated.
The set $L$ at each step can be regarded as an approximation to the current level set of $u$.
When $L$ is structured as a min-heap, updating it incurs an $O(\log(m))$ cost, where $m = |L|$.  The performance of FMM thus depends on the $(d-1)$-dimensional volume of the level sets of $u$; the upper bound for complexity is $O(M\log{M})$.   While the performance of FMM is more robust with respect to changes in the speed function $F$, domain geometry, and grid orientation, this method is not much faster for simpler problems, e.g., when the characteristics are straight lines (the regime where FSM is most efficient).

The HCM is a serial two-scale method that was introduced to combine the strengths of FMM and FSM on different scales.  On the coarse scale, a Fast Marching-like method is used to order the subdomains, and on the fine scale, sweeping (specifically LSM) is used on each subdomain separately; see section \ref{s:HCM} for a thorough description.  The informal motivation for this is that sufficiently zooming in on a portion of the domain reveals that  characteristics are approximately straight lines on that length scale, so sweeping restricted to that portion will converge quickly.  Even though the original purpose of the domain decomposition in HCM was to exploit the structure of the PDE serially, in this paper we show that the parallelization of the HCM is a natural byproduct and proves to be a very effective strategy; see section \ref{s:experiments}.

It is well known that many of the methods for Eikonal equations are directly related to prior algorithms for finding shortest paths on graphs\footnote{
Such graph algorithms are usually called {\em label-setting} and {\em label-correcting}.  To reflect this, we are using terms gridpoint ``value'' and ``label'' interchangeably.}.  Here we simply acknowledge this connection, but it is explored in detail in \cite{ChacVlad}.   In this framework, FMM and the prior non-iterative method in \cite{Tsitsiklis} are analogous to Dijkstra's method \cite{Diks}.
The ``Small Labels First / Large Labels Last" (SLF-LLL) \cite{Bertsekas_SLF} is another fast method originally designed for graph problems but also extended to Eikonal PDEs \cite{BertPolyTsi}; this iterative algorithm was designed to mimic the acceptance-order of nodes in Dijkstra's method while avoiding the costs associated with min-heap data structures.  Even though the worst-case complexity of SLF-LLL is not as good as that of Dijkstra's, in practice its performance can be better on many types of graphs. 

\section{Prior Parallel Methods}
\label{s:prior_parallel}

Several interesting approaches have been used to design parallel methods for Eikonal and related PDEs.
A careful performance/scalability comparison of \emph{all} such methods would be clearly valuable for practitioners but remains outside of scope of the current paper.  Here we give a brief overview of prior approaches primarily to put pHCM in context.  In section 6 we also use one of them as a benchmark for comparison with our own approach.

Two different parallelizations of FSM were introduced in \cite{ZhaoParallel}.  The first parallelization performs a domain decomposition and uses separate processors to run the serial FSM on each subdomain.  Subdomains are pre-assigned to processors and communication takes place along the shared boundaries.  The second approach presented in \cite{ZhaoParallel} does not use domain decomposition and performs all $2^d$ sweeps simultaneously on separate copies of the domain; these copies are then synchronized after each iteration by assigning the minimum value for each gridpoint.

The method of \cite{Detrixhe} is a more recent parallel sweeping technique (which we call ``Detrixhe Fast Sweeping Method" or DFSM) that utilizes the fact that, for the upwind scheme in 3D (eq. \ref{eq:Eik_discr}), gridpoints along certain planar slices through the computational domain do not directly depend on each other.  The planes are given by

\[
\alpha_i i +\alpha_j j + \alpha_k k = C,
\]

for $\alpha_i, \alpha_j,\alpha_k \in \{-1, 1\}$ and $C \in \mathbb{Z}$.  The choice of $\alpha$'s determines one of the $2^3$ sweeping directions; once the $\alpha$'s are fixed, the sweeping is performed by
%COMMENT REMOVED
incrementing $C$ (which corresponds to translating the plane in the sweep direction).  This is a \emph{Cuthill-Mckee} \cite{Saad} ordering of the gridpoints.
Inside any such plane the gridpoint updates are ``embarrassingly parallel'',
%COMMENT REMOVED
but the resulting method is {\em synchronous} since a barrier is required after processing each plane.
%COMMENT REMOVED
%COMMENT REMOVED
%COMMENT REMOVED
%COMMENT REMOVED
Unlike the methods in \cite{ZhaoParallel}, this algorithm requires exactly the same number of sweeps
as the serial FSM and also exhibits much better scalability.
This appears to be the current state-of-the-art in parallel sweeping methods for a shared-memory architecture;
%COMMENT REMOVED
thus, we have chosen to benchmark our results against it in section \ref{s:experiments}.
We note that a similar parallelization approach can also be used with the regular (lexicographic) gridpoint
ordering but with an appropriately extended stencil/discretization.
This idea was previously used in \cite{WeberDevirBronsBronsKimmel} for distance computations on parametric surfaces,
and more recently in \cite{Gillberg} to parallelize the sweeping
for more general (anisotropic) problems.

As for marching approaches, the canonical FMM is inherently serial (as is Dijkstra's method) and relies on a causal ordering of computations.
Several parallelizations of FMM have been developed employing fixed (problem-independent) domain decompositions  and running the serial FMM locally by each processor on preassigned subdomain(s) (e.g., \cite{Herrmann} and \cite{Breuss}).
In the absence of a strictly causal relationship between subdomains, this inevitably leads to erroneous gridpoint values,
which can be later fixed by re-running the FMM whenever the boundary data changes.
A very recent massively parallel implementation for distributed memory architecture in \cite{DonatSeth}
uses coarse grid computations to find a good subdomain preassignment, attempting to exploit non-strict causality to improve the efficiency; the approach is then re-used recursively to create a multi-level framework.

\begin{figure}
\begin{center}
$
\includegraphics[scale = .45]{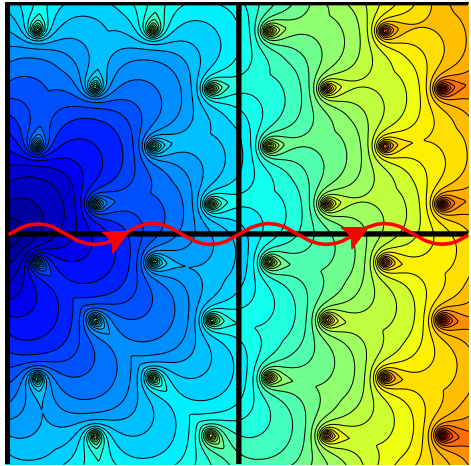}
$
\end{center}
\caption{Level sets for an Eikonal problem in 2D with cell boundaries %in black
and a characteristic
%COMMENT REMOVED
curve.
Since the characteristic repeatedly crosses the subdomain boundary, any method that solves this problem using the given domain decomposition will require a large number of iterations.}
\label{fig:mutual_dependence}
\end{figure}

The main difficulty with making the most effective use of a domain decomposition for the Eikonal equation is that the direction of information flow at subdomain boundaries is not known a priori.  If the domain is decomposed so that there is exactly one subdomain per processor, the loads may not be balanced. Additionally, a problem shared by all algorithms using a fixed domain decomposition is the existence of mutually dependent subdomains with a high degree of dependency; see Figure \ref{fig:mutual_dependence}.  Nevertheless, domain decomposition is often preferred as a %COMMENT REMOVED
parallelization approach %COMMENT REMOVED
to improve the cache locality and to avoid the use of fine-grain mutual exclusion.

A recent approach aims to minimize the inter-domain communications by creating problem-dependent causal domain decompositions. The so-called ``Patchy FMM" developed in \cite{CristianiPatchy} for feedback control systems uses coarse grid computations to build (almost) causal subdomains, which are then processed independently.  The disadvantages of this approach include complicated subdomain geometries, additional errors along subdomain boundaries, and frequent load balancing issues (since the causal subdomains are often very different in size).

In principle, it is also possible to parallelize some prior Eikonal solvers
(e.g., the Dial-like algorithm \cite{Tsitsiklis} and the Group Marching Method \cite{KimGMM})
without resorting to domain decompositions.  But we are not aware of any existing parallel implementations, and
the scalability is likely to be very limited due to the focus on gridpoint-level parallel computations.
For shortest path problems on graphs, examples of asynchronously parallelizable algorithms include the threshold method
and the SLF-LLL method \cite{Bertsekas_LLL}.
The idea in parallelizing the latter is to let each processor run a serial SLF-LLL method on its own local queue,
but with a heuristic used to determine which queue is to receive each graph-node tagged for updating.
A mutex is used for every node to prevent multiple processors from attempting to modify it simultaneously.
This parallel design inspired our own (cell-level) approach in the pHCM.

Several parallel algorithms were also developed for other computer architectures.  One method proposed in \cite{WeberDevirBronsBronsKimmel}, intended for SIMD and GPU architectures, computes shortest geodesic paths on parametric surfaces.  In this ``Parallel Marching Method" (PMM) the subdomains are processed serially with a dynamic ranking procedure similar to that of FMM.  Each time a subdomain is processed, the values of all gridpoints within it are updated using parallel ``raster scans," which are similar to the parallel sweeps in \cite{ZhaoParallel} and \cite{Detrixhe}.

Another method intended for massively parallel (SIMD GPU) architectures is the ``Fast Iterative Method" (FIM) developed in \cite{JeongWhitaker}.  In FIM, an unsorted list $L$ of active gridpoints is maintained, and at each iteration all gridpoints on $L$ are updated in parallel using Jacobi updates.
A variant, the ``Block FIM," maintains blocks of gridpoints on $L$, and all blocks on $L$ are updated in parallel.  New blocks are added based on whether any of their gridpoints received updates.  Blocks are used to take advantage of the SIMD parallelism.

\section{Heap-Cell Method}
\label{s:HCM}
To simplify the exposition, we give the following description of HCM in 2D.  The translation to higher dimensions is straightforward.
We first introduce some relevant notation: \\
$\bullet \,$ $X = \{ \x_1, ..., \x_M\} $, the grid (same as the grid used in FMM or FSM).
This single-subscript notation is meant to emphasize a gridpoint ordering, rather than the geometric
position indicated by the subscripts in formula \eqref{eq:Eik_discr}.
The corresponding gridpoint values are denoted as $V_i = V(\x_i).$\\
$\bullet \,$ $Q' = X \cap Q $, the set of ``exit gridpoints'', whose values are prescribed.\\
$\bullet \,$ $Z = \{c_1, ..., c_J\}$, the set of cells (or ``non-overlapping box-shaped subdomains'').\\
$\bullet \,$ $Q^{c} = \{ c \in  Z \mid  c \cap Q' \neq \emptyset \}.$\\
$\bullet \,$ $N(\x_j),$ the grid neighbors of $\x_j$; i.e., the gridpoints that exist to the north, south, east, and west of $\x_j$.\\
$\bullet \,$ $N^c(c_i)$, the set of neighboring cells of $c_i$; i.e., the cells that exist to the north, south, east, and west of
$c_i$.\\
$\bullet \,$ $N(c_i)$, the grid neighbors of $c_i$; i.e.,
$N(c_i) = \{ \x_j \in X \mid \x_j \not \in c_i \text{ and } N(\x_j) \bigcap c_i \neq \emptyset \}.$\\
$\bullet \,$ $V^{c}$, the cell label.\\
$\bullet \,$ $h_{x}^{c}$ and $h_{y}^{c}$ are the two cell dimensions (assume $h_{x}^{c} = h_{y}^{c} = h^{c}$).\\
$\bullet \,$ $r$, the number of gridpoints per cell-side.

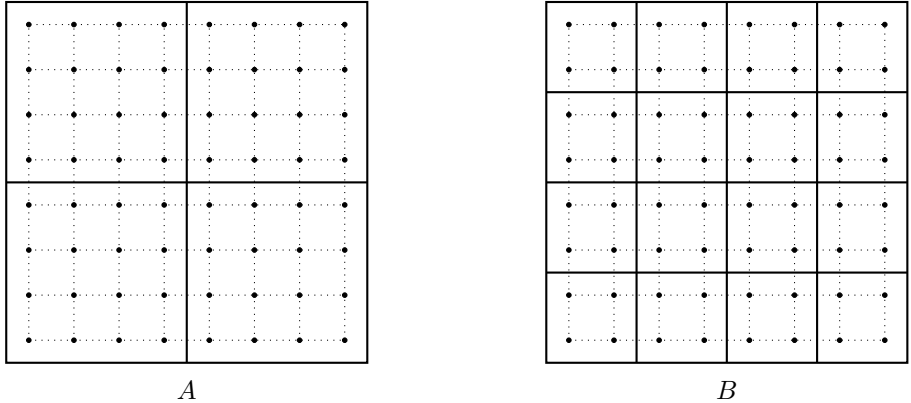
\begin{figure}

\center{
$
\begin{array}{cc}
\begin{tikzpicture}

\draw[dotted](0,0) rectangle +(4.2,4.2);
\foreach \x in {0.6, 1.2, 1.8, 2.4, 3.0, 3.6}
	\draw[dotted](\x,0)--(\x,4.2);
\foreach \y in {0.6, 1.2, 1.8, 2.4, 3.0, 3.6}
	\draw[dotted] (0, \y)--(4.2, \y);
\foreach \x in {0, 0.6, 1.2, 1.8, 2.4, 3.0, 3.6, 4.2}
    \foreach \y in {0, 0.6, 1.2, 1.8, 2.4, 3.0, 3.6, 4.2}
    {
        \filldraw (\x,\y) circle (.03 cm);
    }

\draw[thick](-.3,-.3) rectangle +(4.8,4.8);
\foreach \x in {2.1}
	\draw[thick](\x,-.3)--(\x,4.5);
\foreach \y in {2.1}
	\draw[thick](-.3,\y)--(4.5,\y);

%COMMENT REMOVED

\end{tikzpicture}&
\hspace{2cm}

\begin{tikzpicture}

\draw[dotted](0,0) rectangle +(4.2,4.2);
\foreach \x in {0.6, 1.2, 1.8, 2.4, 3.0, 3.6}
	\draw[dotted](\x,0)--(\x,4.2);
\foreach \y in {0.6, 1.2, 1.8, 2.4, 3.0, 3.6}
	\draw[dotted] (0, \y)--(4.2, \y);
\foreach \x in {0, 0.6, 1.2, 1.8, 2.4, 3.0, 3.6, 4.2}
    \foreach \y in {0, 0.6, 1.2, 1.8, 2.4, 3.0, 3.6, 4.2}
    {
        \filldraw (\x,\y) circle (.03 cm);
    }

\draw[thick](-.3,-.3) rectangle +(4.8,4.8);
\foreach \x in {.9,2.1,3.3}
	\draw[thick](\x,-.3)--(\x,4.5);
\foreach \y in {.9,2.1,3.3}
	\draw[thick](-.3,\y)--(4.5,\y);

%COMMENT REMOVED

\end{tikzpicture}\\

A& \hspace{2cm} B\\
\end{array}
$
}

\caption{
{\footnotesize
Two examples with different domain decompositions.  Both A and B are based on the same grid (dotted), with $M = 8^2$ and $h = 1/7$.  Figure A uses the cell size $h^c = 4/7$, the total number of cells $J = 2^2$, and $r = 4$ gridpoints per cell-side.  Figure B uses $h^c = 2/7$, $J = 4^2$, and $r = 2$.
}}

\label{allDomains}
\end{figure}

  To ensure that each gridpoint belongs to one and only one cell,
the cell boundaries are not aligned with gridlines, and $\Omega^c = \bigcup\limits_{j=1,...,J} \widebar{c_j}$ must be a superset of $\cdomain$; see Figure \ref{allDomains}.

%COMMENT REMOVED
%COMMENT REMOVED
%COMMENT REMOVED
%COMMENT REMOVED

\vspace{.2cm}

The original HCM, presented in \cite{ChacVlad}, is a serial two-scale method.
%COMMENT REMOVED
When the 2D analogue of system (\ref{eq:Eik_discr}) is solved on a cell $c$ (using any method), if the values of $N(c)$ are already correct, then all $\x_i \in c$ will receive their final values $U_i$.  Each cell is therefore dependent on a subset of $N^c(c)$, and the hyperbolic nature of the problem suggests that there is a \emph{preferred order} of processing the cells.
The motivation for employing a serial domain decomposition using sweeping on the grid is that, if the cell sizes are small enough, the characteristics within each cell will be approximately straight lines and sweeping will converge in very few iterations.

The HCM maintains a list $L$ of cells-to-be-processed, initially populated with $Q^c$.  The entire %COMMENT REMOVED
grid is initialized only once, in the same way as it is for LSM\footnote{That is, all $\x_i \not \in Q'$ have $V_i = \infty$; the active flags of gridpoints in \{$\x \in N(\x_i) | \x_i \in Q', \x \not \in Q'\}$ are set to ``active"; the active flags of all other gridpoints are set to ``inactive".}.
At each iteration of the main algorithm, a cell $c$ is chosen from $L$ and equation (\ref{eq:Eik_discr}) is solved by LSM on $X \cap c$ with the boundary conditions specified by the current values on $N(c_i)$.  The order of processing of the cells is determined dynamically based on heuristically assigned and updated cell values.  The name ``Heap-Cell" comes from organizing $L$ as a min-heap data structure.  HCM is designed to mimic FMM on the cell level (though previously processed cells may re-enter $L$; see Algorithm \ref{alg:HCM_main} for the pseudocode).
Since in typical cell-decompositions $J \ll M$, the cost of maintaining the heap $L$ is small compared to the cost of %COMMENT REMOVED
grid computations.  The experimental evidence in \cite{ChacVlad} shows that HCM is very efficient for a wide range of $M$ and $J$ values.

\begin{algorithm}[h]
\caption{Heap-Cell Method main loop.}
\label{alg:HCM_main}
\algsetup{indent=2em}
\begin{algorithmic}[1]
\STATE Initialize cell-values and grid-values
\STATE Add all $c \in Q^c$ cells to $L$
%COMMENT REMOVED
\WHILE{$L$ nonempty}
	\STATE Remove the cell $c$ with the smallest cell value from $L$
	\STATE $V^c(c) \leftarrow +\infty$
	\STATE Perform modified LSM on $c$ until convergence and populate\\ the list $DN$ of \emph{currently downwind} neighboring cells $\qquad$ //see Algorithm \ref{alg:gpUpdate}
	\FOR{each neighbor $c_k \in DN$}
		\STATE Update $V^c(c_k)$, the cell value of $c_k$
		\STATE Add $c_k$ onto $L$ if not already there
		\STATE Update the preferred sweeping directions of $c_k$
	\ENDFOR
\ENDWHILE

\end{algorithmic}
\end{algorithm}

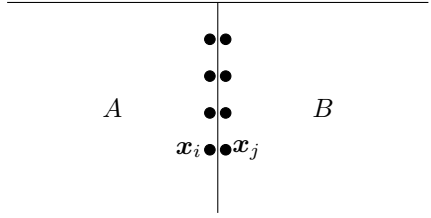
\begin{figure}[hhhh]
\center{
\begin{tikzpicture}
[scale = 0.7]
\draw[-] (0,0) -- (0,-4);
\draw[-] (-4,0) -- (4,0);

\draw(2,-2) node{$B$};
\draw(-2,-2) node{$A$};
\filldraw (.15,-.7) circle (.1 cm);
\filldraw (.15,-1.4) circle (.1 cm);
\filldraw (.15,-2.1) circle (.1 cm);
\filldraw (.15,-2.8) circle (.1 cm);

\filldraw (-.15,-.7) circle (.1 cm);
\filldraw (-.15,-1.4) circle (.1 cm);
\filldraw (-.15,-2.1) circle (.1 cm);
\filldraw (-.15,-2.8) circle (.1 cm);

\draw(.55,-2.8) node{$\x_j$};
\draw(-.55,-2.8) node {$\x_i$};

\end{tikzpicture}
}
\caption{
{\footnotesize
Suppose that as a result of processing the cell $A$ an eastern border value $V_i$ becomes updated.
If $V_i < V_j$ and $\x_j \not \in Q'$, then cell $B$ is \emph{downwind} of cell $A$.  In this case cell $B$ will be added onto $L$ unless already there, its value will be updated, and its preferred sweeping directions will be updated.
}}
\label{fig:cell_inflow_boundary}
\end{figure}

We say that a cell $B$ is \emph{currently downwind} from a cell $A$, if \
{\bf (1)} $A$ was the last processed cell and \ {\bf (2)}
there exist neighboring border gridpoints $\x_i \in A$ and $\x_j \in B$
such that the value of $V_i$ has changed the last time $A$ was processed and
\ {\bf (3)} $V_i < V_j$.
See Figure \ref{fig:cell_inflow_boundary}.
We note that, since this relationship is based on the temporary labels $V$,
it is entirely possible that the same $A$ might be also \emph{downwind} from $B$ at a different stage of the algorithm.
%COMMENT REMOVED
%COMMENT REMOVED

Unfortunately, a good dependency-ordering of cells may not exist even if we could base it on permanent gridpoint labels $U$ or even on the continuous viscosity solution $u(\x)$.
We will say that $B$ \emph{depends} on $A$ if there exists some optimal trajectory crossing the cell boundary from $B$ to $A$ on its way to $Q$.
This allows us to construct a dependency graph on the set of cells.  We will say that a cell-decomposition is \emph{strictly causal} if this dependency graph is acyclic.  A strictly causal decomposition ensures that there exists an ordering of cells such that each of them needs to be processed only once.

%COMMENT REMOVED

Figure \ref{fig:mutual_dependence} shows that, for many generic problems and large $h^c$, neighboring cells $A$ and $B$ are likely to be interdependent, resulting in multiple alternating re-processings of $A$ and $B$.
As $h^c$ decreases, the decomposition becomes \emph{weakly causal} - most cell boundaries become either purely inflow or purely outflow.
Additionally, if the ordering is such that most dependents are processed after the cells they depend on, the average number
of times each cell is processed becomes close to one.
As confirmed by the numerical evidence in \cite{ChacVlad}, weakly causal domain decompositions are very useful in decreasing the computational costs of serial numerical methods.

{\bf Processing cells by using Fast Sweeping Methods}:  Sweeping using LSM \cite{Renzi} is performed on the cell $c$ by using the neighboring %COMMENT REMOVED
grid values as boundary data.  Precisely, the domain for processing $c$ is $\tilde{c} = c \cup N(c)$, with the boundary conditions defined as
	 $\ \tilde{q}(\x_i) = q(\x_i)$  on $c \cap Q'$ and
	 $\ \tilde{q}(\x_i) = V_i$  on $N(c)$.
The sweeping processes gridpoints one at a time, with the
gridpoint update procedure detailed in Algorithm \ref{alg:gpUpdate}.

\begin{algorithm}[h]
\caption{Modified LSM update at a gridpoint $\x_i$.}
\label{alg:gpUpdate}
\algsetup{indent=2em}
\begin{algorithmic}[1]

\IF {$\x_i$ is inactive}
	\STATE Do nothing
\ELSE
	\STATE Set $\x_i$ inactive
	\STATE  Compute a possible new value $\widetilde{V}$ for $\x_i$ by solving equation \eqref{eq:Eik_discr}
	\IF{$\widetilde{V} < V(\x_i)$}
		\STATE $V(\x_i) \leftarrow \widetilde{V}$
		\FOR{each $\x_j \in N(\x_i) \backslash Q'$}
			\IF{$V(\x_j) > V(\x_i)$}
				\STATE Set $\x_j$ active			
				\IF {$\x_j$ is in a different cell from $\x_i$}
					\STATE Tag that cell as part of the list $DN$ of \emph{currently downwind} cells
				\ENDIF
			\ENDIF
		\ENDFOR
	\ENDIF
\ENDIF
\end{algorithmic}
\end{algorithm}

As in the usual LSM, we loop through different sweeping directions, using a new one in each iteration.
However, by the time a cell $B$ needs to be processed, the boundary information from its previously processed neighboring cells can be used to determine the preferred directions to start sweeping, with the likely effect of reducing the total number of sweeps needed to converge in $B$.
This is accomplished by having each cell maintain a list of boolean \emph{preferred-sweep-direction} flags, and by LSM beginning sweeping only from the directions marked {\tt TRUE}.  If the convergence is not achieved
after performing sweeps in these preferred directions we revert back to a standard loop
(i.e., in 2D the default standard loop would be SW, SE, NE, NW).
%COMMENT REMOVED
After a cell is processed, all sweep-direction flags are set
to {\tt FALSE}.  A sweep-direction flag of a cell $B$ is updated to {\tt TRUE}
only at the time a neighboring cell $A$ tags $B$ as \emph{downwind}.   The directions that are
updated depend on the location of $A$ relative to $B$.  For example, if $B$ is downwind from $A$ as in Figure \ref{fig:cell_inflow_boundary}, then both $A$-relevant sweep-direction flags
in $B$ (i.e., both NW and SW) will be set to {\tt TRUE}.

%COMMENT REMOVED
%COMMENT REMOVED
%COMMENT REMOVED
In principle,
the actual values of the border gridpoints could also be used to further restrict
the list of preferred sweep-directions (with the goal of avoiding unnecessary sweeping).
The {\em Fast Heap-Cell Method} introduced in \cite{ChacVlad} uses one such acceleration technique by checking the ``monotonicity'' of boundary data.  Since this technique is more costly in $R^3$, we are not using it in the current implementation of HCM.
%COMMENT REMOVED

{\bf Assigning Cell Values}:
Cell values are computed heuristically and intended to capture the direction of information flow.
If a cell $B$ depends on a cell $A$, then ideally $V^c (A) < V^c (B)$ should hold to ensure that $A$ is processed earlier.
We emphasize that the choice of a particular cell value heuristic {\bf does not} affect the final output of the HCM
(see \cite{ChacVlad} for a proof of convergence), but may affect the method's overall efficiency.
%COMMENT REMOVED
An ideal heuristic would reflect the inherent causal structure.  For example, if the cell decomposition is strictly causal, using a good cell-value heuristic would result in exactly $J$ heap removals.  For weakly causal cell decompositions (attained for all problems once $h^c$ becomes sufficiently small), a good cell-value heuristic ensures that the average number of heap removals per cell becomes closer to 1; see \cite{ChacVlad} and sections \ref{ss:serial_performance}, \ref{ss:parallel_performance} of the current paper for experimental evidence.

In this paper, our treatment of the cell value is different from the one in \cite{ChacVlad} in two ways:  1)  whenever a cell $B$ is removed from $L$, we reset $V^c(B)$ to $+\infty$, and 2) we assign $V^c(B)$ as the smallest of the newly updated gridpoint values in $N(B)$; see formula \eqref{eq:cell_value}. The logic is that cells should be ranked by the currently most upwind inflow.
We reset $V^c(B)$ so that if $B$ is to be processed again, the later time-of-processing will be determined only by new inflow information.  This heuristic appears to be very efficient for a variety of examples and easily generalizes to higher dimensions.  Most importantly, it seems to be effective at handling discontinuities in the speed function that do not align with the cell boundaries, which was a weakness of the cell value in \cite{ChacVlad}.

\begin{equation}
\begin{split}
\widetilde{V}^c(B) \leftarrow \min_{j \in A_{new}} V(\x_j)  \qquad
V^c(B) \leftarrow \min(V^c(B), \widetilde{V}^c(B)) \\\\
\end{split}
\label{eq:cell_value}
\end{equation}

\noindent where $A_{new}$ is the set of newly updated ``inflow for $B$" gridpoints of $A$  along the relevant cell border;
i.e., $A_{new} = \{ \x_i \in N(B) \cap A \, \mid \, \text{recently updated } U_i < U_j
\text{ for some } \x_j \in B \cap N(\x_i) \}.$
An efficient implementation of this heuristic %COMMENT REMOVED
relies on
updating the current minimum border value of $B$ at line 12 of Algorithm \ref{alg:gpUpdate}.
%COMMENT REMOVED
The equivalent of formula \eqref{eq:cell_value} was also previously used to
determine the processing order of large ``charts'' in \cite{WeberDevirBronsBronsKimmel}.

Finally, we use a natural initialization of cell values before the main loop of the algorithm:
$$
V^c(c) \; \leftarrow \;
\begin{cases}
\min \{ V(\x_j) \, | \, \x_j \in c \cap Q' \}, & \text{ if } c \in Q^c;\\
+\infty, & \text{ otherwise.}
\end{cases}
$$
%COMMENT REMOVED
%COMMENT REMOVED
\vspace{.4cm}

\iffullversion
\noindent
Performance comparison with the heuristic originally used in \cite{ChacVlad} is also included in section \ref{ss:otherCellVals}.
%COMMENT REMOVED
\fi 

\section{Parallelization}

\label{s:pHCM}

There are several different approaches one can take to parallelize HCM.  It is possible, for instance, to parallelize the sweeping scheme within an individual cell.  Our choice for pHCM was to have multiple subdomains processed simultaneously.
Each processor $p$ essentially performs the serial HCM on its own local cell-heap $L_p$, but with one important difference:  when a cell $c$ is tagged for re-processing,
we attempt to add it to the heap $L_j$ with the lowest current number of cells.
Except for some modifications explicitly described below,
most of the subroutines of the serial HCM can be directly reused in pHCM as well.  In algorithm \ref{alg:pHCM}, all data is shared unless stated otherwise.

\begin{algorithm}[h]
\caption{Parallel Heap-Cell Method pseudocode.}
\label{alg:pHCM}
\algsetup{indent=2em}
\begin{algorithmic}[1]
\STATE Cell Initialization: same as in HCM (divide cells $Q^c$ evenly among all heaps $L_p$)
\STATE Fine Grid Initialization: same as in HCM
\STATE $P \leftarrow$ number of threads
\STATE {\it activeCellCount} $\leftarrow |Q^c|$
%COMMENT REMOVED
\STATE PARALLEL SECTION
%COMMENT REMOVED
\WHILE{{\it activeCellCount} $> 0$ }
	\WHILE{$L_{p}$ is nonempty}
		\STATE Lock heap $L_{p}$
		\STATE \hspace{\algorithmicindent} Position-lock cell $c$ at the top of $L_{p}$
		\STATE \hspace{\algorithmicindent} \hspace{\algorithmicindent} Remove $c$ from $L_{p}$
		\STATE \hspace{\algorithmicindent} \hspace{\algorithmicindent} $V^c(c) \leftarrow +\infty$
		\STATE \hspace{\algorithmicindent} Position-unlock $c$
		\STATE Unlock $L_{p}$
		\STATE Compute-Lock $c$
		\STATE \hspace{\algorithmicindent} Perform modified LSM on $c$ and populate the (local) list $DN$
        \STATE \hspace{\algorithmicindent} of \emph{currently downwind} neighboring cells $\qquad$ //see Algorithm \ref{alg:gpUpdate}
		\STATE \hspace{\algorithmicindent} Set all preferred sweeping directions of $c$ to {\tt FALSE}
		\STATE Compute-Unlock $c$

		\FOR{each $c_k \in DN$}
		\vspace{.05cm}
			\STATE Compute a possible new (local) cell value $\widetilde{V}$ for $c_k$
			\vspace{.1cm}
			\IF{$\widetilde{V} < V(c_k)$}
				\STATE {\bf Set Cell Value ($c_k, \widetilde{V}$)} $\qquad$ //see Algorithm \ref{alg:set_cell_val}
			\ENDIF

			\IF{$c_k$ is not on a heap}
			\STATE {\bf Add Cell ($c_k)$} $\qquad$ //see Algorithm \ref{alg:add_cell}
	
			\ENDIF
			\STATE Update sweeping directions of $c_k$ based on location of $c$
		\ENDFOR
		\STATE {\it activeCellCount} $--$   $\qquad$ {\bf (atomic)}

	\ENDWHILE
	
%COMMENT REMOVED

\ENDWHILE
\end{algorithmic}
\end{algorithm}

%COMMENT REMOVED

\begin{algorithm}[h]
\caption{Set Cell Value ($c_k, \widetilde{V}$).}
\label{alg:set_cell_val}
\algsetup{indent=2em}
\begin{algorithmic}[1]

\STATE {\it success} $\leftarrow $ {\tt FALSE}

\WHILE{{\it success} == {\tt FALSE}}

\IF {$c_k$ is not on a heap}
	\STATE Position-lock $c_k$
		 \IF  {$c_k$ is still not on a heap}
		\STATE  $V(c_k) \leftarrow \min (\widetilde{V},V(c_k))$
		  \STATE    {\it success} $\leftarrow$ {\tt TRUE}
		 \hspace{\algorithmicindent} \ENDIF
	\STATE  Position-unlock $c_k$
\ELSE
	\STATE $j \leftarrow $ index of the heap of $c_k$
	\STATE Lock $L_j$
	\STATE Position-lock $c_k$
	\IF{$c_k$ is still on $L_j$}
		\STATE  $V(c_k) \leftarrow \min (\widetilde{V},V(c_k))$
	\STATE Heap-sort $L_j$
	\STATE   {\it success} $\leftarrow$ {\tt TRUE}
	\ENDIF
	\STATE Position-unlock $c_k$
	\STATE Unlock $L_j$

\ENDIF

\ENDWHILE
\end{algorithmic}
\end{algorithm}

\begin{algorithm}[h]
\caption{Add Cell ($c_k)$.}
\label{alg:add_cell}
\algsetup{indent=2em}
\begin{algorithmic}[1]
\STATE $j \leftarrow $ index of heap with fewest elements (no locking; counts may be outdated during search);

\STATE {\it testCount} $\leftarrow$ 0
\WHILE {Lock $L_{(j + \it{testCount}) \% P }$ can not be immediately obtained}
	\STATE {\it testCount}++
\ENDWHILE

\STATE

\STATE  Position-Lock $c_k$
\IF{$c_k$ is still not on a heap}
	\STATE  Add $c_k$ onto  $L_{(j + \it{testCount})\% P}$
	\STATE {\it activeCellCount} ++   $\qquad$ {\bf (atomic)}
%COMMENT REMOVED
\ENDIF
\STATE Position-Unlock $c_k$

\STATE Unlock $L_{(j + \it{testCount})\%  P}$

\end{algorithmic}
\end{algorithm}

The described algorithm gives rise to occasional (benign) data race conditions.
%COMMENT REMOVED
But before explaining why they have no impact on correctness/convergence,
we highlight several main
design decisions:
%COMMENT REMOVED
%COMMENT REMOVED
\begin{itemize}

\item To ensure efficiency/scalability, there is no synchronization mechanism at the gridpoint level.

\item Unlike many other parallel Eikonal solvers, pHCM is {\em asynchronous}; i.e.,
no barriers are used to %COMMENT REMOVED
ensure that cells are processed in some specific order.

\item There are two separate ``individual cell operations'' that must be serialized: %COMMENT REMOVED
     1) the movement of a cell onto/ off / within a heap and 2) the update of gridpoint values within that cell.  However, both of these %COMMENT REMOVED
    can be safely performed simultaneously.  Thus, each cell maintains both a ``compute" lock and a ``position" lock to allow for the overlapping of these operations.

\item Adding a cell onto the heap with fewest elements ensures good load balancing.
But if that heap is currently locked, waiting for the lock to be released might have the opposite effect on the method's performance.  Since we can assign the cell to another heap without drastically altering the balance, we attempt to obtain the lock using the {\bf \tt omp\_test\_lock} subroutine, and move on to the next heap if that attempt was unsuccessful; see Algorithm \ref{alg:add_cell}.
Profiling shows that this approach always results in better performance
than using the {\bf \tt omp\_set\_lock}.

\item The {\it activeCellCount} is decremented on line 29 of Algorithm \ref{alg:pHCM} (rather than around line 10) to prevent other threads from quitting prematurely.

\item The cell update (lines 15-17 of Algorithm \ref{alg:pHCM}) is exactly the same sweeping procedure as in HCM.  Just as in HCM, any other method that solves system \eqref{eq:Eik_discr} within a cell $c$ may be substituted in place of LSM.  However, if the
    grid-value
    updates
    inside $c$
    also involve updating any grid-level data in $N^c(c)$, the potential race conditions must be handled carefully.
    %COMMENT REMOVED
    Below we explain how this issue is handled in LSM for the active flag updates across cell-boundaries.

\end{itemize}

%COMMENT REMOVED
\subsection{Efficiency and data race conditions}
\label{ss:phcm_features}

%COMMENT REMOVED

There is always a delicate trade-off between performance-boosting heuristics in the serial realm and the synchronization penalty they would incur in the parallel %COMMENT REMOVED
implementation.  The serial HCM has several features (the use of LSM within cells, the use of preferred sweeping directions, the accuracy of cell values at predicting information flow) that
%COMMENT REMOVED
could cause contention when parallelized.  In this section we describe how we chose to handle those features in designing pHCM.  Since there is no synchronization at the gridpoint level, we have actually allowed several data races to be present in the algorithm.  We first check the convergence of the algorithm in the presence of these data races.

For all of the following arguments we assume a {\em sequentially consistent} memory model, meaning that the instructions in Algorithm \ref{alg:pHCM} are executed in the order they appear.  On modern platforms it is possible that compilers or hardware will reorder the program's instructions.  While these optimizations are innocuous in serial codes, in a multi-threaded environment this can lead to unexpected results\footnote{Indeed, in our implementation it was actually necessary to explicitly
%COMMENT REMOVED
prevent such reordering of certain lines of code (using Open MP's ``flush" pragma).  %COMMENT REMOVED
}.

%COMMENT REMOVED

Consider first a more basic version of pHCM that uses FSM within cells instead of LSM.  There is still a possibility of data races along the boundary of each cell: updating a border gridpoint by Eq. \eqref{eq:Eik_discr} requires reading information in a neighboring cell.  But it is easy to see that the monotonicty of gridpoint value updates makes such data races harmless.
%COMMENT REMOVED
Suppose two cells $A$ and $B$ are being simultaneously swept by processors $p_A$ and $p_B$ respectively (see Figure \ref{fig:cell_inflow_boundary}). Suppose also that $B$ undergoes its final sweep.  First, the most obvious outcome is that

a.  $p_A$ updates $\x_i$ (and writes $V_i$).

b.  $p_A$ checks $V_j$ and finds $V_i < V_j$, $\Rightarrow$ tags $B$ to be added onto a heap.

\noindent So, $B$ will have a chance to use the new boundary information $V_i$ the next time it is processed.  Now, suppose neighbors $\x_i$ and $\x_j$ are updated simultaneously (i.e., Algorithm \ref{alg:gpUpdate} is executed in parallel at $\x_i$ and $\x_j$ by the different processors).  Suppose also that the final sweep in $A$ leaves $V_i < V_j$.  Then either

a.  $p_A$ writes $V(\x_i)$.

b.  $p_B$ writes $V(\x_j)$.

c.  $p_A$ checks $V_j$ and finds $V_i < V_j$, $\Rightarrow$ tags $B$ to be added onto a heap.

d.  $p_B$ checks $V_i$ and finds $V_i < V_j$, $\Rightarrow$ does nothing.

\noindent or

a.  $p_B$ writes $V(\x_j)$ .

b.  $p_B$ checks $V_i$ and finds $V_j < V_i$, $\Rightarrow$ tags $A$ to be added onto a heap.

c.  $p_A$ writes $V(\x_i)$.

d.  $p_A$ checks $V_j$ and finds $V_i < V_j$, $\Rightarrow$ tags $B$ to be added onto a heap.

\noindent In the latter case the cell $A$ is unnecessarily added onto a heap, but this redundancy does not impact the convergence.  Therefore, a cell with new inflow boundary information %COMMENT REMOVED
is always
guaranteed to be reprocessed at some later point.

But our reliance on the Locking Sweeping technique introduces an additional issue: it is also necessary to ensure that all relevant boundary gridpoints in that yet-to-be-reprocessed cell will be marked as ``active'' -- since otherwise the first cell-sweep will not touch them.  Recall that $p_A$ will only set the gridpoint values within $A$, but because of LSM, it might also change the active flags of gridpoints in $N(A) \cap B$.
What if $\x_i$ and $\x_j$ are updated simultaneously, $p_A$ makes $\x_j$ active, but $p_B$ immediately resets it as inactive
and $V_j$ is never recomputed based on the new value of $V_i$?
The order of operations in Algorithm \ref{alg:gpUpdate} makes this scenario impossible, since setting a gridpoint inactive is immediately {\em followed} by the re-computation of that gridpoint's value.

Finally, there is an one additional design choice we have made that causes a race condition at the cell-level when setting the cell's preferred sweeping direction flags.
After processing a cell $A$, we typically need to update the preferred sweeping directions of its neighboring cells.
If one of these neighboring cells $B$ is simultaneously processed using LSM, the preferred directions data might be overwritten. We could avoid this scenario by obtaining $B$'s computation lock before updating its preferred directions.
Our implementation does not use this idea because the preferred directions only reduce the number of sweeps without affecting the convergence, and because the additional contention would dominate the savings for most $M/J$ ratios.  Since all other access to cell-level data is lock protected, pHCM converges.

%COMMENT REMOVED
%COMMENT REMOVED
%COMMENT REMOVED
%COMMENT REMOVED
%COMMENT REMOVED

%COMMENT REMOVED
%COMMENT REMOVED
%COMMENT REMOVED
%COMMENT REMOVED
%COMMENT REMOVED
%COMMENT REMOVED
%COMMENT REMOVED
%COMMENT REMOVED
%COMMENT REMOVED
%COMMENT REMOVED
%COMMENT REMOVED
%COMMENT REMOVED
%COMMENT REMOVED
%COMMENT REMOVED
%COMMENT REMOVED
%COMMENT REMOVED
%COMMENT REMOVED
%COMMENT REMOVED

\section{Numerical Experiments}
\label{s:experiments}

In this section we present and compare the performance of FMM, FSM, LSM, HCM, DFSM (a parallel sweeping method), and pHCM on three qualitatively different examples.
Our primary goal is to test the ``strong scalability" of pHCM with various cell decompositions.
Sections \ref{ss:serial_performance} and \ref{ss:parallel_performance} provide a more detailed performance analysis of the serial and parallel methods respectively.
Our source code and scripts for all methods and examples in this paper are publicly available from
\url{http://www.math.cornell.edu/~vlad/papers/pHCM/}.
\\

\noindent {\bf Benchmark problems}

We consider three Eikonal examples %COMMENT REMOVED
with an exit set $\{ (0.5, 0.5, 0.5) \}$
on a unit cube domain $\cdomain = [0,1]\times[0,1]\times[0,1]$.
In all three cases, the boundary conditions are $q=0$ in the center and $q=+\infty$ on the boundary of the cube.
Since the center of the computational domain is not a gridpoint (i.e., $M$ is even), we have initialized $U$ on the set $Q$ of the 8 gridpoints closest to the center.  Since $J$ values are also even, the set $Q^c$ contains 8 cells in all of the examples.

The speed functions are:

\begin{enumerate}
\item $F \equiv 1$.

\item $F(x,y,z) = 1 + .5\sin{(20\pi x)}\sin{(20\pi y)}\sin{(20\pi z)}$.

\item $F(x,y,z) = 1 + .99\sin{(2\pi x)}\sin{(2\pi y)}\sin{(2\pi z)}$.

\end{enumerate}

\vspace{.2cm}
These examples are ``representative'' in the sense that their respective viscosity solutions are qualitatively very different.
%COMMENT REMOVED
In example 1, all characteristics are straight lines.  In example 2, the characteristics are highly oscillatory and might weave through cell boundaries many times.  The third example has more moderate behavior, with curved characteristics that do not oscillate rapidly.\\
%COMMENT REMOVED

\noindent {\bf Experimental setup and implementation details}

All experiments
\iffullversion
(except for those in subsection \ref{ss:octopus})
\fi
were performed on the Texas Advanced Computing Center's ``Stampede" computer, using a single Dell PowerEdge R820 node with four E5-4650 8-core 2.7 GHz processers and 1TB of DDR3 memory.
We implemented all methods in C++ and compiled with the {\tt -O2} level of optimization using the Intel Composer XE compiler v13.0.  All solutions
\iffullversion
(except for those in subsection \ref{ss:single_prec})
\fi
were computed and stored using \emph{double} precision.  The speed $F(x,y,z)$ was computed by a separate function call as needed, instead of precomputing and storing it for every gridpoint.  HCM and pHCM use Locking Sweeping, which is experimentally always much faster than regular Fast Sweeping.
\iffullversion
\else
In all iterative methods,  the sweeps were continued as long as some gridpoints received updated values.
\fi
In benchmarking all parallel methods,
we have used one thread per core, up to a total of 32 cores. %COMMENT REMOVED
In addition, for some $r$ values, the performance of pHCM may be significantly influenced by both system-level background processes
and variations in the effective speed of the cores.  To fully reflect this, each pHCM test was performed
30 times and we report both the median values and the max/min ``error bars".

We compare our methods' performance/scaling to a parallelization of the sweeping methods.  Our implementation largely follows the method described in \cite{Detrixhe}, but with two exceptions:

\begin{itemize}
\item Detrixhe et. al. have not tested a ``locking sweeping'' version of their method; our implementation of DLSM is based on a straightforward substitution of LSM-updates for FSM-updates.

\item Our implementation of DFSM and DLSM use the default Open MP static loop scheduling (``omp for") to divide the work amongst threads instead of the manual load balancing procedure described in \cite{Detrixhe}.
\end{itemize}

\iffullversion
In all iterative methods,  the sweeps were continued as long as some gridpoints received updated values; in subsection \ref{ss:earlyTerm} we separately investigate the performance improvements due to an ``early termination''.
In subsection \ref{ss:single_prec} we explore the influence of memory footprint by storing/computing values in \emph{single} precision.
In subsection \ref{ss:octopus} we provide additional benchmarking results on a different shared memory architecture.
Subsections \ref{ss:checkerboard} and \ref{ss:shellMaze} contain results for additional examples (with piecewise-constant $F$).
Finally, in subsection \ref{ss:otherCellVals} we provide data for performance with a different cell value heuristic.
\else

\vspace*{1mm}

\noindent
We have also conducted many additional numerical tests including
\begin{itemize}
\item benchmarking with {\em single} precision data/computations;
\item benchmarking with {\em ``early termination criterion''} for sweeping methods;
\item benchmarking on a different shared memory architecture;
\item benchmarking with a different cell value heuristic;
\item additional test problems (with piecewise-constant $F$).
\end{itemize}
The results (included in the expanded version of this paper
\cite{pHCM_expanded})
%COMMENT REMOVED
are sufficiently similar qualitatively, and we omit them here for the sake of brevity.
\fi

%COMMENT REMOVED

\vspace*{2mm}
\noindent {\bf Layout of experimental results}

The HCM tests were run using $J$ = $M/2^3, M/4^3, M/8^3, M/16^3$, and $M/32^3$, so there are $r$ = 2/4/8/16/32 gridpoints per cell side.    ``HCM$r$" and ``pHCM$r$" in the legends mean HCM and pHCM with $J = M/r^3$.
(This notation emphasizes the amount of work per cell, but it is different from the format previously used in \cite{ChacVlad, ChacVlad_expanded},
where the table headings directly stated $J$ rather than $r$.)
On each test problem the performance of pHCM depends on 3 problem parameters: $M$, $r$, and $P$, the number of processors.  The performance/scaling plots for pHCM2 are omitted to improve the readability of all figures.

\begin{figure}
\hspace{-.6 in}
$
\begin{array}{ccc}
\includegraphics[scale = .45]{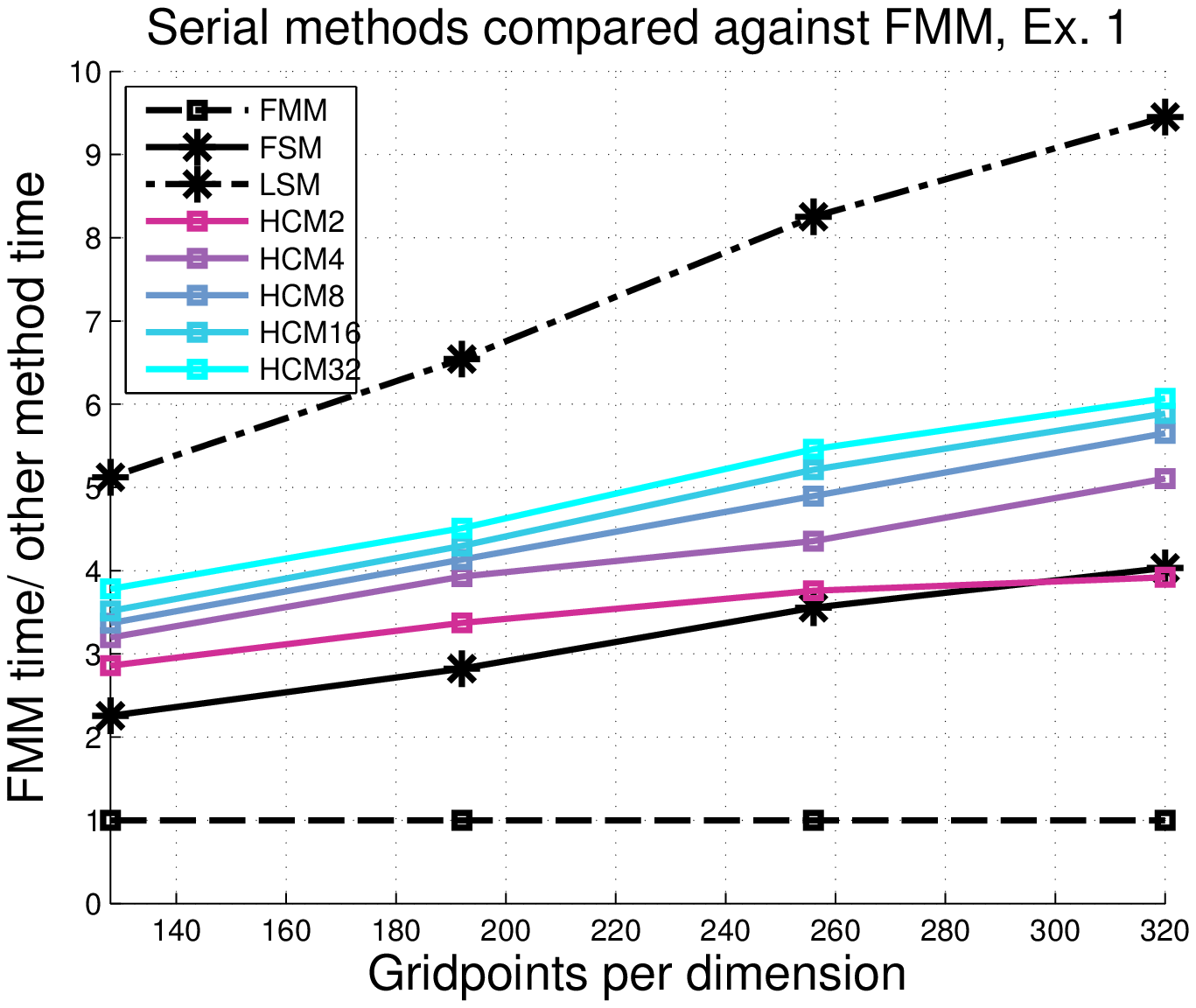}&
\hspace{-.2in}
\includegraphics[scale = .45]{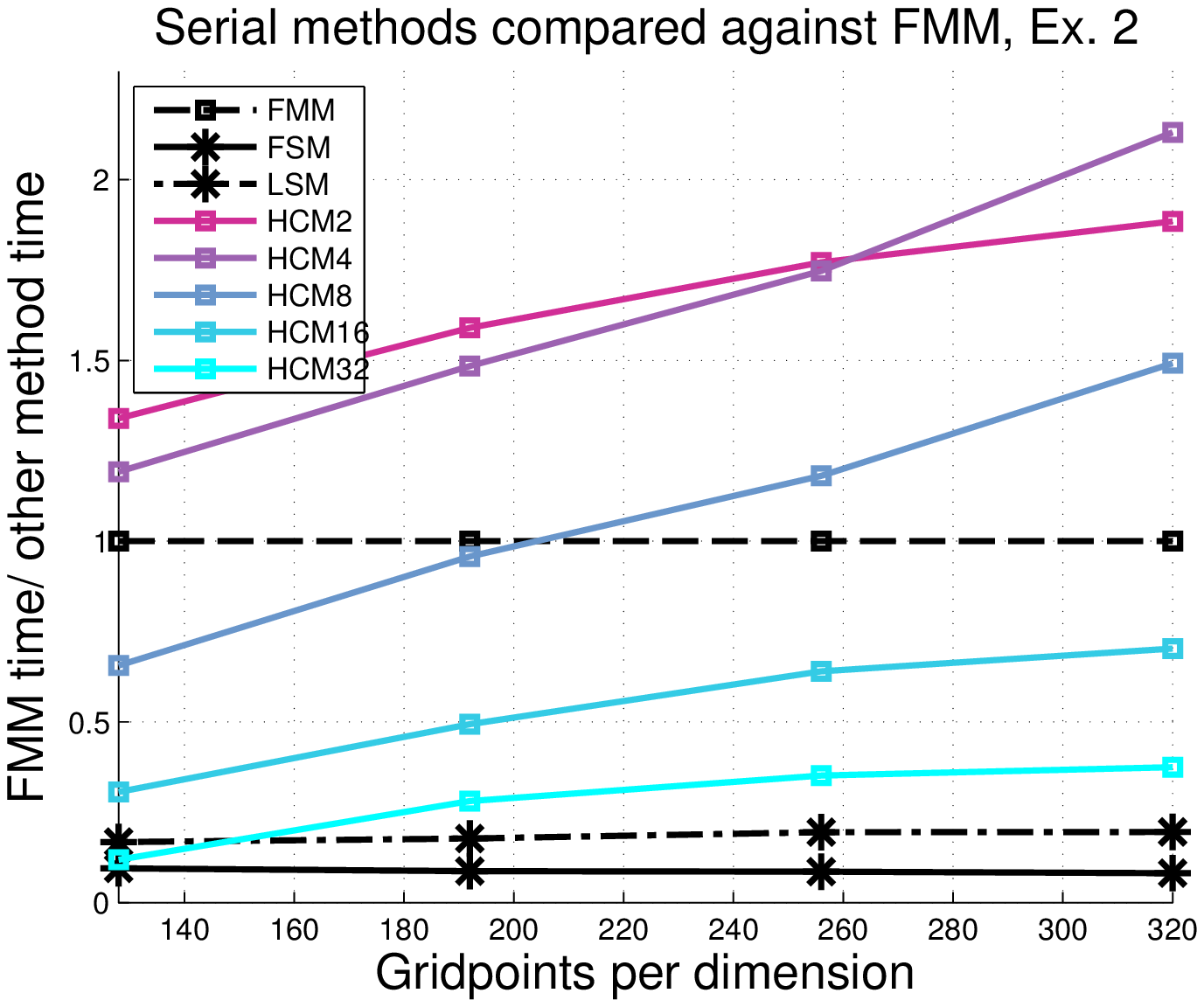}&
\hspace{-.2in}
\includegraphics[scale = .45]{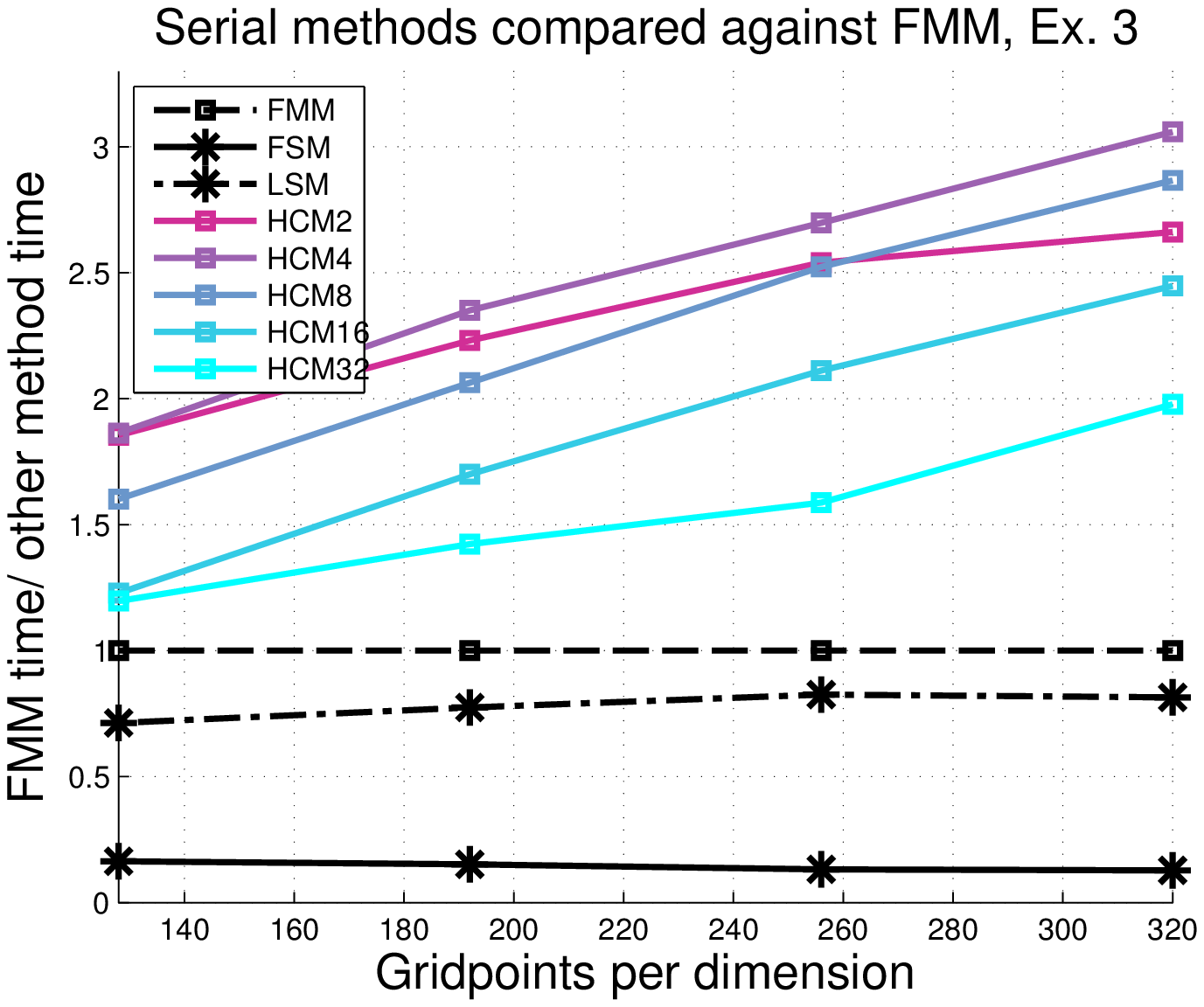}\\

A&B&C\\
\end{array}
$

\caption{\footnotesize Performance of the serial methods for different $M$.  The first chart has $F \equiv 1$, the second has $F = 1 + .5\sin{(20\pi x)}\sin{(20\pi y)}\sin{(20\pi z)}$, and the third has $F = 1 + .99\sin{(2\pi x)}\sin{(2\pi y)}\sin{(2\pi z)}$.  The data is given as a ratio of  FMM's CPU time to the times of all other method.}
\label{fig:serial_charts}
\end{figure}

\begin{figure}
\hspace{-.8 in}
$
\begin{array}{ccc}
\includegraphics[scale = .45]{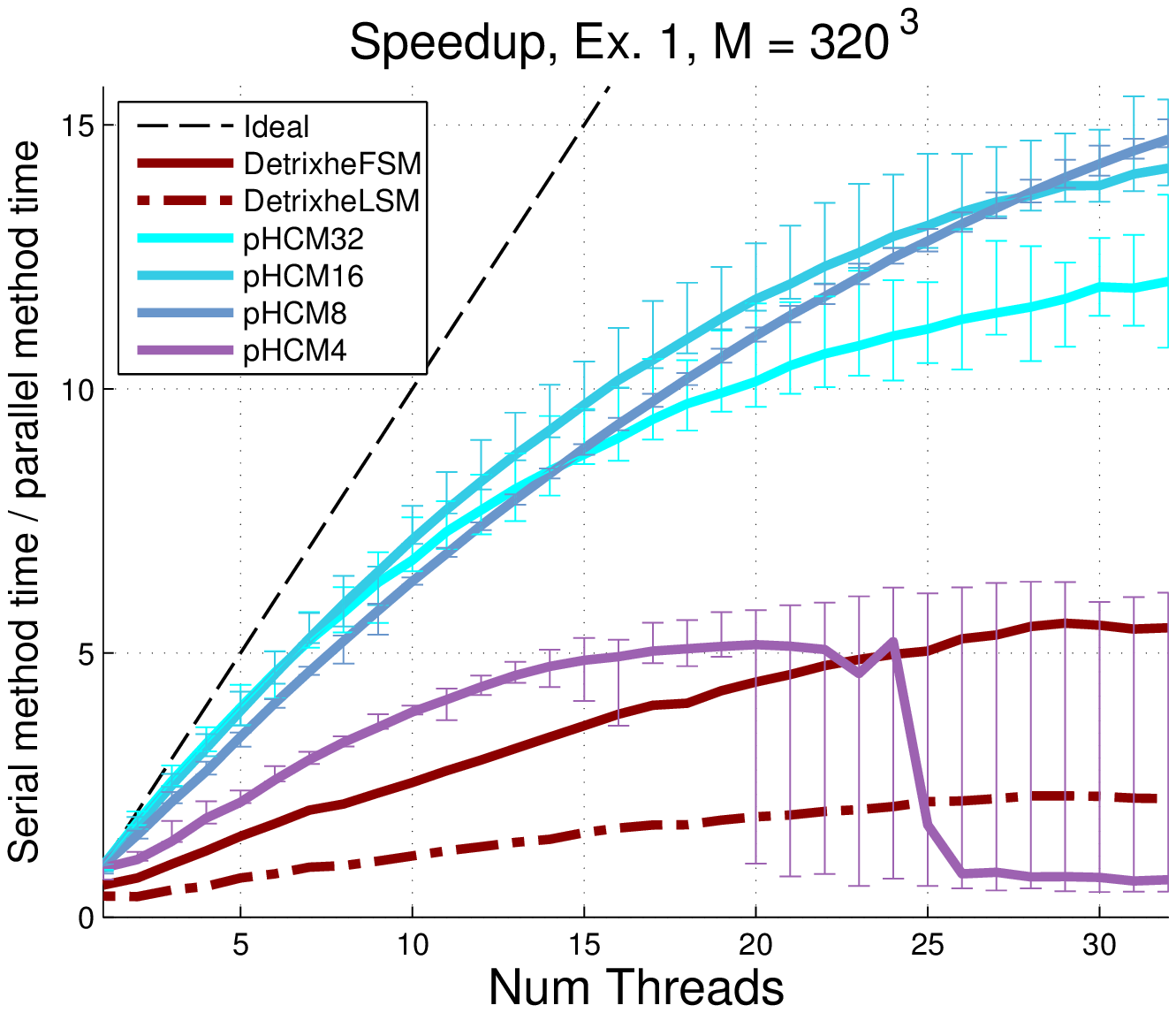}&
\hspace{-.2in}
\includegraphics[scale = .45]{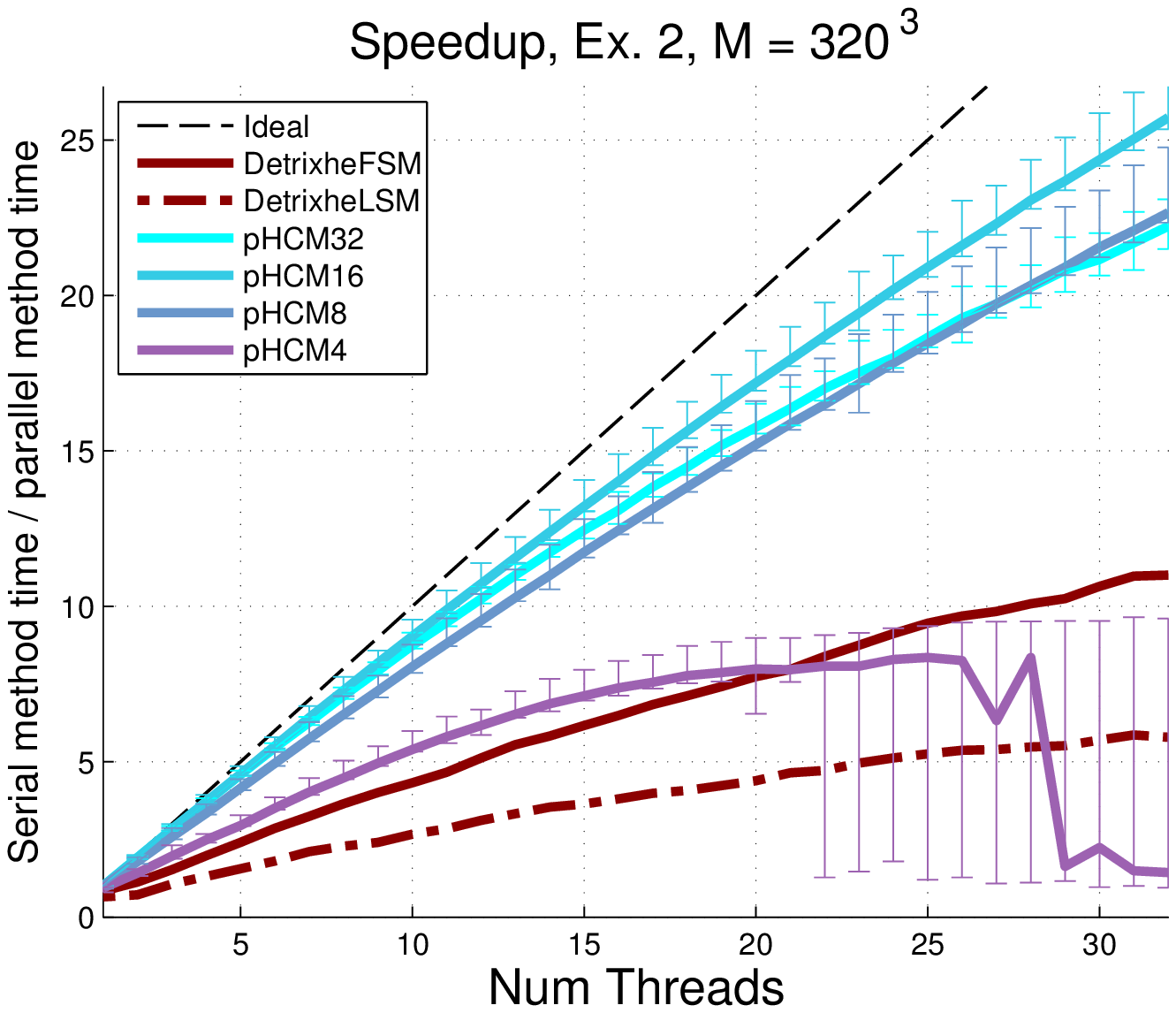}&
\hspace{-.2in}
\includegraphics[scale = .45]{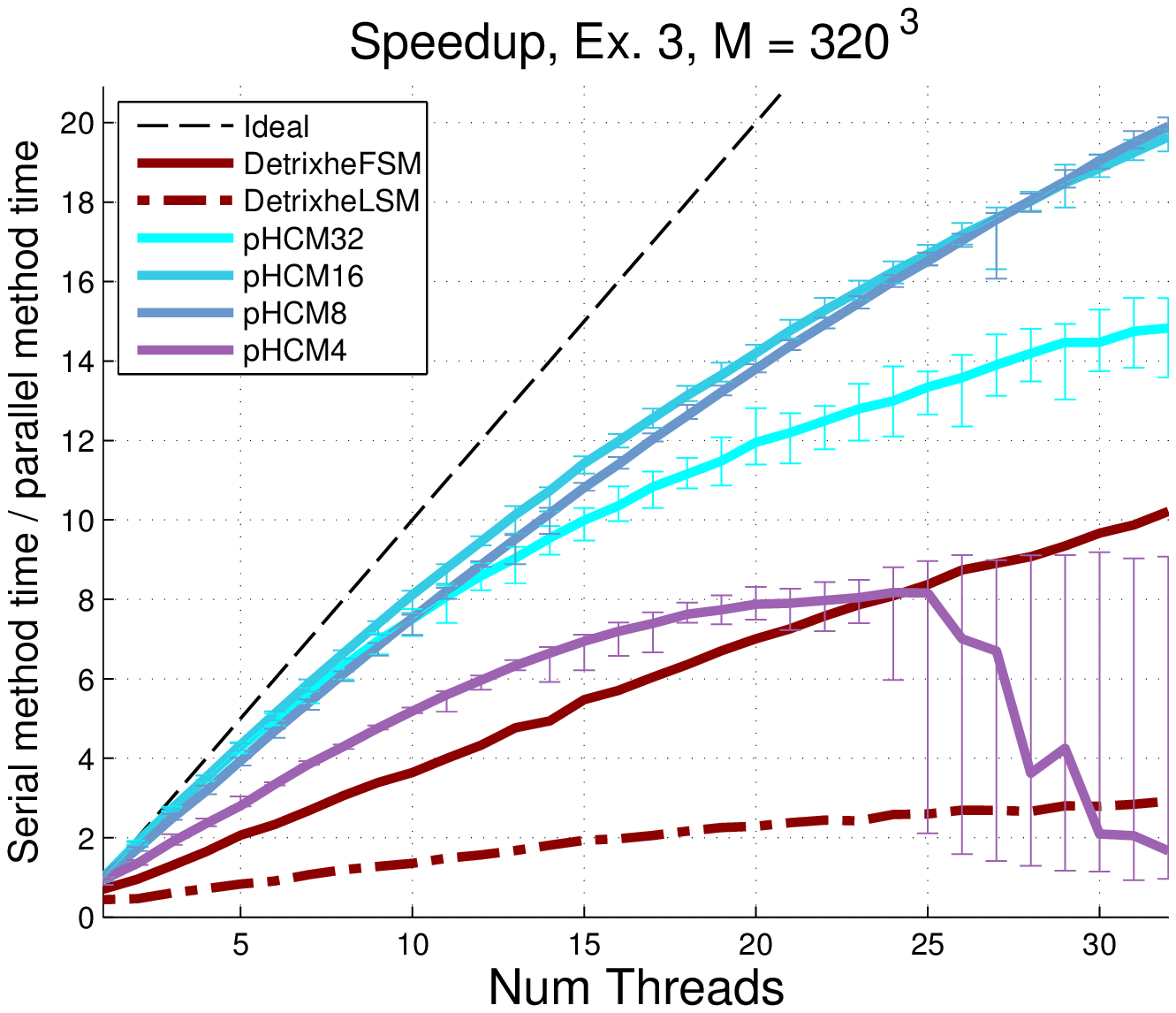}\\

A&B&C\\

\includegraphics[scale = .45]{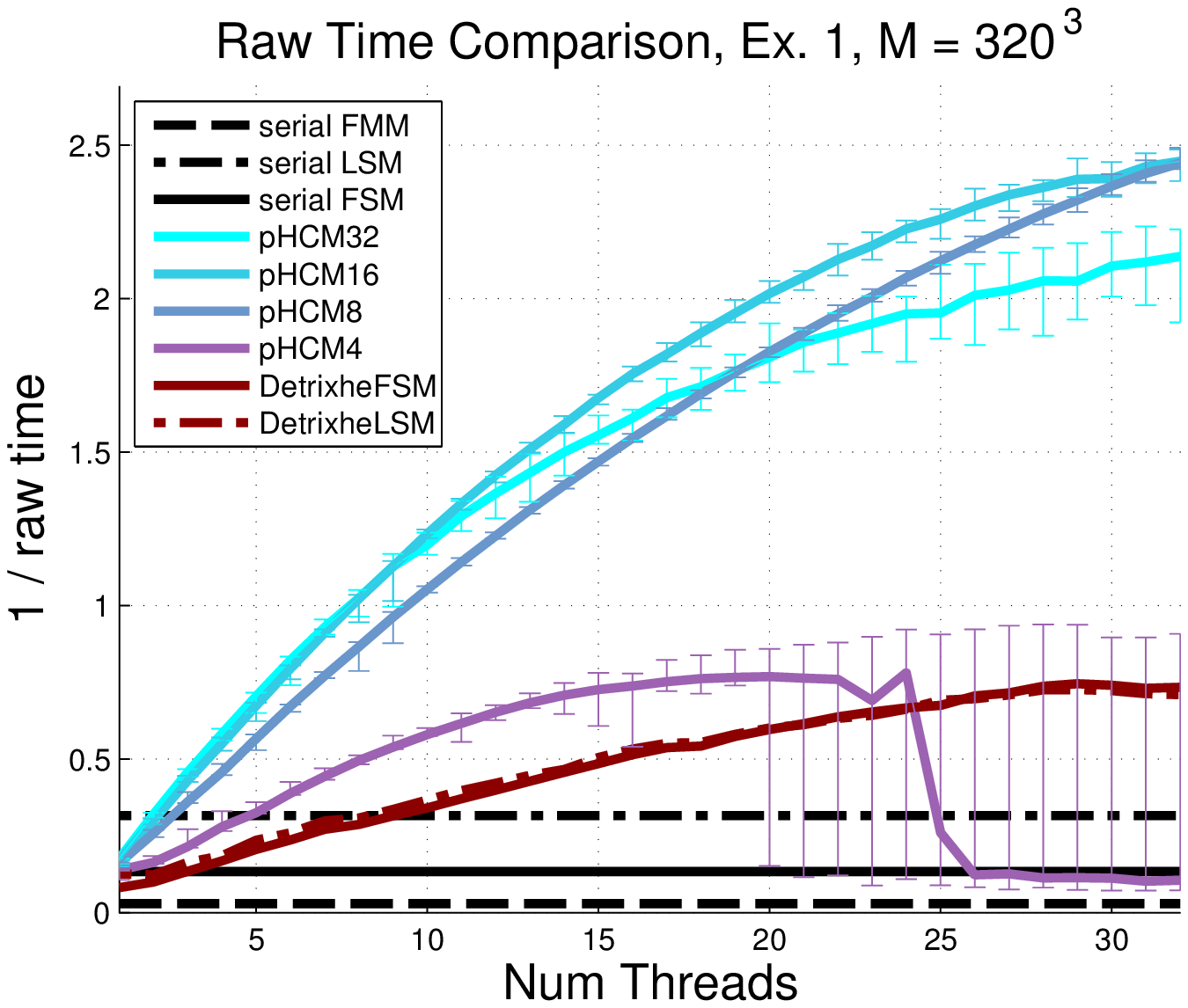}&
\hspace{-.2in}
\includegraphics[scale = .45]{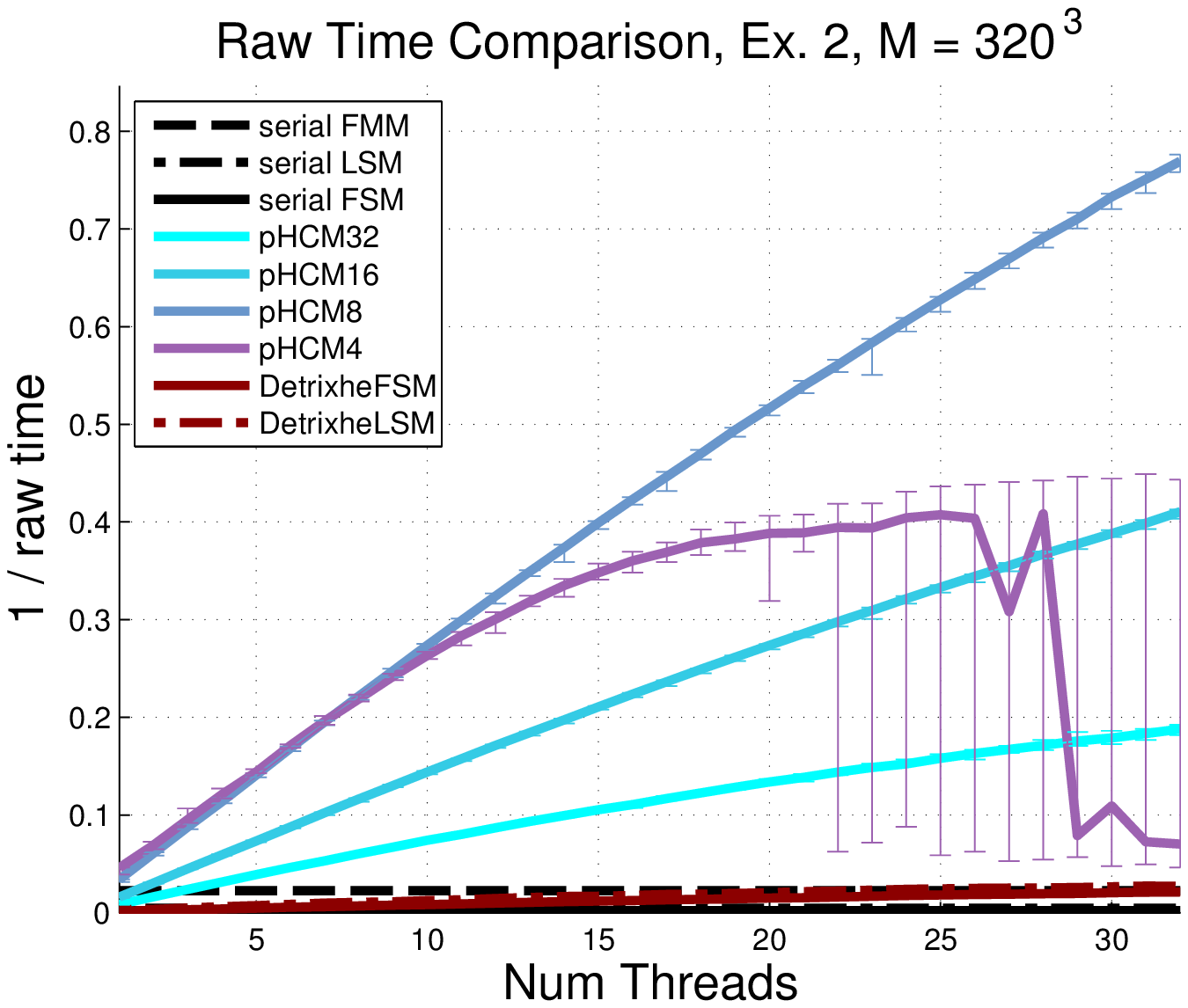}&
\hspace{-.2in}
\includegraphics[scale = .45]{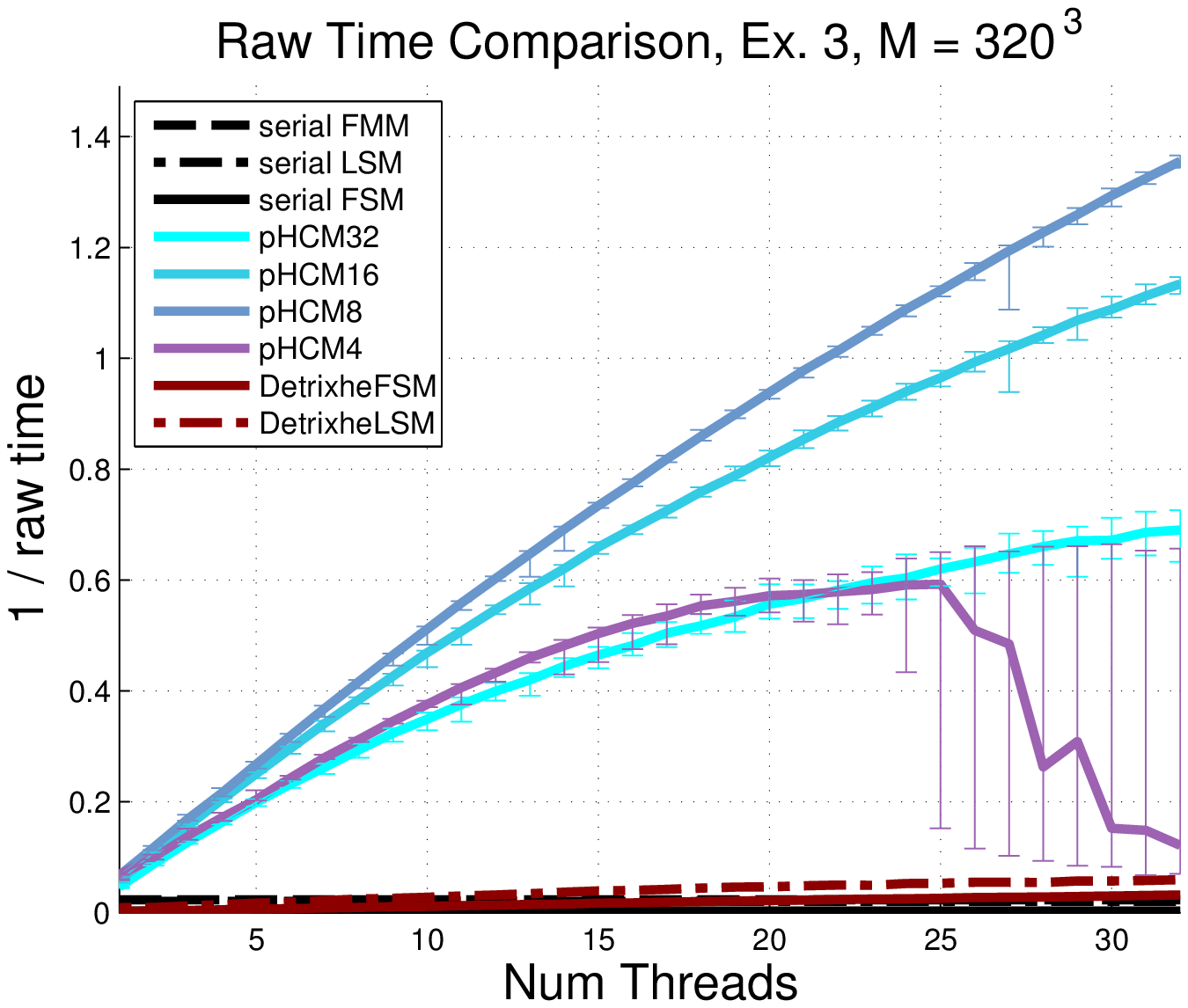}\\

D&E&F

\end{array}
$

\caption{\footnotesize Scaling and performance for pHCM at $M = 320^3$.  The first column has $F \equiv 1$, the second has $F = 1 + .5\sin{(20\pi x)}\sin{(20\pi y)}\sin{(20\pi z)}$, and the third has $F = 1 + .99\sin{(2\pi x)}\sin{(2\pi y)}\sin{(2\pi z)}$.}

\label{fig:parallel_charts}
\end{figure}

Figures \ref{fig:serial_charts}, \ref{fig:parallel_charts}, \ref{fig:overheads_320}, and \ref{fig:parallel_128} are organized so that columns present different examples and rows give different comparison metrics\footnote{Table versions of the same benchmarking results are also included in the ``supplementary materials''.}.  Figure \ref{fig:serial_charts} compares the performance of serial methods by  plotting the ratio of FMM CPU-time to other methods' times for $M = 128^3, 192^3, 256^3,$ and $320^3$.
Since we are interested in strong scalability, we test pHCM$r$ with a fixed problem size while varying $P$.
In Figure \ref{fig:parallel_charts}, \bm{M} {\bf is frozen at \bm{320^3}}.  The first row reports the speedup factors of the parallel methods over the serial methods; these are (HCM$r$ time / pHCM$r$ time), (FSM time / DFSM time), and (LSM time / DLSM time).  The second row of Fig. \ref{fig:parallel_charts} provides the performance comparison of all serial and parallel methods.
The growth of parallel overhead and the change in total work (as functions of $P$) are presented for each pHCM$r$ in Figure \ref{fig:overheads_320}.  Plots similar to Figure \ref{fig:parallel_charts} but computed for $M = 128^3$ are presented in subsection \ref{ss:parallel_performance}.

%\iffullversion
%\pagebreak
%\fi

\noindent

\noindent

\noindent {\bf Main observations:}

\begin{enumerate}

\item LSM significantly outperforms FMM on example 1 (Fig. \ref{fig:serial_charts}$A$) and its advantage grows with $M$. FMM greatly outperforms LSM on example 2 (Fig. \ref{fig:serial_charts}$B$) for all values of $M$. Their performance is more comparable on the third example (Fig. \ref{fig:serial_charts}$C$).

\item The performance ranking among serial HCM$r$ methods is problem-dependent (Fig. \ref{fig:serial_charts}$A$-\ref{fig:serial_charts}$C$).

\item Figures \ref{fig:parallel_charts}$D$-\ref{fig:parallel_charts}$F$ demonstrate that pHCM has a large advantage over all serial methods for most $r$ and $P$ combinations.  On the three examples with $M = 320^3$, the median performance for pHCM8 on 32 threads was between 34 and 84 times faster than FMM, between 7.7 and 166 times faster than LSM, and between 18.4 and 436 times faster than FSM.

\item Generally, the pHCM speedup over HCM is greater when there is more work per cell.  We see in Figures \ref{fig:parallel_charts}$A$-\ref{fig:parallel_charts}$C$ that the experiments with higher gridpoints-per-cell number $r$ exhibit better parallelization, and the speedup of pHCM4 is always the worst.

\item   In Figure \ref{fig:parallel_charts} the position of each curve relative to its error bar reveals the most likely outcome.  For example,  the pHCM4 scaling plummets in the worst cases and plateaus in the best cases.  At 32 threads, since the median is near the bottom of the error bar in all examples, the good cases are relatively rare.

\item Based on Figure \ref{fig:parallel_charts}, for most $r$ values pHCM scales much better than DFSM/DLSM.
 %COMMENT REMOVED
     Since DFSM is a synchronous parallel algorithm, it comes as no surprise that using the Locking Sweeping does not boost performance significantly -- LSM only reduces the amount of work performed by a subset of the threads.
    %COMMENT REMOVED
    %COMMENT REMOVED
    %COMMENT REMOVED
    Better scaling in DLSM would likely be achieved if
    %COMMENT REMOVED
    it were possible
    to apply a special load balancing procedure
    %COMMENT REMOVED
    based on the set of currently ``active" gridpoints.

\end{enumerate}

%COMMENT REMOVED

\subsection{Further comments on performance of serial methods}
\label{ss:serial_performance}

\begin{enumerate}

\item {\it Tradeoffs between FMM and LSM.} It is well known that Marching and Sweeping methods are each advantageous on their own subsets of Eikonal problems.  The exact delineation remains a matter of debate.  The readers can find careful comparative studies in \cite{GremaudKuster,HysingTurek} and partly in \cite{ChacVlad}.
In each example (Figs. \ref{fig:serial_charts}$A$-\ref{fig:serial_charts}$C$) we observe that, as $M$ increases, the ratio of FMM time to LSM time increases due to the greater cost of each heap-sort operation.  However, FMM's performance is much more robust to the qualitative differences in the solution; FMM's raw times for $M = 320^3$ ranged between 32s (Ex. 1) and 51s (Ex. 2), while the LSM times were between 3s (Ex. 1) and 363s (Ex. 2).  FMM is also usually much more efficient on problems with complicated domain geometry (e.g., on domains containing multiple impenetrable obstacles).

\item {\it Grid memory layout and caching issues.} Large grids, particularly common in higher dimensional problems, present an additional challenge for all (serial and parallel) methods implemented on a shared memory architecture.  Solving equation \eqref{eq:Eik_discr} requires accessing the $U$ values for all gridpoints neighboring $\x_{ijk}$, but the geometric neighbors can be far apart in memory when the higher-dimensional grid is stored lexicographically.  This results in frequent cache-swapping, ultimately impacting the computational cost.  More detailed profiling (not included here) confirms the resulting slow-down in all serial methods, including LSM.  In other applications space-filling curves have been successfully used to alleviate this problem (e.g., \cite{spaceCurve}), but we are not aware of any successful use in Eikonal solvers.
    We believe that allocating the fine grid separately per-cell would be advantageous for a robust extension of HCM/pHCM to higher dimensions.  However, our current implementation of heap-cell methods does not take advantage of this idea.
    %COMMENT REMOVED

\item{\it FMM scaling in $M$.}  Since the length of the heap increases with $M$, the number of flops per heap operation increases too.  On top of this, FMM is affected by additional caching issues: the time per heap-related memory access increases, since the parent/child relationships of heap entries do not translate to memory proximity of the corresponding gridpoints.  Profiling shows that the cache miss rate increases noticeably with $M$.

\item{\it HCM scaling in $M$.} For most cell decompositions, when $J \ll M$, the heap maintenance is negligible.
    As $J$ becomes large (e.g., for $r=2$), HCM$r$ is affected by the same issues described for FMM above.

\item {\it Optimal $J$ in HCM.}   As cell sizes decrease, the causality among cells becomes stronger (see the end of Section \ref{s:HCM}) and our cell value heuristic does a better job of capturing the dependency structure; the average number of times each cell is processed tends to 1.  Additionally, the characteristics within each cell become approximately straight lines, so the per-cell LSM converges quickly.  On the other hand, if $J$ is large enough, the overhead due to heap maintenance becomes significant; this is quantified in Tables \ref{tab:const_HCM_analysis}, \ref{tab:sin1_HCM_analysis}, and \ref{tab:sin2_HCM_analysis} (``Heap Maintenance \%" means the percentage of execution time spent outside of sweeping cells).
Turning to individual examples:

\begin{enumerate}
\item Ex.1:  HCM with larger cell sizes performs better.  See Figure \ref{fig:serial_charts}$A$ and Table \ref{tab:const_HCM_analysis}.  This %COMMENT REMOVED
    is due to a very special property of $F \equiv 1$: %COMMENT REMOVED
    since there is exactly one heap removal per cell regardless of $J$, the maintenance of the heap is the dominant factor affecting the performance.  Correspondingly, LSM performs the best.  (LSM is equivalent to HCM using only one cell.)
\begin{table}[h] \footnotesize
\begin{center}
\caption{Performance analysis of HCM on Ex. 1, $M =320^3$.}\vspace*{2mm}
\begin{tabular}{|c|c|c|c|c|c|}
\hline
& \textbf{HCM32} & \textbf{HCM16} & \textbf{HCM8} & \textbf{HCM4} & \textbf{HCM2}\\
\hline
\textbf{Avg. Sweeps per Cell}  &4.84 & 4.92 & 4.96 & 4.98 & 4.12\\
\hline
\textbf{Heap Maintenance \%}  & 1.09 & 1.12 & 1.66 & 5.88 & 33.9\\
\hline
\end{tabular}
\vspace*{.5cm}

\label{tab:const_HCM_analysis}
\end{center}
\end{table}

\item Ex. 2: Due to the %COMMENT REMOVED
        oscillatory nature of characteristics, HCM performs better with \emph{smaller} cell sizes.  The ranking among HCM$r$ methods is more or less the reverse of that for example 1, and the sweeping methods are the slowest.  See Figure \ref{fig:serial_charts}$B$ and Table \ref{tab:sin1_HCM_analysis}.
\begin{table}[h] \footnotesize
\begin{center}
\caption{Performance analysis of HCM on Ex. 2, $M =320^3$.}\vspace*{2mm}
\begin{tabular}{|c|c|c|c|c|c|}
\hline
& \textbf{HCM32} & \textbf{HCM16} & \textbf{HCM8} & \textbf{HCM4} & \textbf{HCM2}\\
\hline
\textbf{Avg. Sweeps per Cell}  &223 & 100 & 31.1 & 12.9 & 6.97\\
\hline
\textbf{Heap Maintenance \%}  & 0.076 & 0.214 & 0.954 & 4.95 & 30.6\\
\hline
\end{tabular}
\vspace*{.5cm}

\label{tab:sin1_HCM_analysis}
\end{center}
\end{table}

\item Ex. 3: Figure \ref{fig:serial_charts}$C$ and Table \ref{tab:sin2_HCM_analysis} show that the performance among the HCM$r$ methods is %COMMENT REMOVED
    qualitatively different from the previous examples.  A weakly causal ordering already exists here for moderately-sized cells.
\begin{table}[h] \footnotesize
\begin{center}
\caption{Performance analysis of HCM on Ex. 3, $M =320^3$.}\vspace*{2mm}
\begin{tabular}{|c|c|c|c|c|c|}
\hline
& \textbf{HCM32} & \textbf{HCM16} & \textbf{HCM8} & \textbf{HCM4} & \textbf{HCM2}\\
\hline
\textbf{Avg. Sweeps per Cell}  &29.3 & 14.6 & 9.37 & 7.14 & 5.02\\
\hline
\textbf{Heap Maintenance \%}  & 0.292 & 0.424 & 0.914 & 4.55 & 28.5\\
\hline
\end{tabular}
\vspace*{.5cm}

\label{tab:sin2_HCM_analysis}
\end{center}
\end{table}

\end{enumerate}

\end{enumerate}

\iffullversion
\else
%COMMENT REMOVED
%COMMENT REMOVED
The sweeping methods can be accelerated by stopping at the iteration where the maximum change over gridpoint values is less than or equal to a certain threshold $\kappa \geq 0$.  If $\kappa > 0$, the method will terminate ``early'', and the output will be different than the true solution of the discretized system (\ref{eq:Eik_discr}).
Ideally, $\kappa$ should be chosen based on the $L_{\infty}$-norm discretization error, but since the latter is a priori unknown, a common practical approach is to use a small heuristically selected constant (e.g., \cite{Zhao}).  We note that, for a fixed $\kappa>0$, the number of needed iterations can be quite different for different $h$,
and there is currently no proof that the early-terminated numerical values are within $\kappa$ from the correct solution.
See also the discussion and experimental results in \cite{ChacVlad_expanded}.
The results reported here were obtained with $\kappa=0$.
Testing with $\kappa=10^{-8}$, reduces the number of sweeps but the scaling behavior remains largely the same; see \cite{pHCM_expanded}.
%COMMENT REMOVED
%COMMENT REMOVED
\fi

\subsection{Detailed performance analysis of parallel methods}
\label{ss:parallel_performance}

Two key factors that affect the speedup of parallel methods are the amount of parallel overhead (contention, inter-thread communication, etc.) and the change in the amount of work performed from serial to parallel.  In this section we focus on both the overhead analysis and the algorithmic differences between pHCM and HCM.  The overhead is the sum of the parallel overhead and the ``base" heap maintenance.  The latter is given above in Tables \ref{tab:const_HCM_analysis}, \ref{tab:sin1_HCM_analysis}, and \ref{tab:sin2_HCM_analysis}.

We define:

$\bullet$ AvS = $\Sigma_{p=0} ^{P-1}$ (Total number of sweeps performed by processor $p$) $/ J$.

$\bullet$  Cell Comp \% = percent of total time spent on sweeping cells alone.

$\bullet$ Overhead \% = 100\% - Cell Comp \%, i.e., percent of total time spent beyond sweeping cells.

\noindent

\begin{figure}

\hspace{-.4 in}
%COMMENT REMOVED
$
\begin{array}{ccc}

\includegraphics[scale = .4]{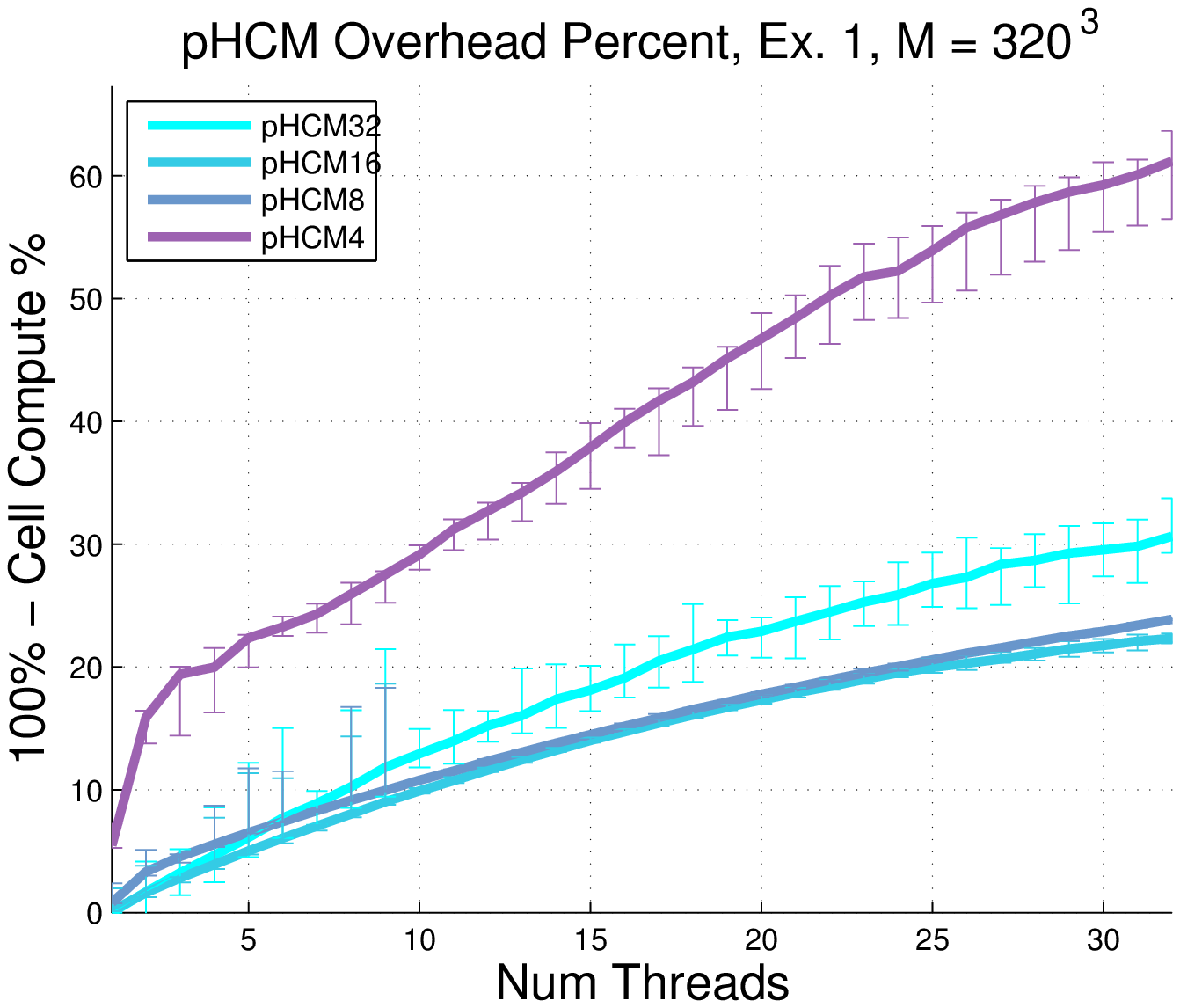}&
\hspace{-.2in}
\includegraphics[scale = .4]{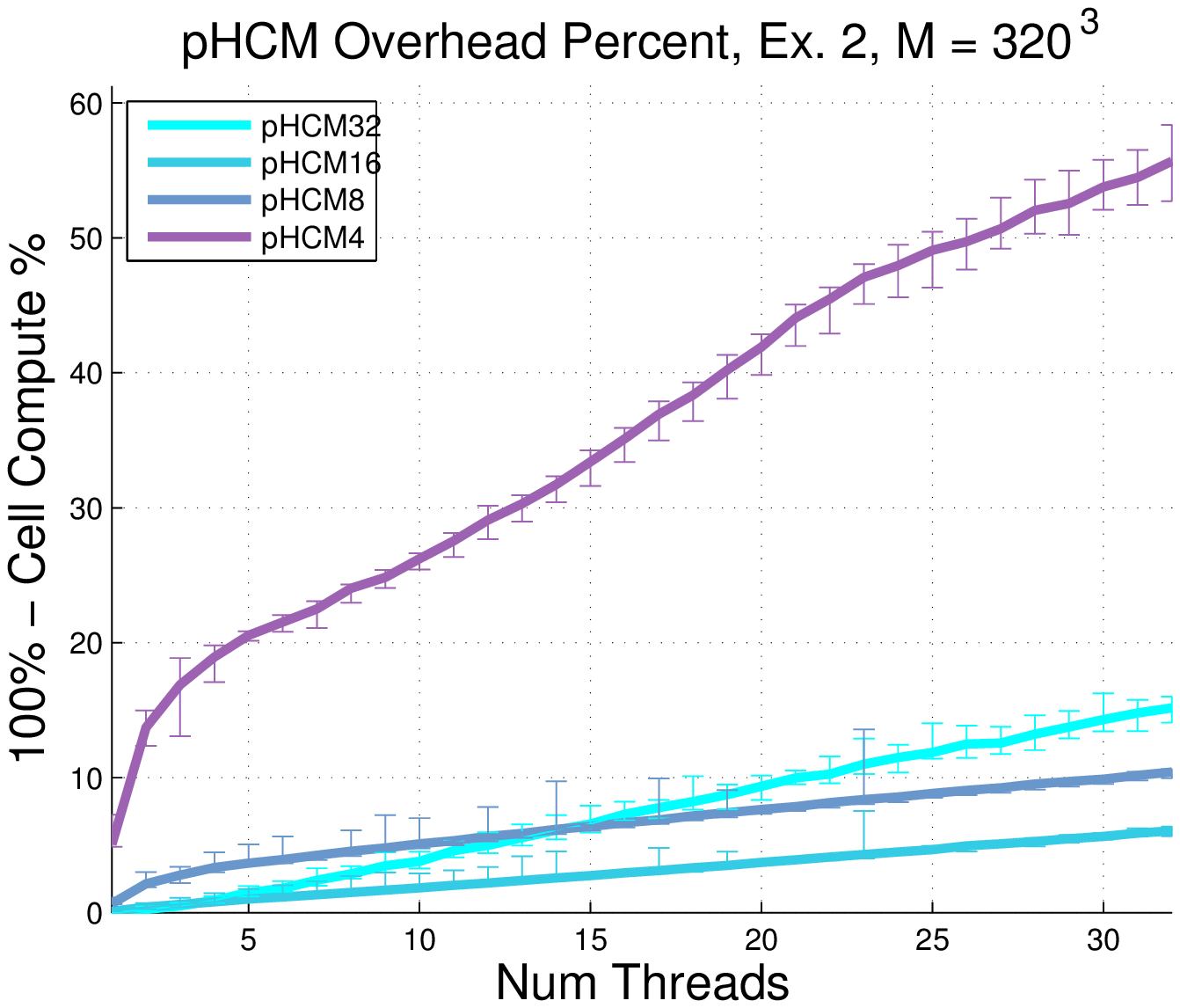}&
\hspace{-.2in}
\includegraphics[scale = .4]{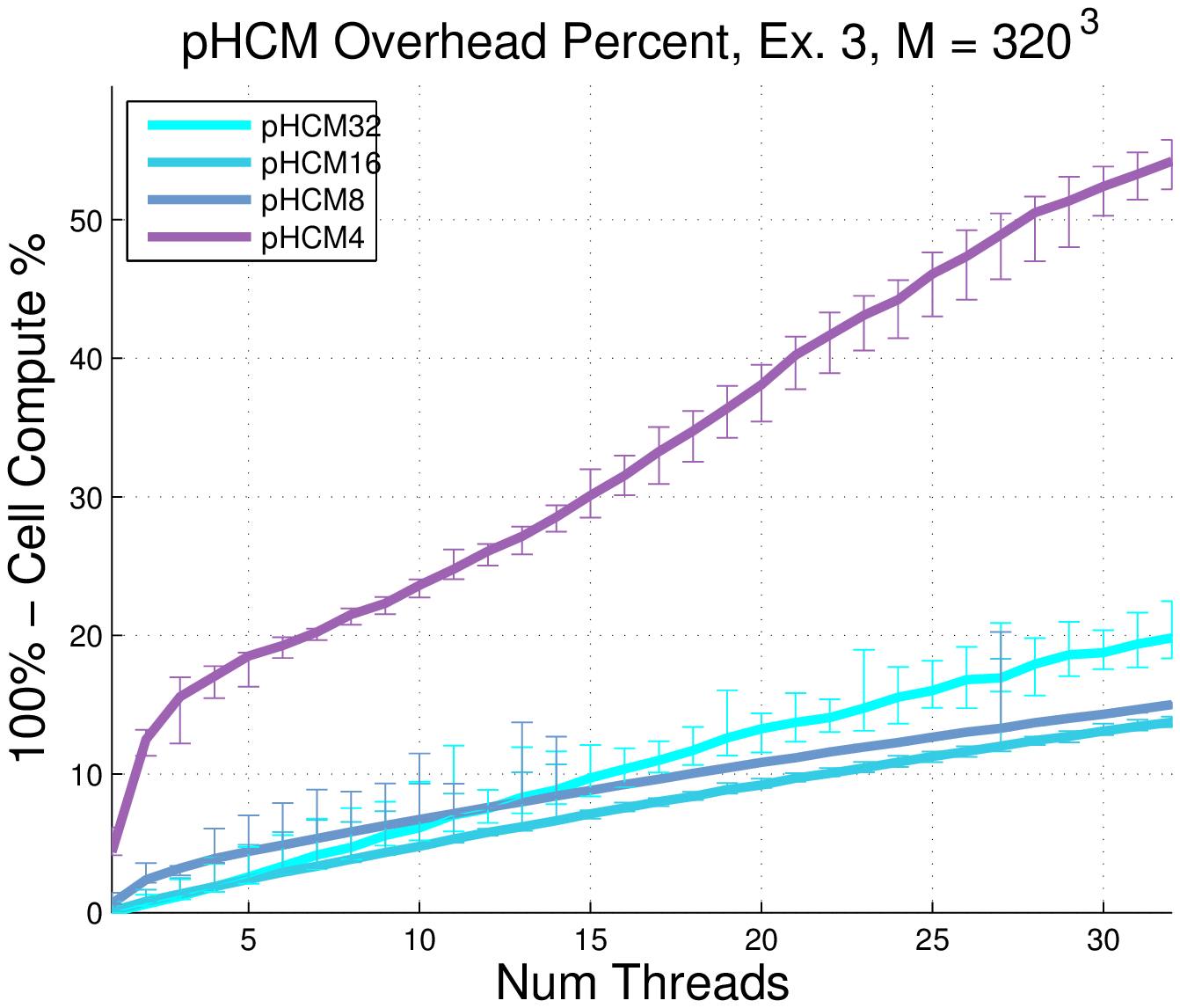}\\

A&B&C\\

\includegraphics[scale = .4]{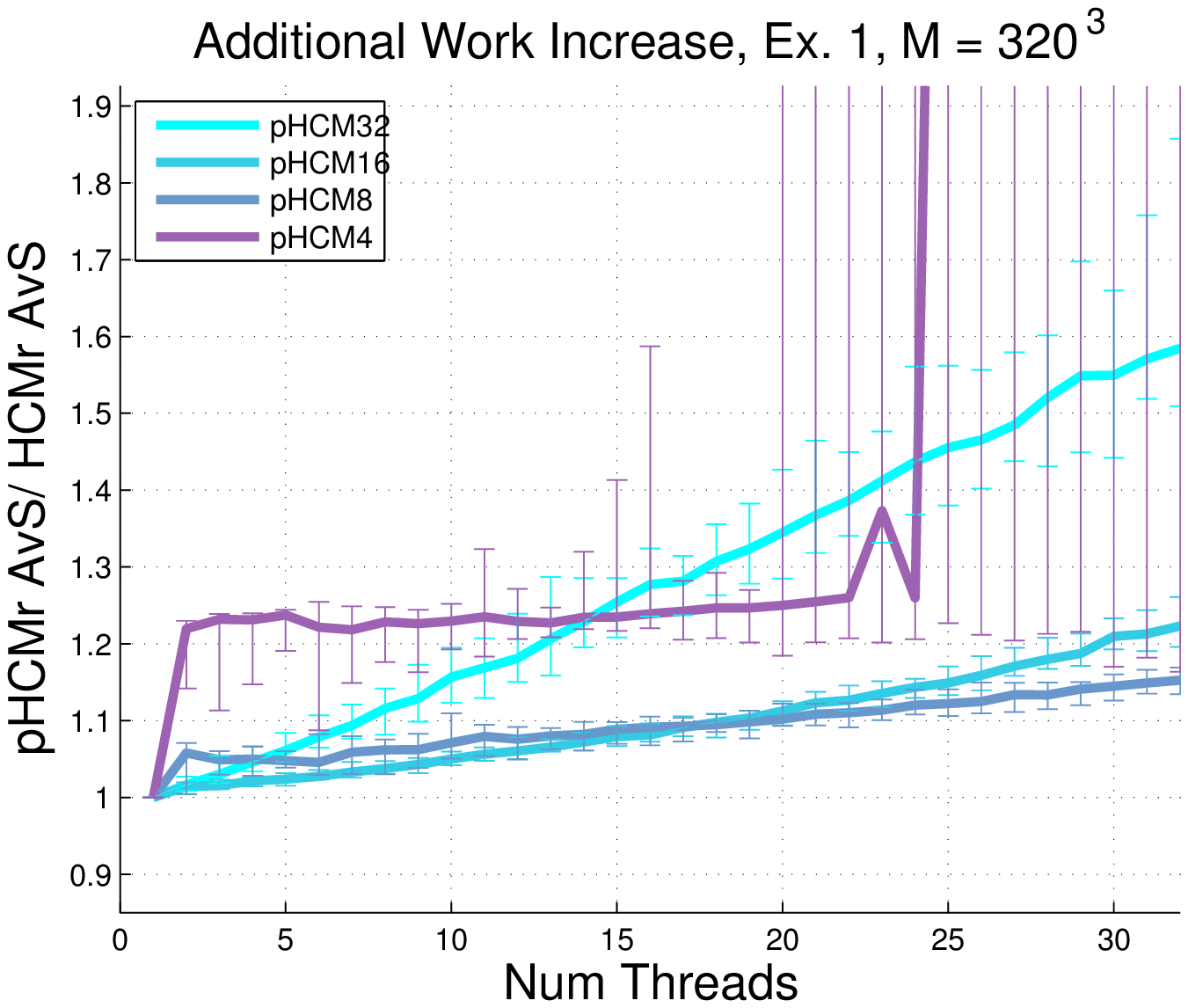}&
\hspace{-.2in}
\includegraphics[scale = .4]{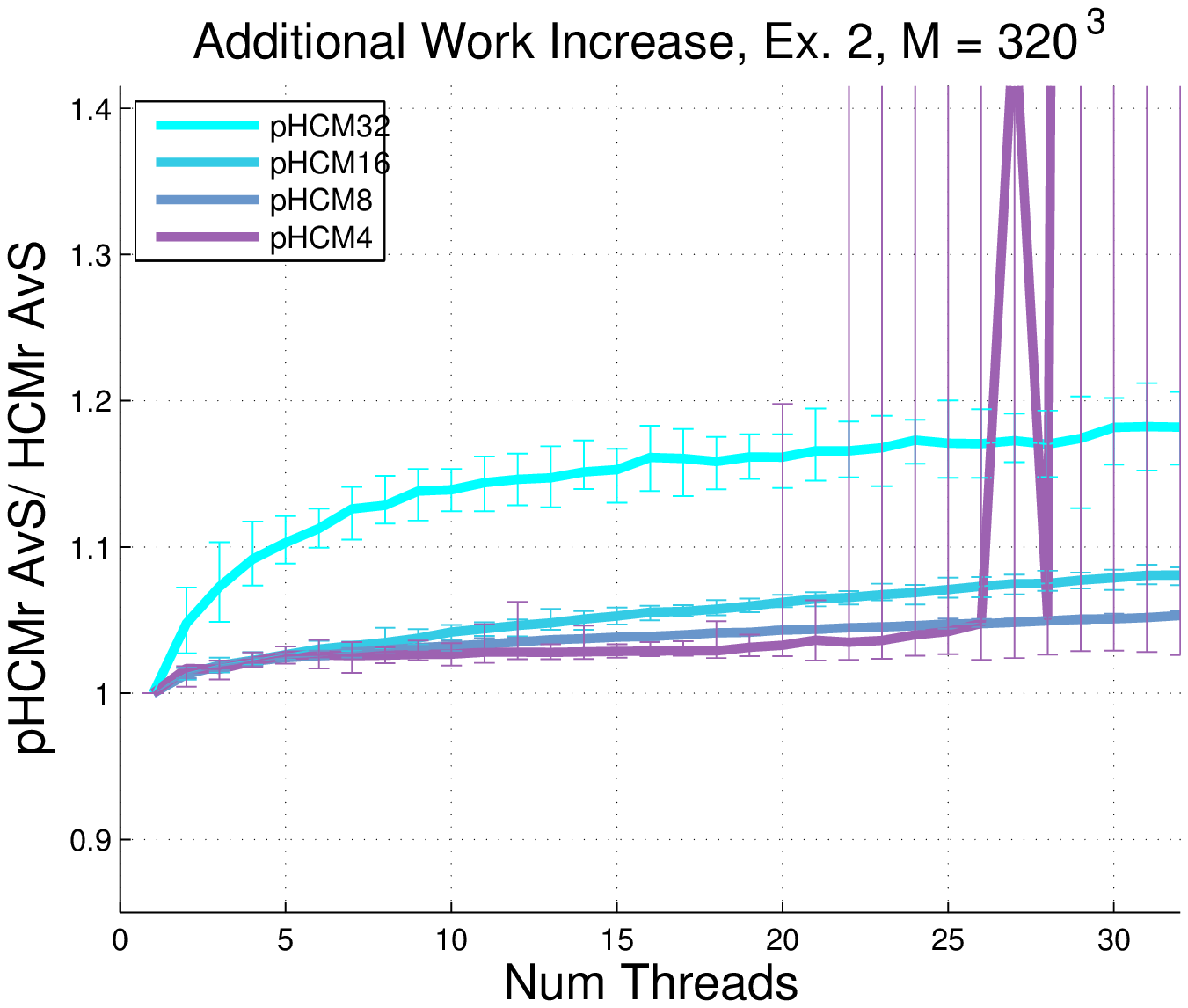}&
\hspace{-.2in}
\includegraphics[scale = .4]{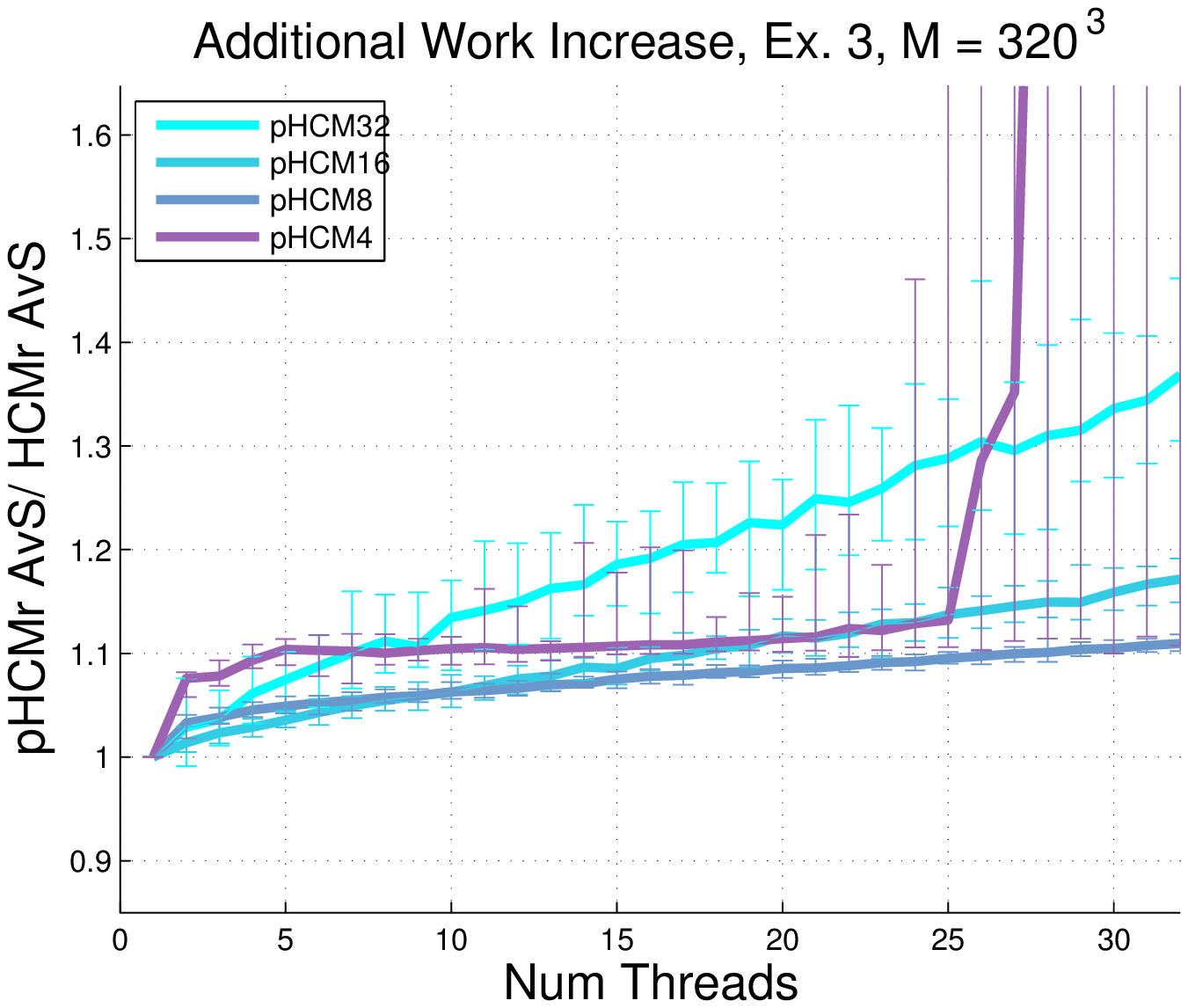}\\

D&E&F\\
\end{array}
$

\caption{\footnotesize Overhead percentages and additional work in pHCM$r$ for different $P$ for the three examples, with $M = 320^3$.  In figures A, B, and C the value at Num Threads $= 1$ of each curve approximately gives the part of the overhead accounted for by heap maintenance alone; the \emph{parallel} overhead would be given approximately by subtracting it from each curve.}

\label{fig:overheads_320}
\end{figure}

\begin{enumerate}

%COMMENT REMOVED
\item{\it Effects of $P$ on overhead.} As $P$ increases, contention and network communication increase.  If more threads are used for a given cell discretization, it is more likely for a processor $\hat{p}$ to wait to obtain a lock (e.g., as in line 8 of algorithm \ref{alg:pHCM}).

\item{\it Effects of $J$ on overhead.} The overhead percentage can be large if either 1) $J$ is large, so processors spend more time doing heap sorts and contending with each other to obtain locks to shared data structures, or 2) $J$ is small and $P$ is large, so there is not enough total work to be divided among the processors.  In this case a processor may spend a significant amount of time outside the main loop just waiting for work.

\item{\it Effect of a strong causal structure.}  The order of processing the cells is different for pHCM and HCM.  On Ex. 1 (Figure \ref{fig:serial_charts}$A$) there is a strict causal relationship among cells, resulting in exactly 1 heap removal per cell in HCM.  For pHCM the AvS is larger since cells are not generally processed in their strict causal order.   In fact, on \emph{any} problem for which HCM has exactly one heap removal per cell, pHCM will almost surely see an increase in the total number of heap removals.  A cell is added to some processor's heap when one of its neighbors updates a gridpoint on the inflow boundary. But with a strictly causal cell decomposition, this may well result in avoidable/premature computations if that cell is actually processed before {\em all} of its inflow boundary data is finalized.    This situation is particularly common when $P$ is large and $J$ is relatively small.

\item  {\it Effects of multiple caches.} Even by comparing only the time spent on cell-level sweeping (and accounting for differences in the total AvS) one sees that the speedup factor is closer to $P$ but not exact.
%COMMENT REMOVED
When $P$ is larger it is more likely that adjacent cells will be processed simultaneously, a situation whereby individual sweeps may become slower than their serial counterparts.  Referring back to Figure \ref{fig:cell_inflow_boundary},
suppose in the process of updating a border gridpoint $\x_i \in A$ the value of its neighbor $\x_j \in B$ is loaded into the cache of the local processor $p_A$.  If $\x_j$ changes value as a result of sweeps on cell $B$, the value stored in $p_A$ will either need to be invalidated or have the new value communicated to it \cite{ompBook}.  This operation is orders of magnitude slower than simply updating a cached value without communication.

\item {\it Robustness of pHCM performance.}  There is a possibility of the total amount of work increasing significantly if processor speeds vary. %COMMENT REMOVED
    %COMMENT REMOVED
    Suppose processor $\hat{p}$ is slow or has become slow and is processing a high-priority cell $A$.  The other fast processors will not be able to do useful work on cells downwind from $A$.  What is more, there is a cascade effect: cells downwind from the downwind neighbors of $A$ will need to be readded, etc.  This effect is more commonplace for small cells, as seen in Fig. \ref{fig:overheads_320}$D$ - \ref{fig:overheads_320}$F$.  The non-robust performance of pHCM4 appears to be due entirely to this effect - the error bars for the work are large while those for the overhead are small.  Not surprisingly, pHCM2 (omitted in this paper) shows even less robustness than the reported pHCM$r$.  For small cells and large $P$, a synchronous parallel implementation may be a wiser choice.
    %COMMENT REMOVED
%COMMENT REMOVED

\item {\it Coarser grids.}   The charts in Figure \ref{fig:parallel_128} present the same information as in Figure \ref{fig:parallel_charts}, but for $M = 128^3$.  The speedup of the parallel methods here is expectedly worse than for $M = 320^3$.  %COMMENT REMOVED
    Indeed, for a fixed $r$ and $P$, a smaller $M$ yields a smaller number of cells $J$.
    For larger values of $P$, smaller $J$ results both in an increased overhead and premature processing of cells;
    see items 2 and 3 above.
    A good illustration of this is the pHCM32 curve in Fig. \ref{fig:parallel_128}$A$ and \ref{fig:parallel_128}$C$.  Since here $M = 128^3$, the cell decomposition for pHCM32 is only 4 cells per domain side; the scaling plateaus at a low number of threads.

\iffullversion
\item{\it Possible decrease in work.}  The total amount of work performed by pHCM may also actually \emph{decrease} compared to HCM in cases where the cell heuristic poorly predicts the dependency structure of the cells.  See subsection \ref{ss:otherCellVals}.
%COMMENT REMOVED
\fi

\item{\it Parallel Sweeping.}  As reported in \cite{Detrixhe}, the algorithmic complexity of Detrixhe Sweeping is constant in the number of threads; for DFSM and DLSM, charts like \ref{fig:overheads_320}$D$-\ref{fig:overheads_320}$F$  would all show a constant value of 1.  Unfortunately, the performance is also affected by the fact that memory access patterns are more complicated for DFSM/DLSM than for FSM/LSM, which %COMMENT REMOVED
    may prevent the compiler from taking advantage of data locality.
    %COMMENT REMOVED
%COMMENT REMOVED
%COMMENT REMOVED
Based on our own OpenMP implementation %COMMENT REMOVED
on a shared memory architecture, the scalability is also sensitive to hardware properties of the specific
%COMMENT REMOVED
\iffullversion
%COMMENT REMOVED
platform; see also subsection \ref{ss:octopus}.
%COMMENT REMOVED
\else
platform.
%COMMENT REMOVED
In \cite{pHCM_expanded} we have also observed somewhat better scaling of DFSM/DLSM on a different hardware system
with a smaller ratio of memory bandwidth to CPU speed.
\fi
We note that the authors of \cite{Detrixhe} have also implemented their method in lower-level memory languages (MPI, CUDA) to alleviate this sensitivity.\\
%COMMENT REMOVED
\end{enumerate}

\iffullversion
\begin{figure}[h]
\else
\begin{figure}
\fi
\hspace{-.8 in}
%COMMENT REMOVED
$
\begin{array}{ccc}

\includegraphics[scale = .45]{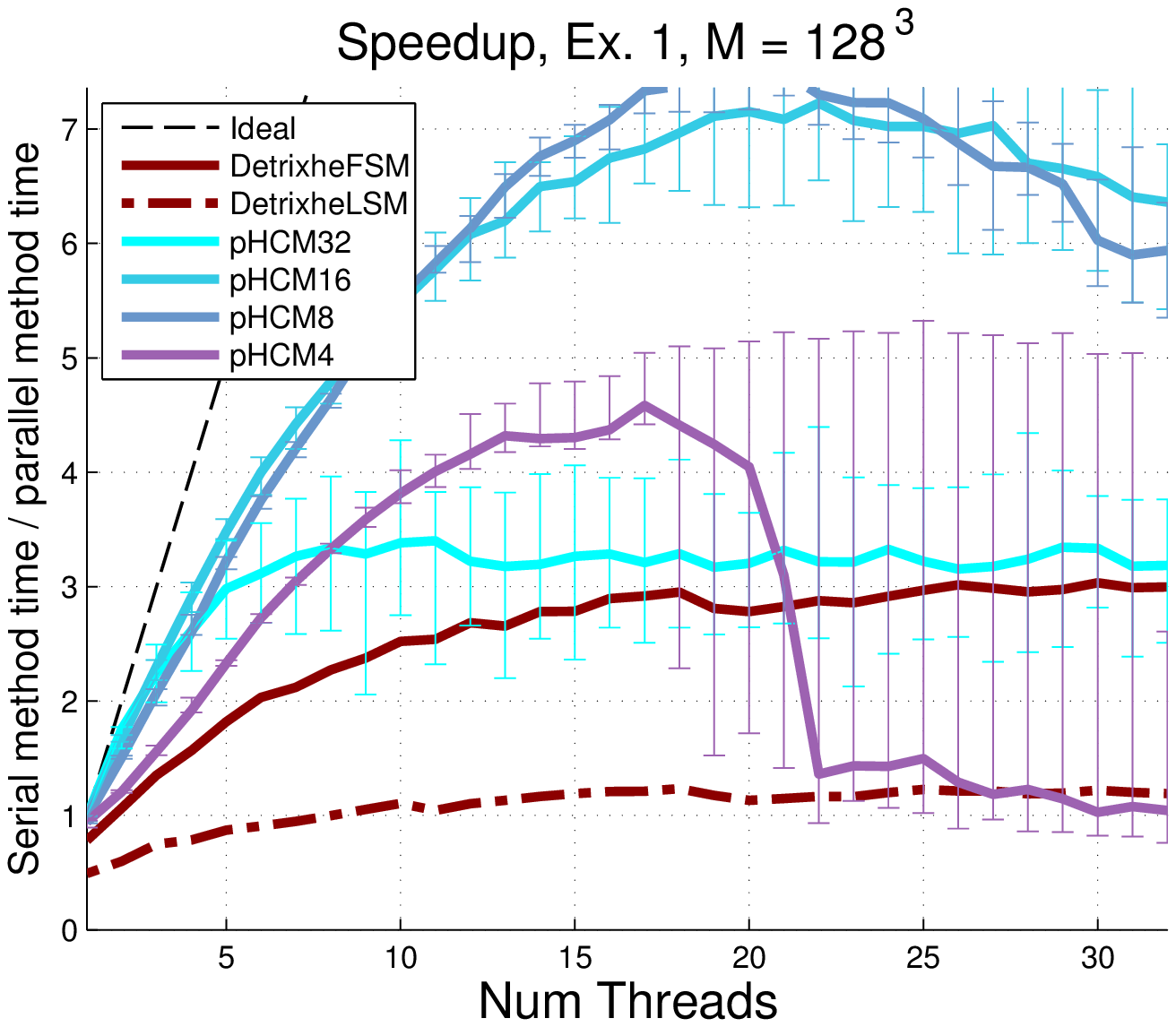}&
\hspace{-.2in}
\includegraphics[scale = .45]{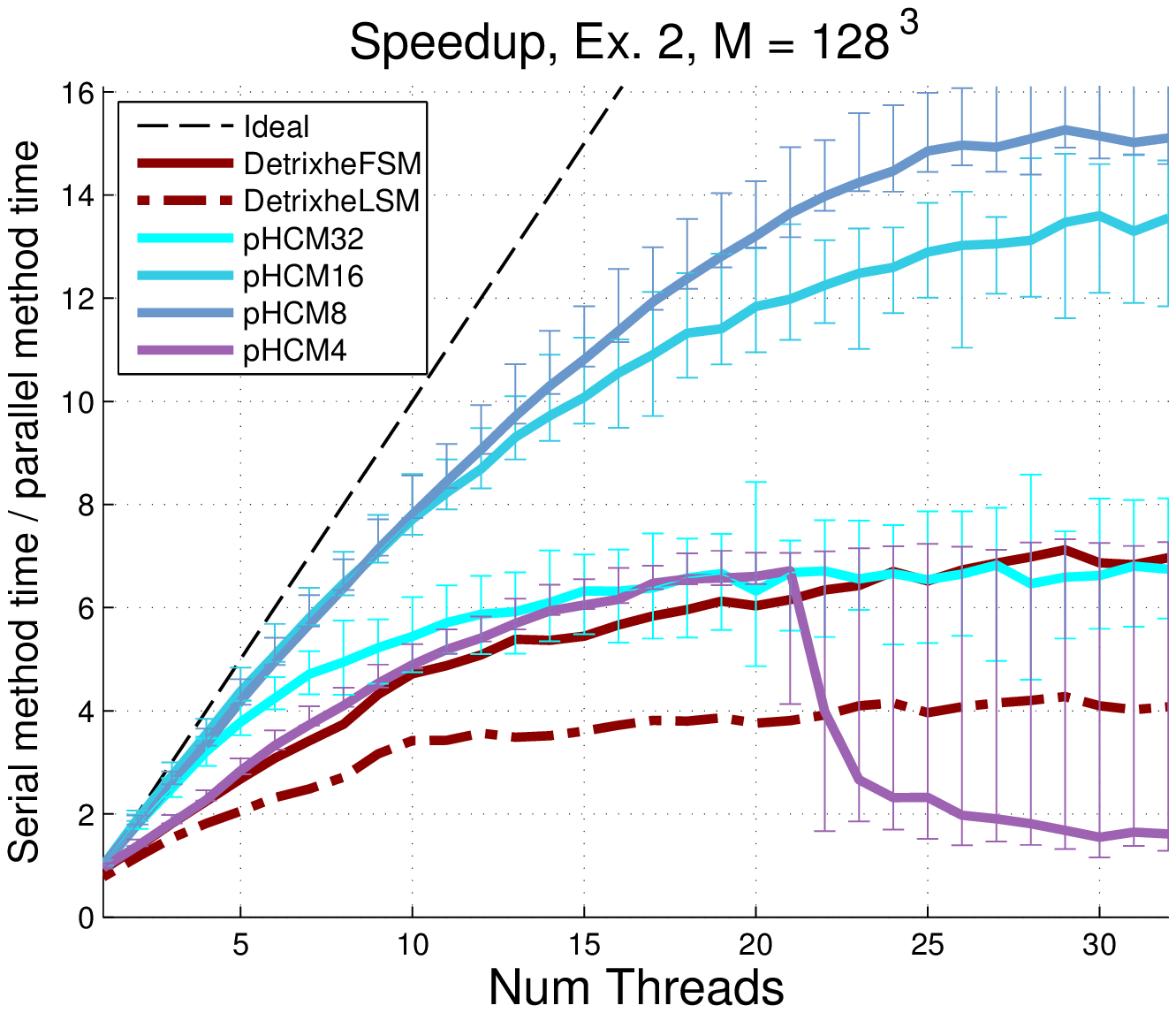}&
\hspace{-.2in}
\includegraphics[scale = .45]{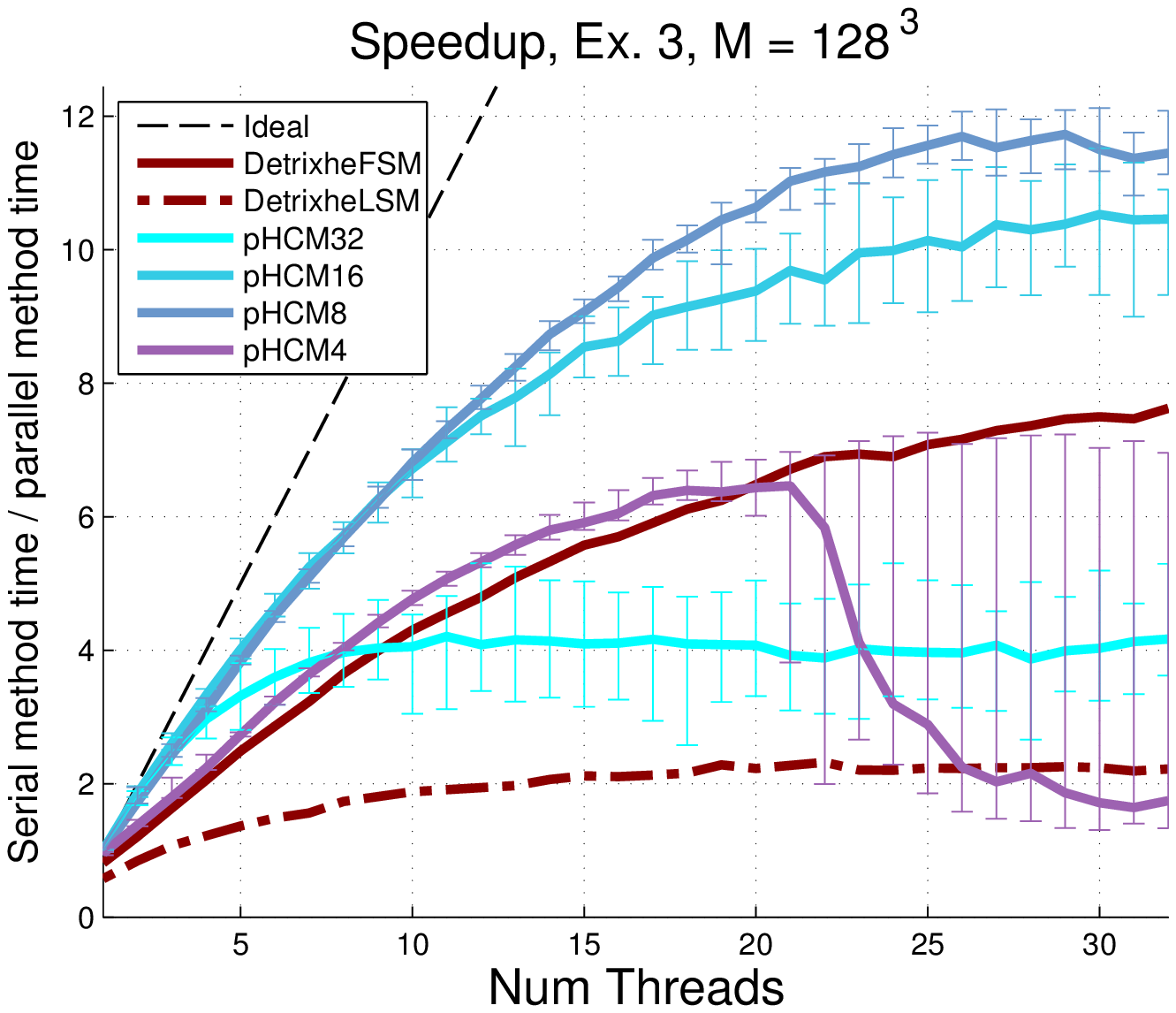}\\

A&B&C\\

\includegraphics[scale = .45]{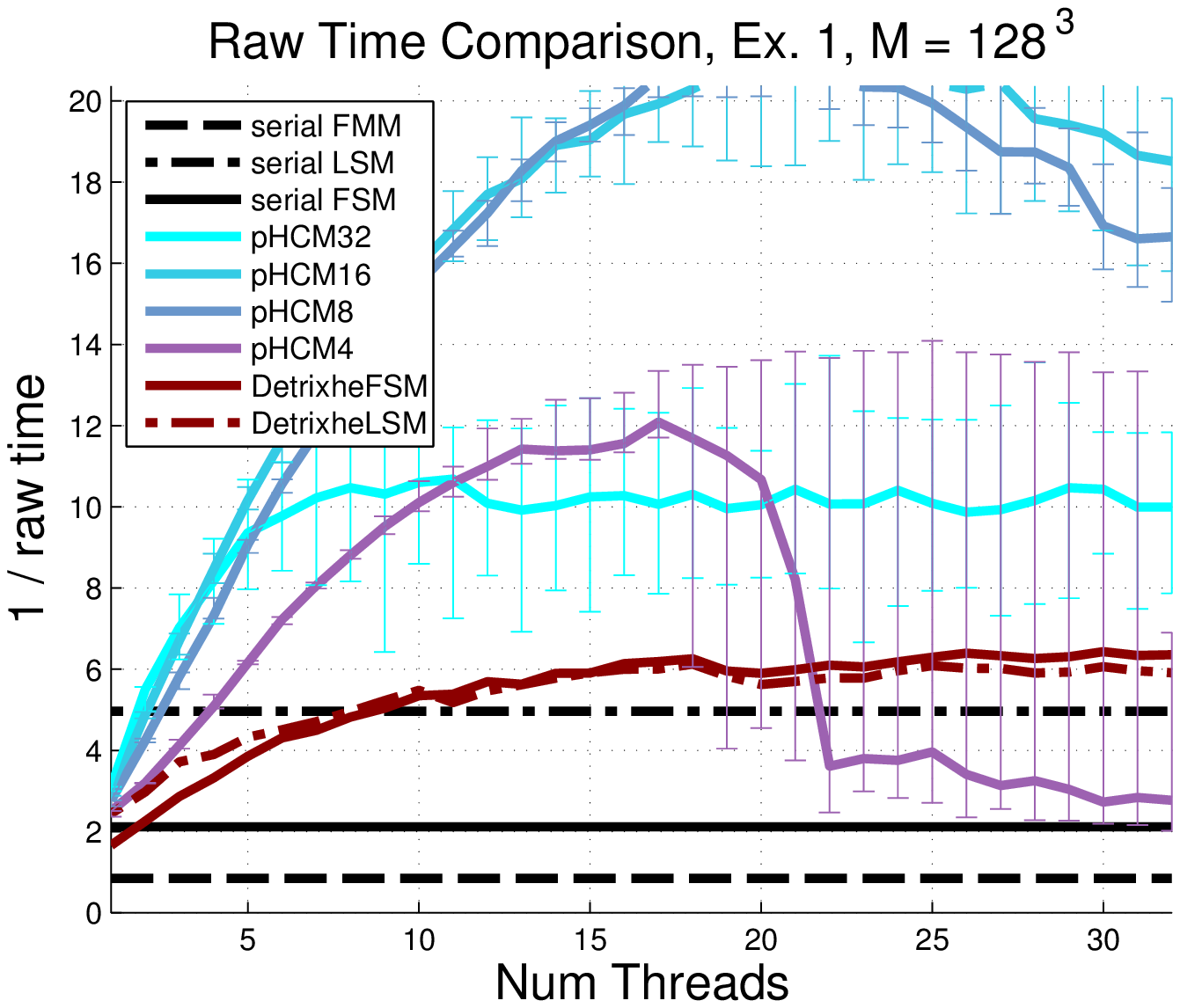}&
\hspace{-.2in}
\includegraphics[scale = .45]{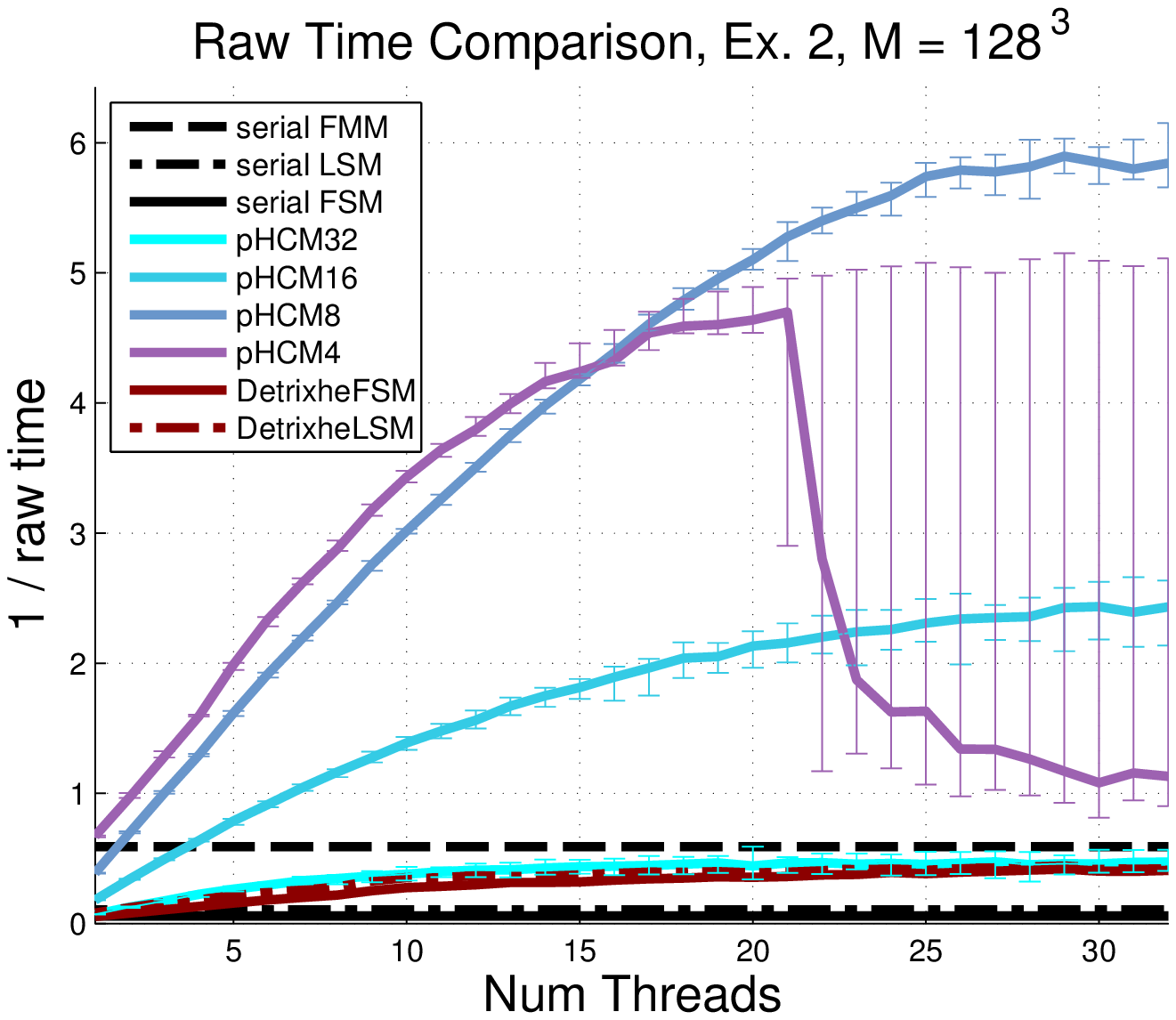}&
\hspace{-.2in}
\includegraphics[scale = .45]{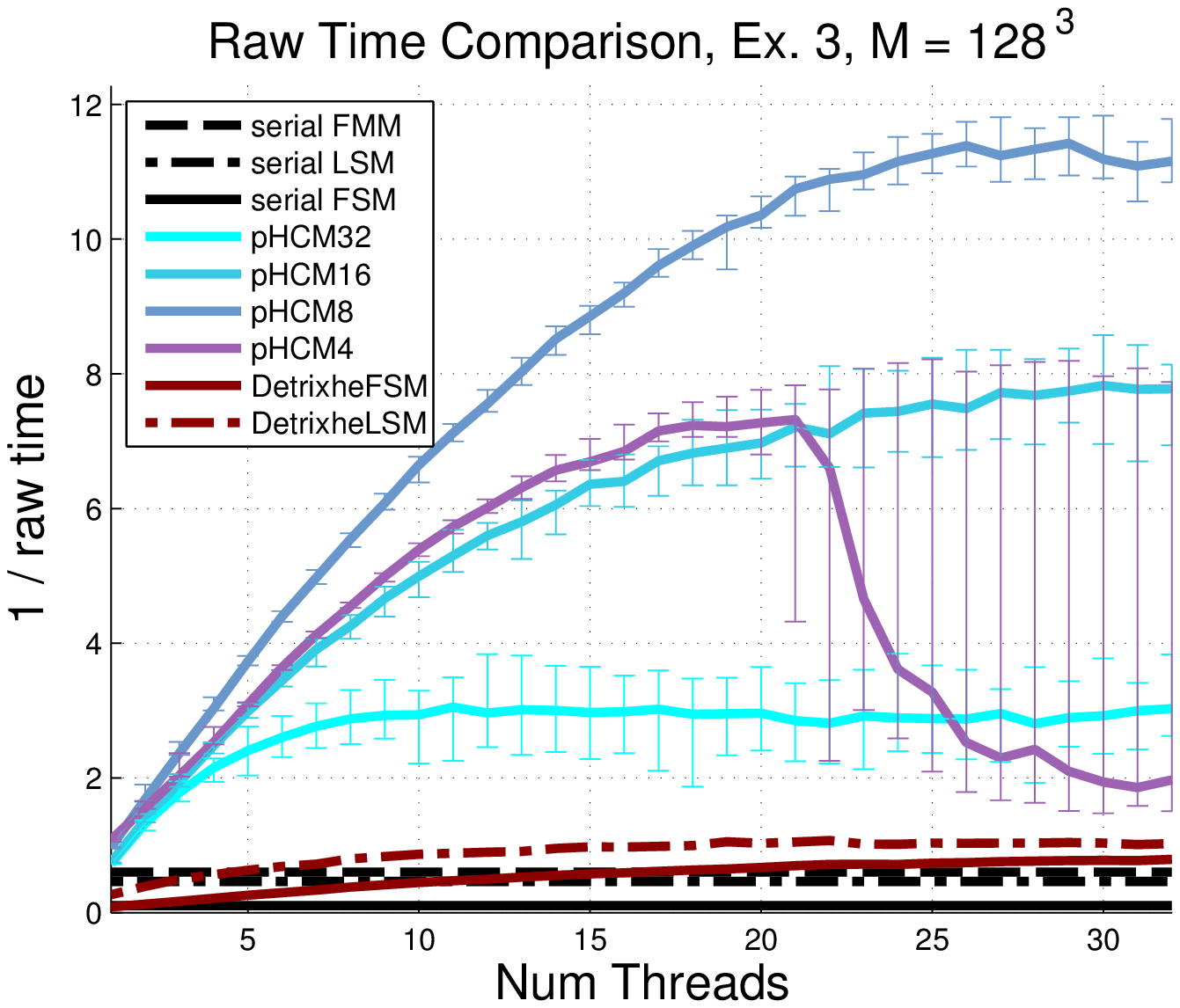}\\

D&E&F\\

\end{array}
$
\caption{\footnotesize Scaling and performance for pHCM at $M = 128^3$.  The first column has $F \equiv 1$, the second has $F = 1 + .5\sin{(20\pi x)}\sin{(20\pi y)}\sin{(20\pi z)}$, and the third has $F = 1 + .99\sin{(2\pi x)}\sin{(2\pi y)}\sin{(2\pi z)}$.}

\label{fig:parallel_128}
\end{figure}

%COMMENT REMOVED
%COMMENT REMOVED
Choosing the \emph{optimal} cell decomposition for a given problem and grid resolution remains a difficult problem
even for the serial HCM.  But luckily, as shown in Fig. \ref{fig:serial_charts} and in \cite{ChacVlad}, a wide range of medium-sized cells exhibits good serial performance {\bf and} parallelizes sufficiently well (Figures \ref{fig:parallel_charts} and \ref{fig:parallel_128}).  In all cases, the parallelization is better when there is more work per cell (e.g., $r$ is large) and there are enough active cells to keep all processors busy.

\noindent

%COMMENT REMOVED
%COMMENT REMOVED
%COMMENT REMOVED
%COMMENT REMOVED
%COMMENT REMOVED
%COMMENT REMOVED
%COMMENT REMOVED
%COMMENT REMOVED
%COMMENT REMOVED
%COMMENT REMOVED
%COMMENT REMOVED
%COMMENT REMOVED
%COMMENT REMOVED
%COMMENT REMOVED
%COMMENT REMOVED
%COMMENT REMOVED
%COMMENT REMOVED
%COMMENT REMOVED

\iffullversion

\pagebreak
\subsection{Performance with ``early sweep terminations"}
\label{ss:earlyTerm}
All sweeping methods can be accelerated by stopping the iterations once the maximum change over gridpoint values is less than or equal to a certain threshold $\kappa \geq 0$.  If $\kappa > 0$, the method will terminate ``early'', and the output will be different than the true solution of the discretized system (\ref{eq:Eik_discr}).
Ideally, $\kappa$ should be chosen based on the $L_{\infty}$-norm discretization error, but since the latter is a priori unknown, a common practical approach is to use a small heuristically selected constant (e.g., \cite{Zhao}).  We note that, for a fixed $\kappa>0$, the number of needed iterations can be quite different for different $h$,
and there is currently no proof that the early-terminated numerical values are within $\kappa$ from the correct solution; see \cite{ChacVlad_expanded}.
%COMMENT REMOVED

All results reported in previous subsections were obtained with $\kappa=0$, but on a computer with finite precision the iterations stop when the gridpoint value changes fall below the machine epsilon.  I.e., for ``double precision'' computations this is equivalent to using $\kappa=2^{-52} \approx 2.2 \times 10^{-16}.$

Here we repeat the same 3 examples but with $\kappa= 10^{-8}$ to force an early sweeping termination, keeping all other parameters the same as in subsections \ref{ss:serial_performance}-\ref{ss:parallel_performance}.  As expected, this modification results in faster termination for FSM, LSM, DFSM, and DLSM (see Figure \ref{fig:early}).
For a fair comparison, in HCM/pHCM we now terminate the sweeping within a cell when the maximum change in a gridpoint's value is less than $\kappa$.
We also add an additional condition on line 11 of Algorithm \ref{alg:gpUpdate}: if a gridpoint value changes by less than $\kappa$, then the procedure on line 12 will not be executed (i.e., the adjacent cell will not be marked for update).  For most $r$ values and on most examples, the number of ``updates per gridpoint" done by HCM$r$ decreases when $\kappa = 10^{-8}$ -- yielding the expected decrease in CPU times. However, we have also observed a surprising (and as of now unexplained) work {\em increase} for HCM32 on Example 2 with $M=320^3$.

For the parallel methods, the scaling is about the same (e.g., Figures \ref{fig:early}D and \ref{fig:early}F) or slightly worse (e.g. Figure \ref{fig:early}E) than it was before with $\kappa = 0.$  For pHCM this is not surprising, since there is effectively less work per cell.
However, for most $r$ values, the improvement in HCM still results in faster pHCM execution times (compared to those in Figure \ref{fig:parallel_charts}).

An experimental study of additional errors due to early termination can be found in \cite{ChacVlad_expanded}.

\begin{figure}[H]
\hspace{-.6 in}
$
\begin{array}{ccc}
\includegraphics[scale = .45]{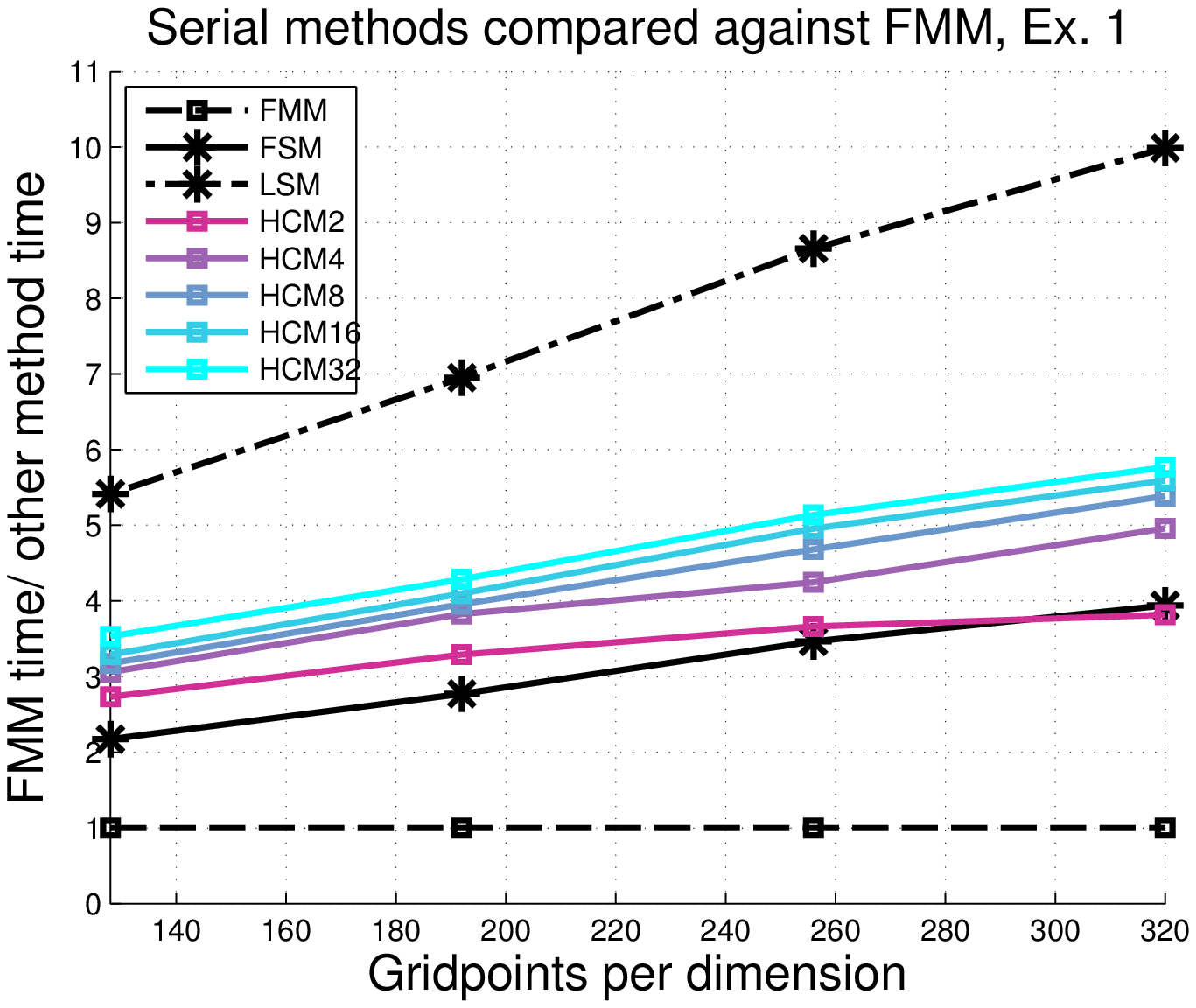}&
\hspace{-.2in}
\includegraphics[scale = .45]{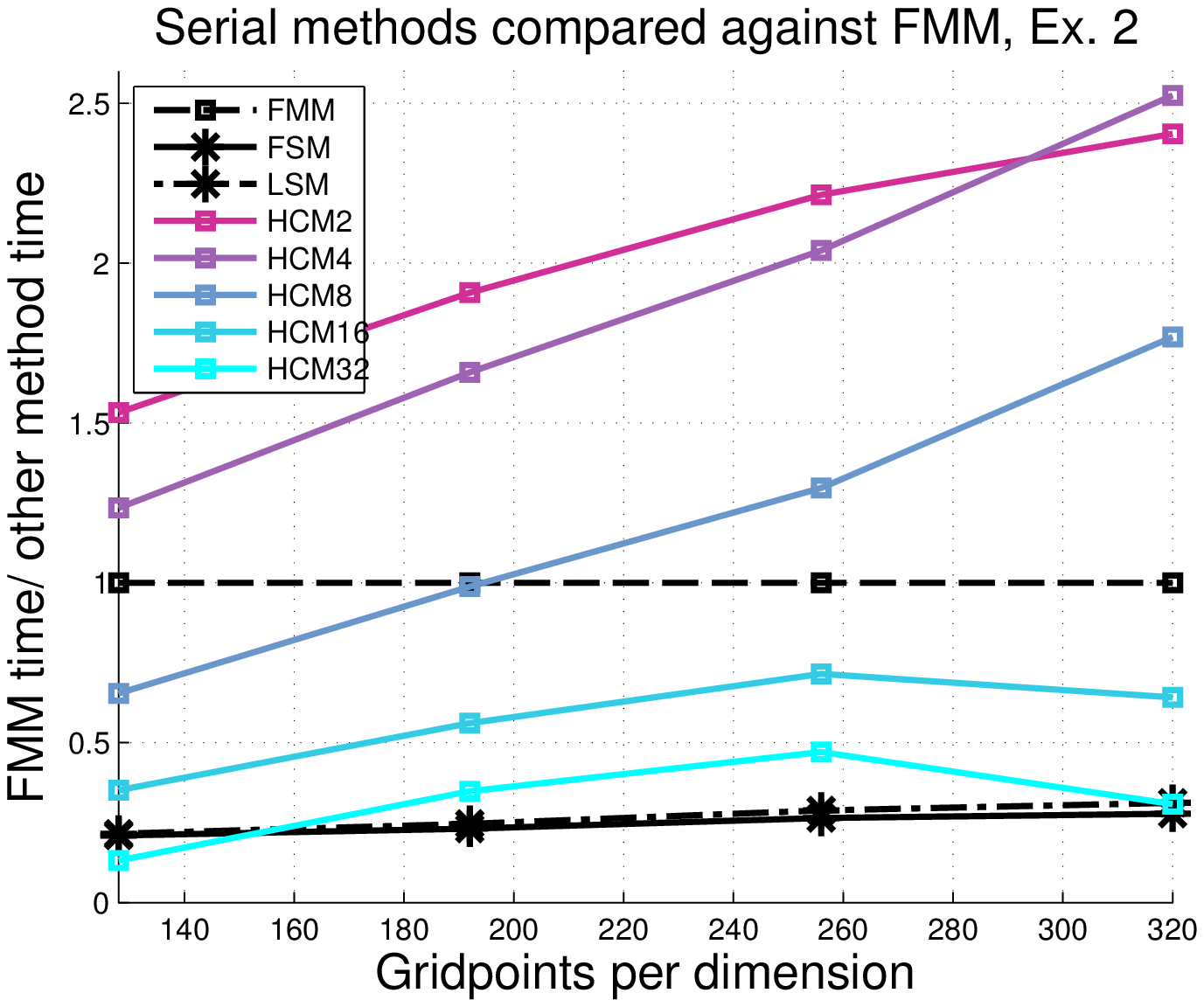}&
\hspace{-.2in}
\includegraphics[scale = .45]{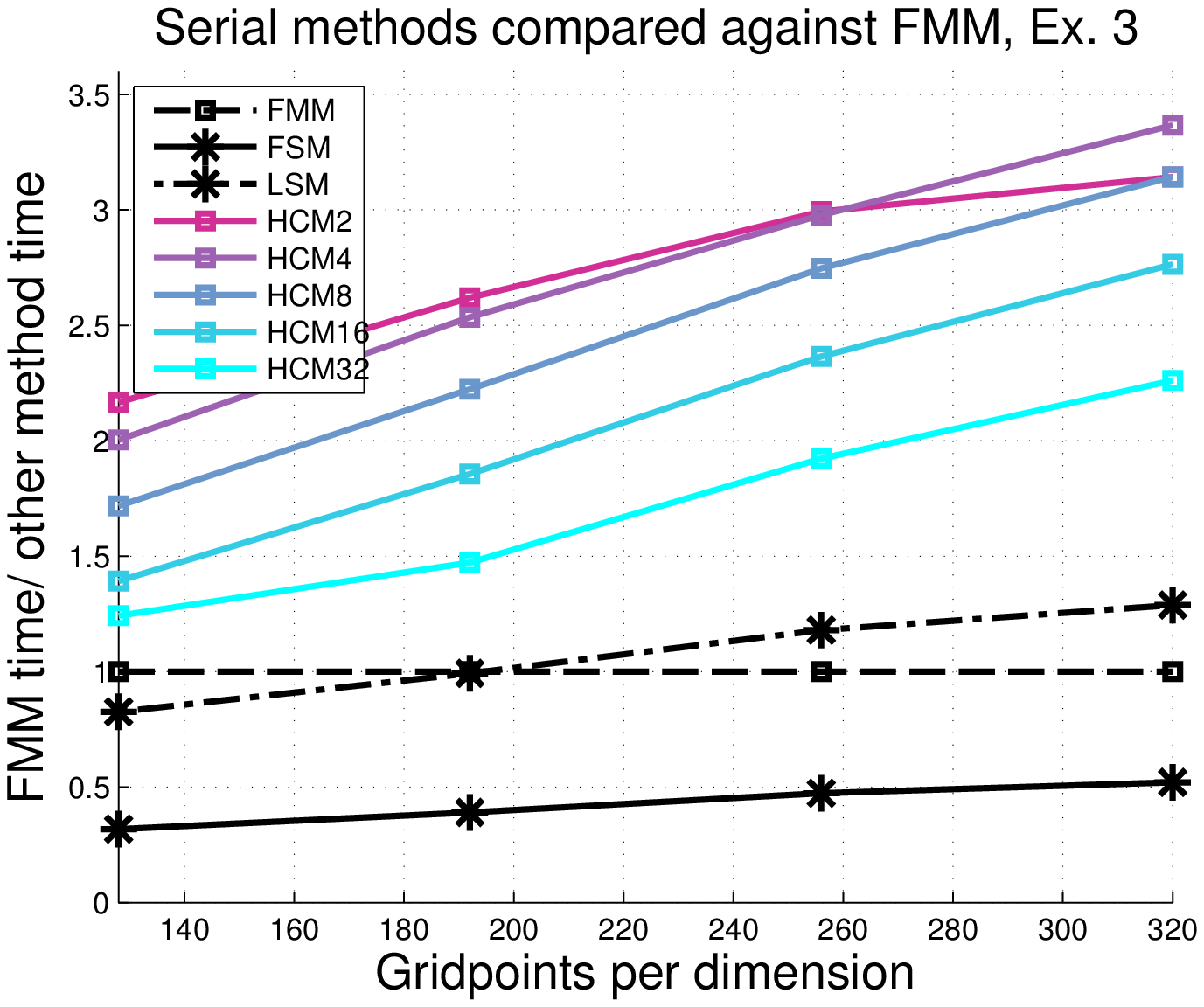}\\

A&B&C\\

\includegraphics[scale = .45]{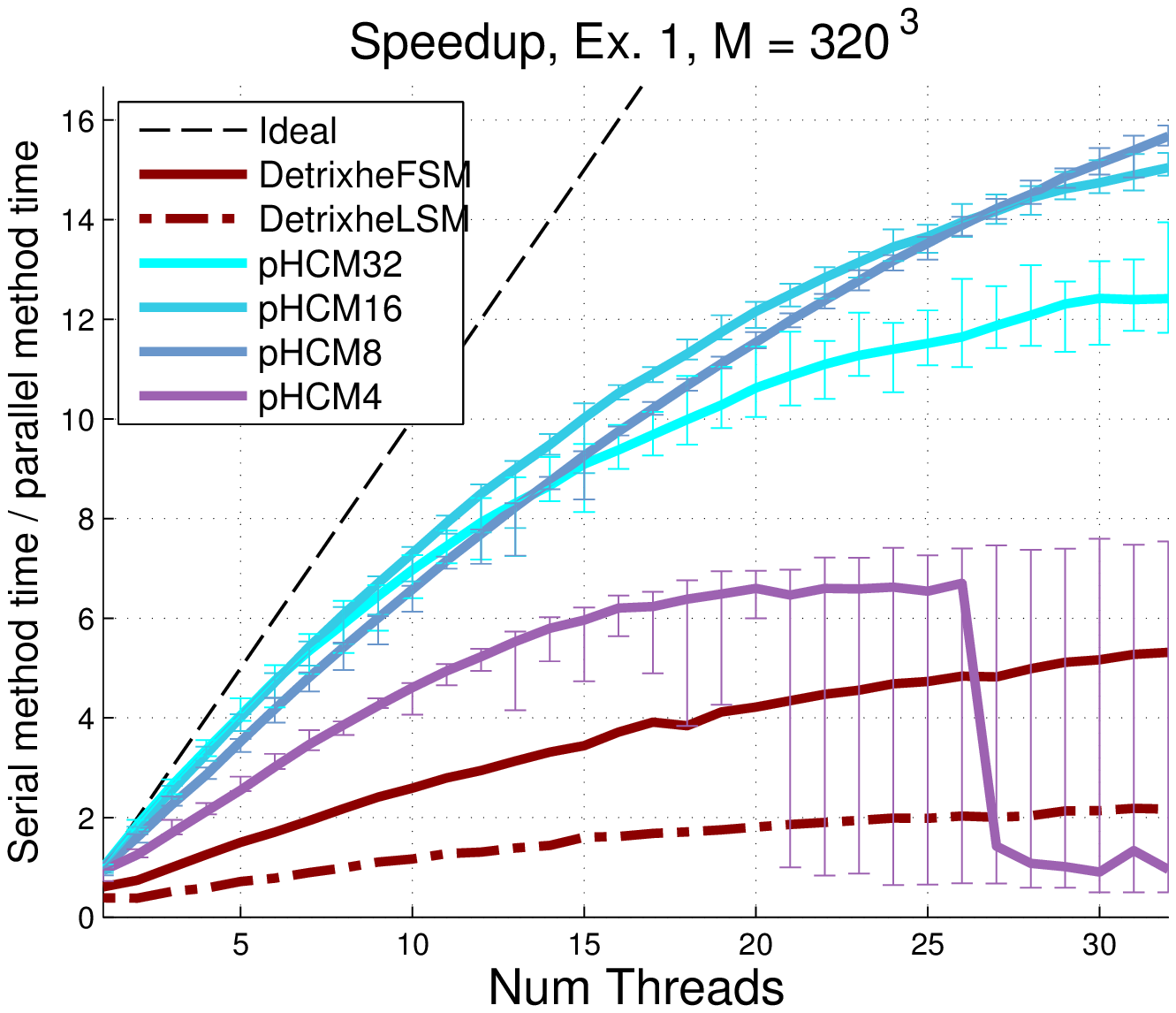}&
\hspace{-.2in}
\includegraphics[scale = .45]{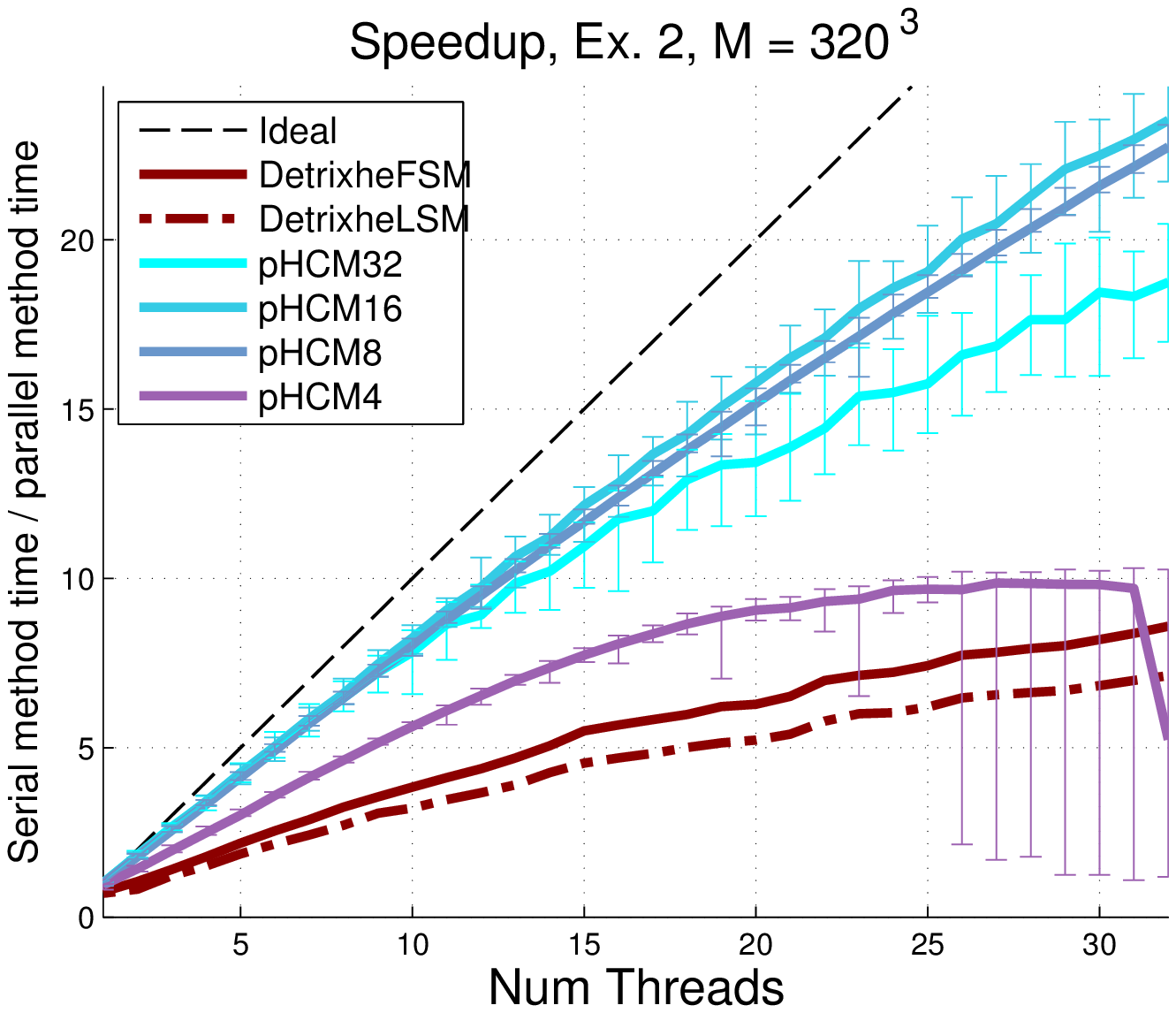}&
\hspace{-.2in}
\includegraphics[scale = .45]{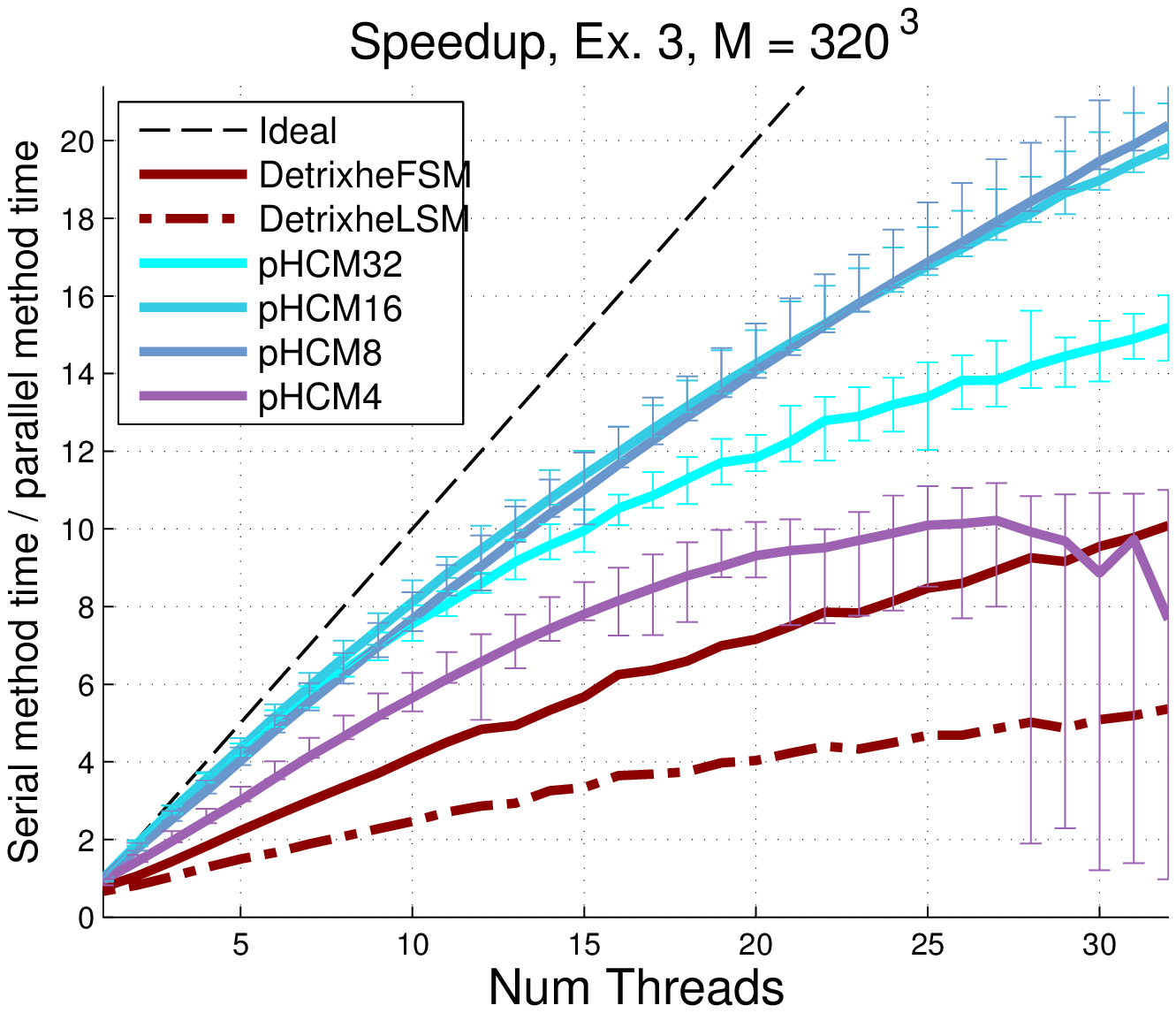}\\

D&E&F\\

\includegraphics[scale = .45]{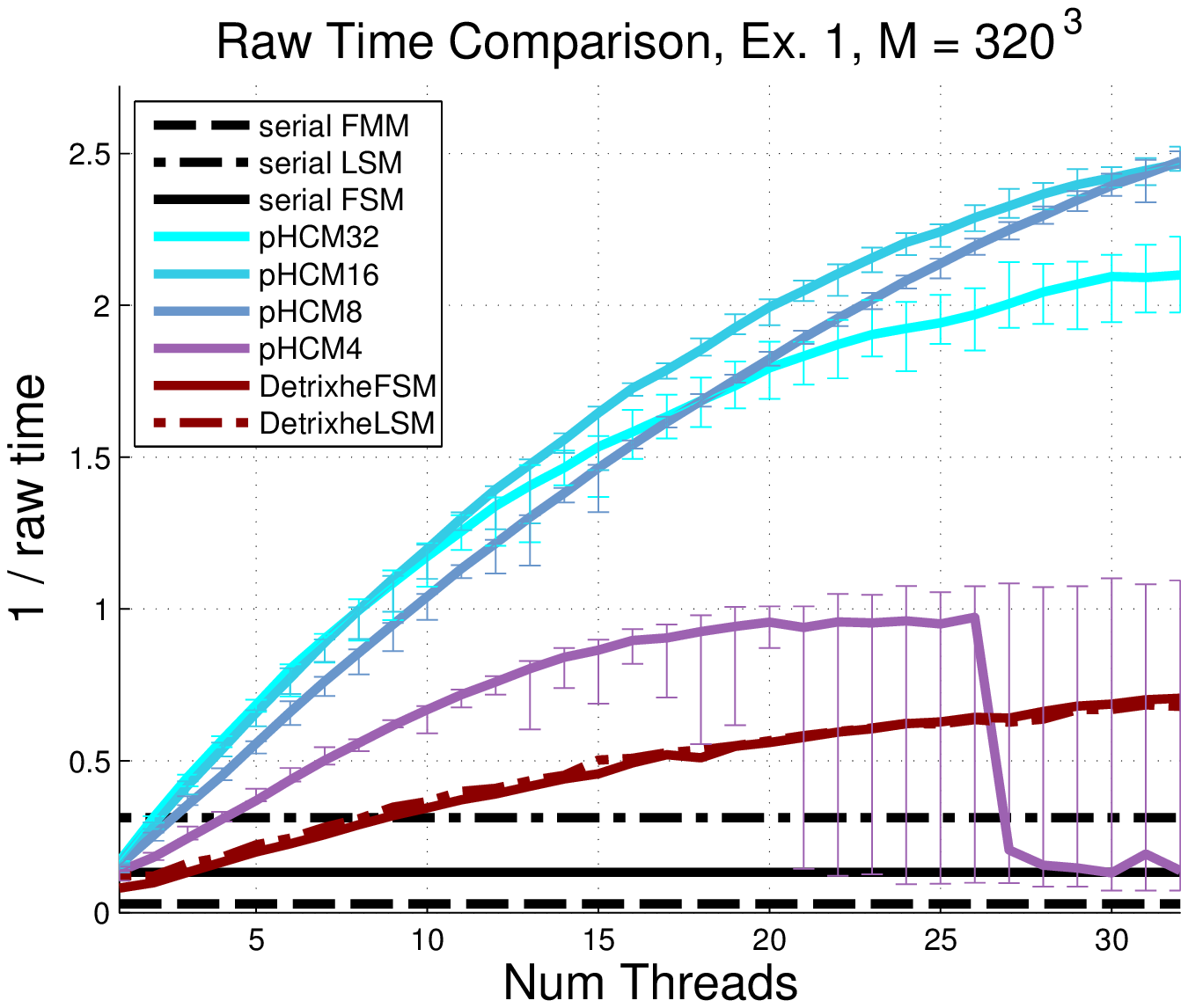}&
\hspace{-.2in}
\includegraphics[scale = .45]{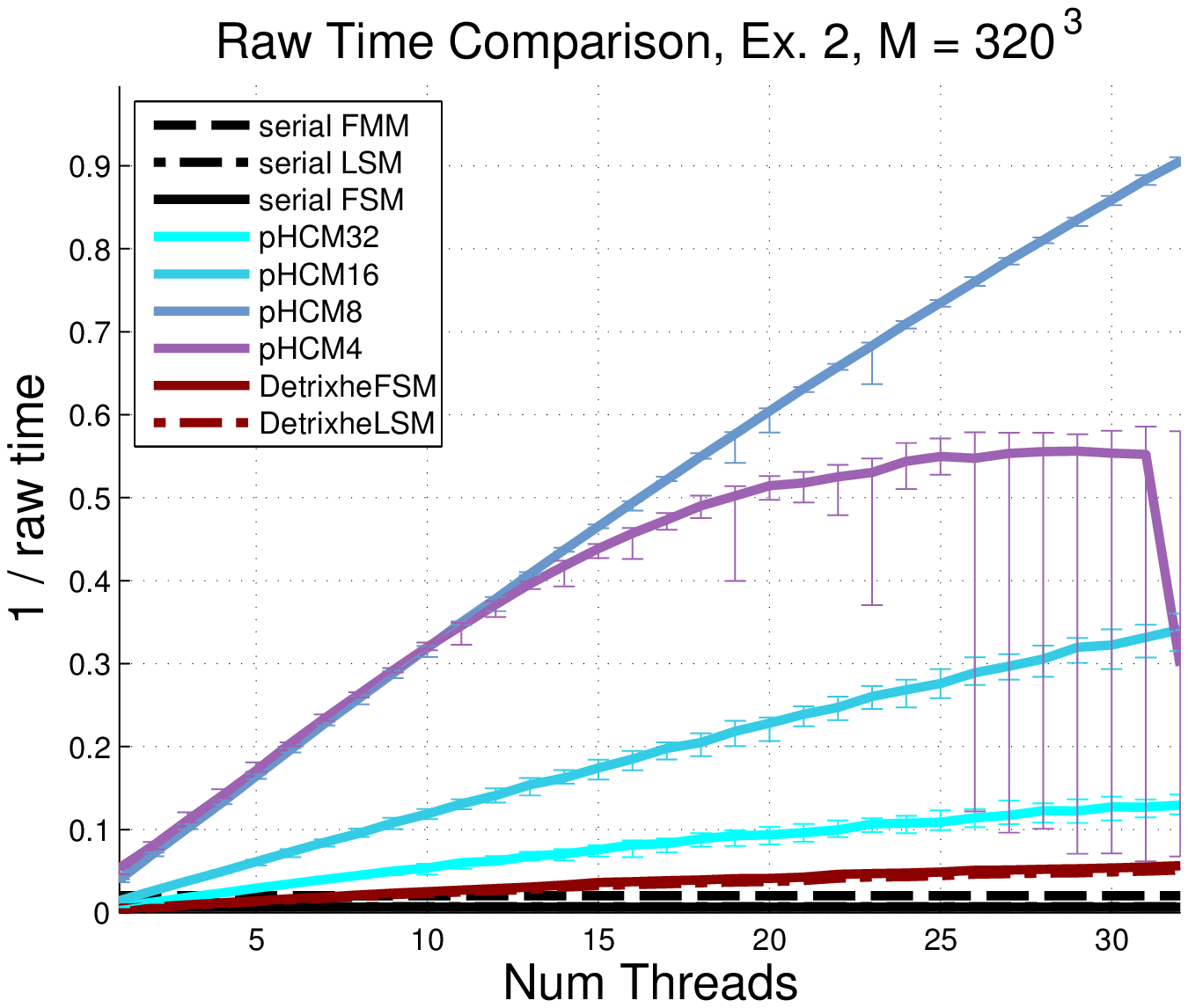}&
\hspace{-.2in}
\includegraphics[scale = .45]{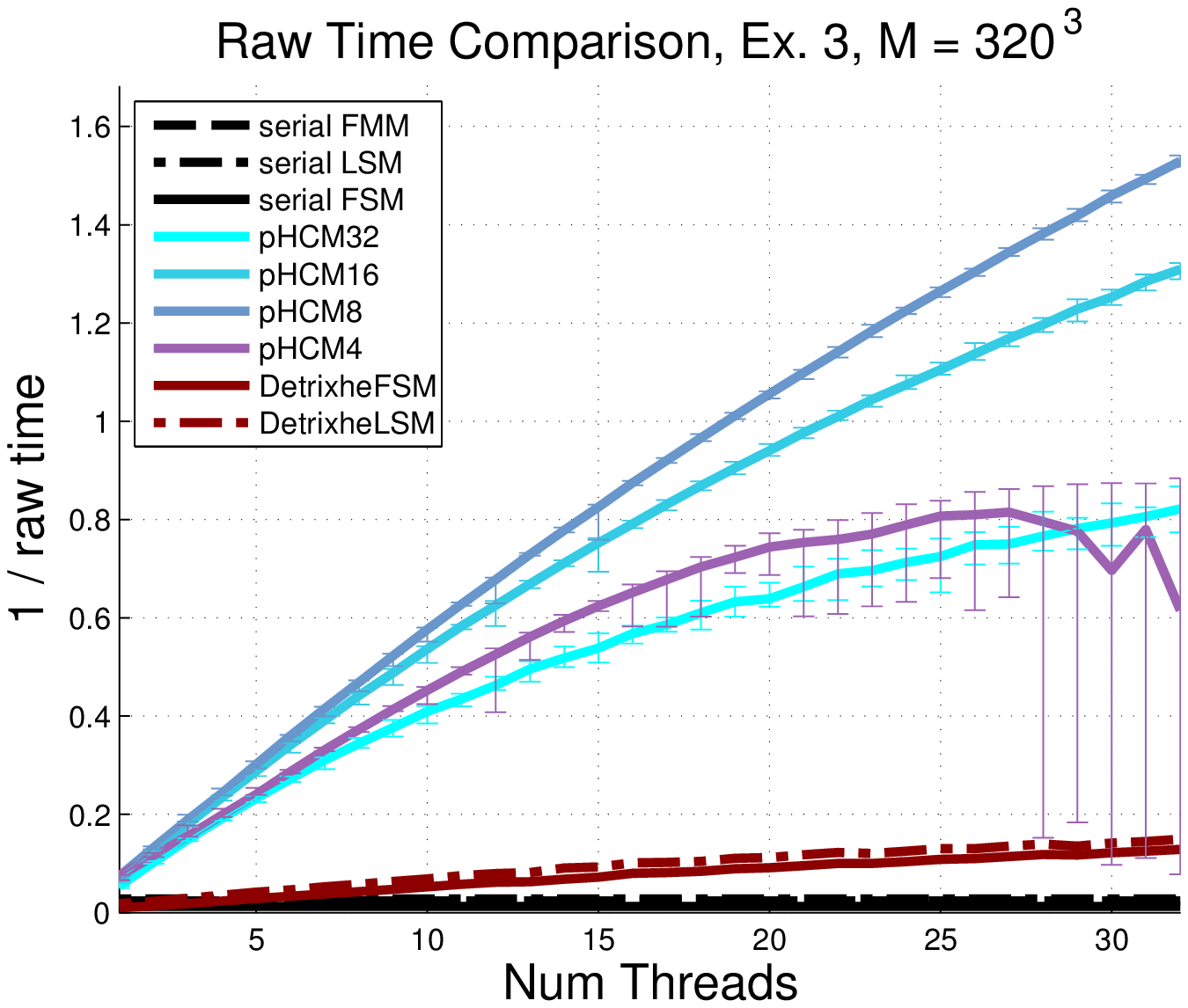}\\

G&H&I

\end{array}
$

\caption{\footnotesize Early termination testing (subsection \ref{ss:earlyTerm}).
Top row: performance of serial methods for different $M$; compare with Figure \ref{fig:serial_charts}.
Two bottom rows: scaling and performance for pHCM at $M = 320^3$; compare with Figure \ref{fig:parallel_charts}.
}
\label{fig:early}
\end{figure}

\pagebreak

\subsection{Performance with ``single precision" data}
\label{ss:single_prec}
In this subsection we repeat the same three experiments but storing/computing
the numerical solution in single precision.  This implementation uses ``float''
instead of ``double'' variables throughout the C++ code.  The results are presented in Figure \ref{fig:float}.

We would expect that in single precision a smaller data footprint would have advantages for high-level cache operations and scaling.  This is mostly true, as illustrated best for DFSM and pHCM on Example 3 (Figure \ref{fig:float}F).  It is also natural to expect that switching to single precision should substantially decrease the total number of needed iterations to convergence, because the iterations stop when the maximum change in values is less than machine epsilon (i.e.,  we are effectively using $\kappa = 2^{-23} \approx 1.2 \times 10^{-7}$).  Tables \ref{tab:sweepConvergence1} and \ref{tab:sweepConvergence2} are a side-by-side comparison of sweeping-convergence data
for Example 2 with $M = 64^3$ under single and double precision.
%COMMENT REMOVED
%COMMENT REMOVED
Based on Table \ref{tab:sweepConvergence2}, it is natural to expect that sweeping in single precision should converge in about 33 sweeps.  Table \ref{tab:sweepConvergence1} shows that this is {\bf not} the case: 53 sweeps are in fact required for convergence.  The reason for this discrepancy is that intermediate computations are also conducted in single precision.  In fact, Table \ref{tab:sweeps} shows that on Ex. 3 with $M = 320^3$, the number of sweeps to convergence is actually higher in single than in double precision.  This helps explain the downward-sloping LSM curve in Figure \ref{fig:float}C.

We note that Table \ref{tab:sweeps} also shows a growth in the number of iterations-to-convergence with $M$ for the sweeping methods on examples 2 and 3 in either single or double precision.

%COMMENT REMOVED

%COMMENT REMOVED

%COMMENT REMOVED

\vfill

\begin{table}[H]
\caption{\footnotesize Number of sweeps for different values of $M$ in double and single precision.}
\centering
\begin{tabular}{|l|l|c|c|c|c|c|}
\hline
& & $64^3$ & $128^3$ & $192^3$ & $256^3$ & $320^3$ \\
\hline
\multirow{2}{*}{Ex. 1} & double & 9 & 9 & 9 & 9 & 9 \\
&single& 9 & 9 & 9 & 9 & 9\\
\hline
\multirow{2}{*}{Ex. 2} & double & 69 & 99 & 131 & 164 & 191 \\
&single& 53 & 88 & 116 & 144 & 173\\
\hline
\multirow{2}{*}{Ex. 3} & double & 42 & 58 & 77 & 107 & 121 \\
&single& 36 & 56 & 89 & 97 & 129\\
\hline

\end{tabular}

\label{tab:sweeps}
\end{table}

\vfill

\begin{figure}[H]
\hspace{-.6 in}
$
\begin{array}{ccc}
\includegraphics[scale = .45]{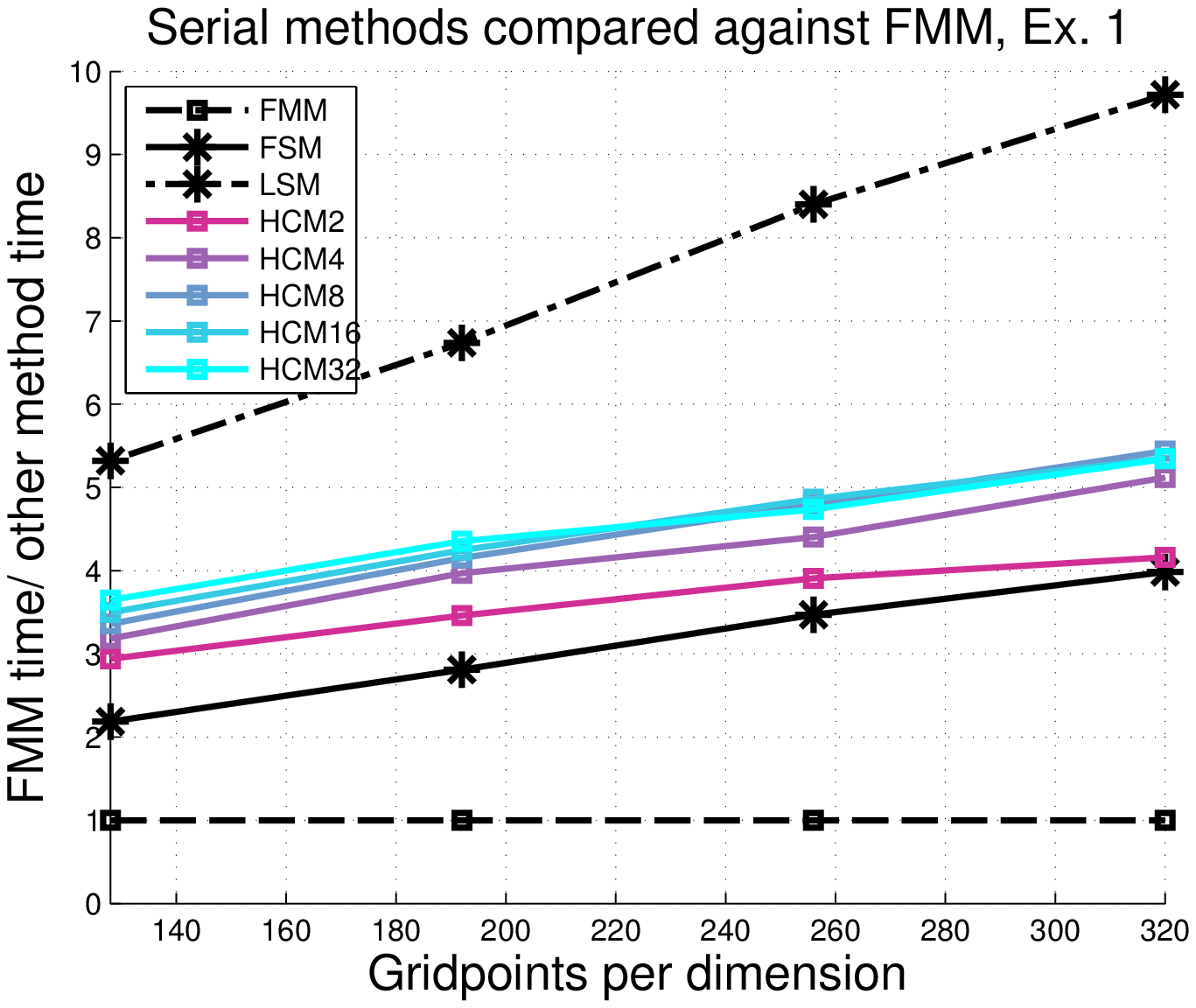}&
\hspace{-.2in}
\includegraphics[scale = .45]{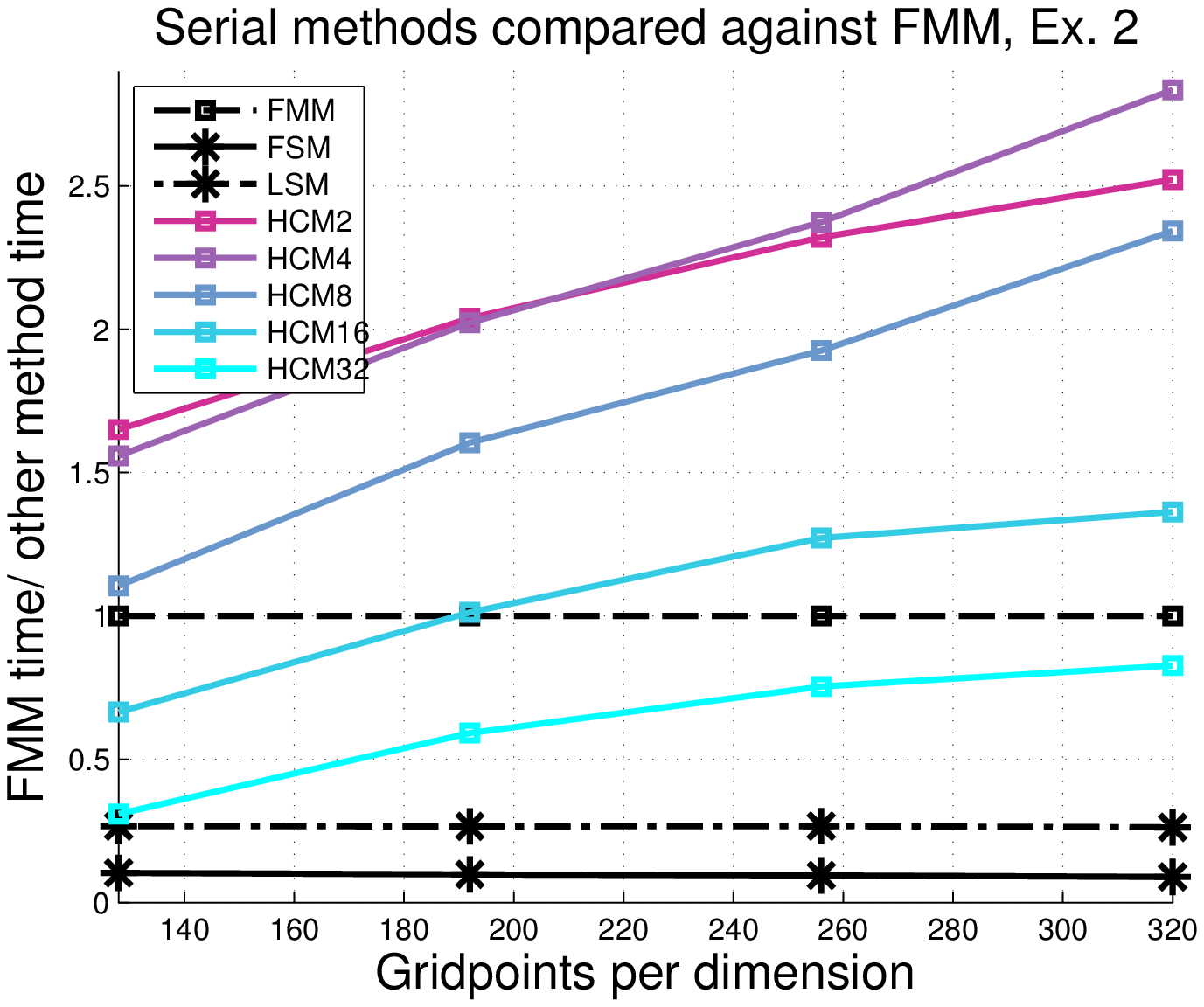}&
\hspace{-.2in}
\includegraphics[scale = .45]{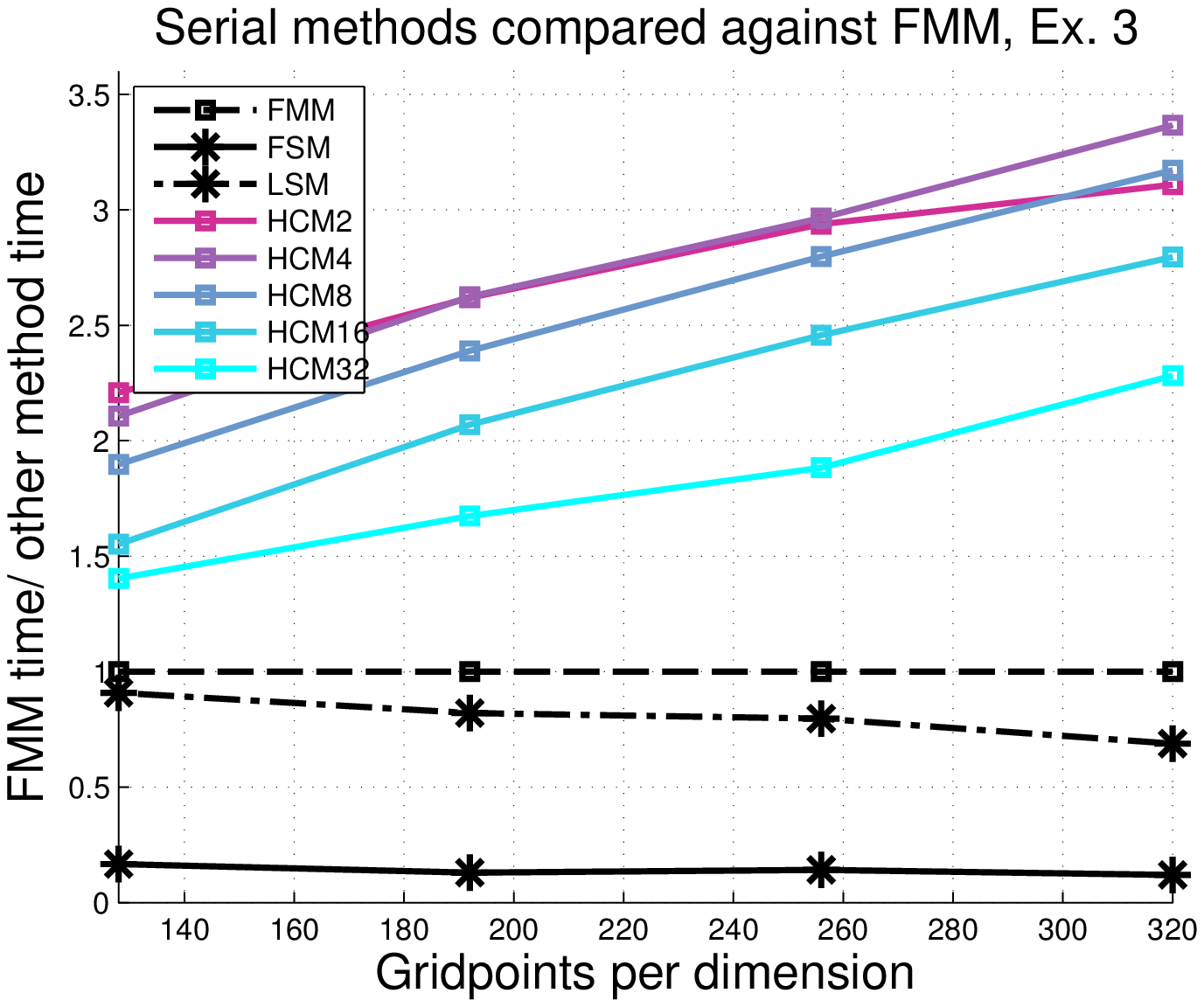}\\

A&B&C\\

\includegraphics[scale = .45]{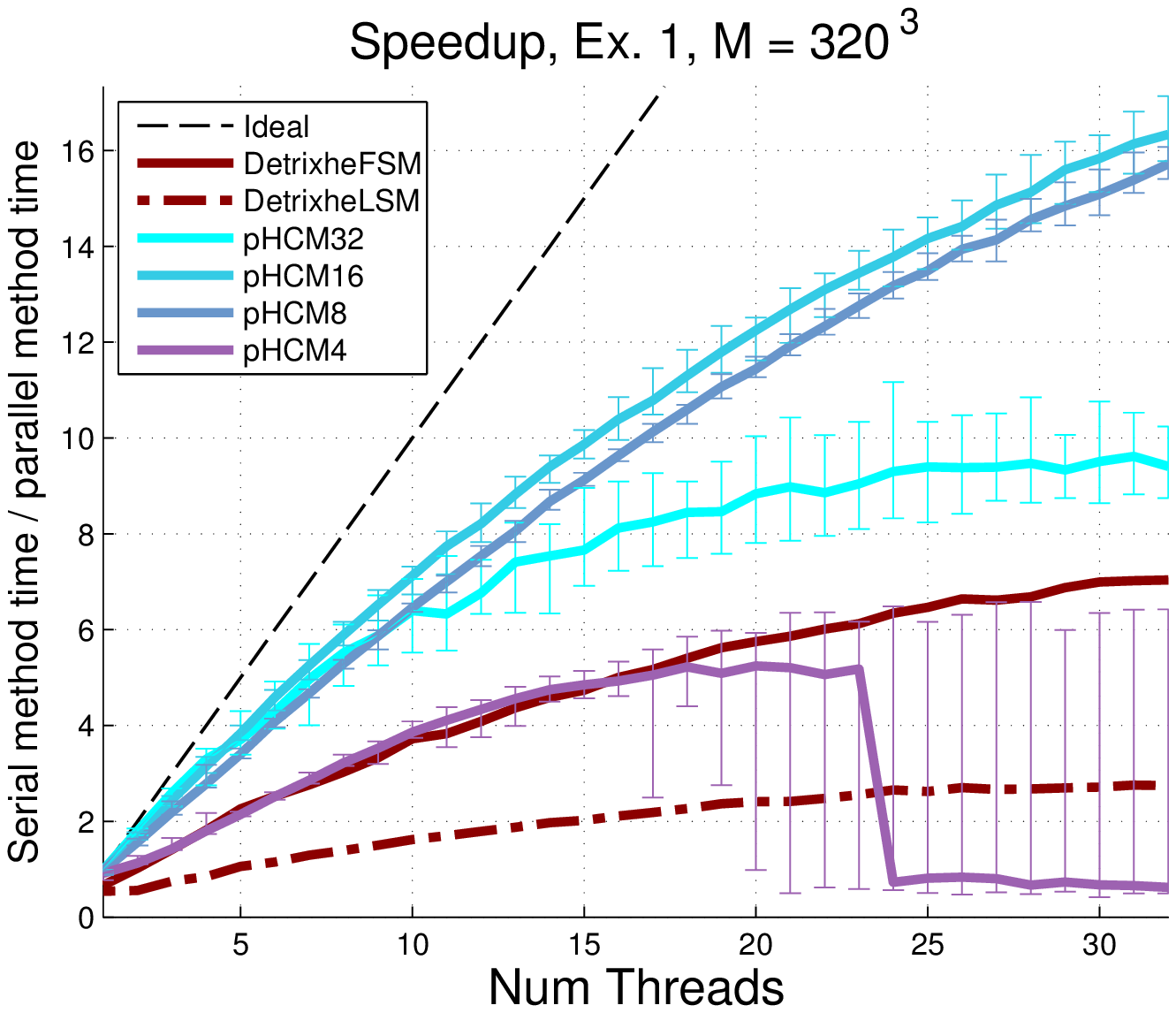}&
\hspace{-.2in}
\includegraphics[scale = .45]{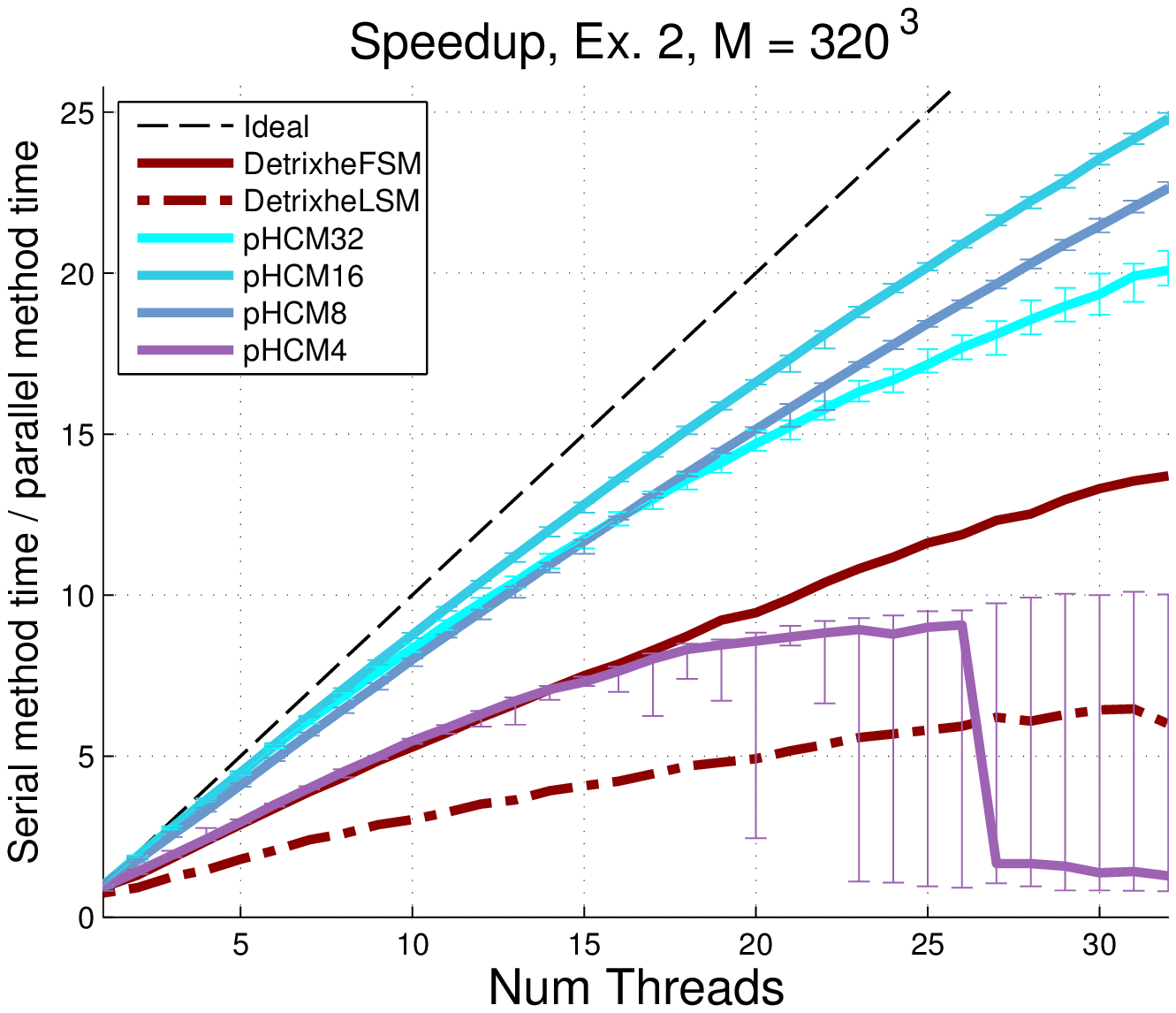}&
\hspace{-.2in}
\includegraphics[scale = .45]{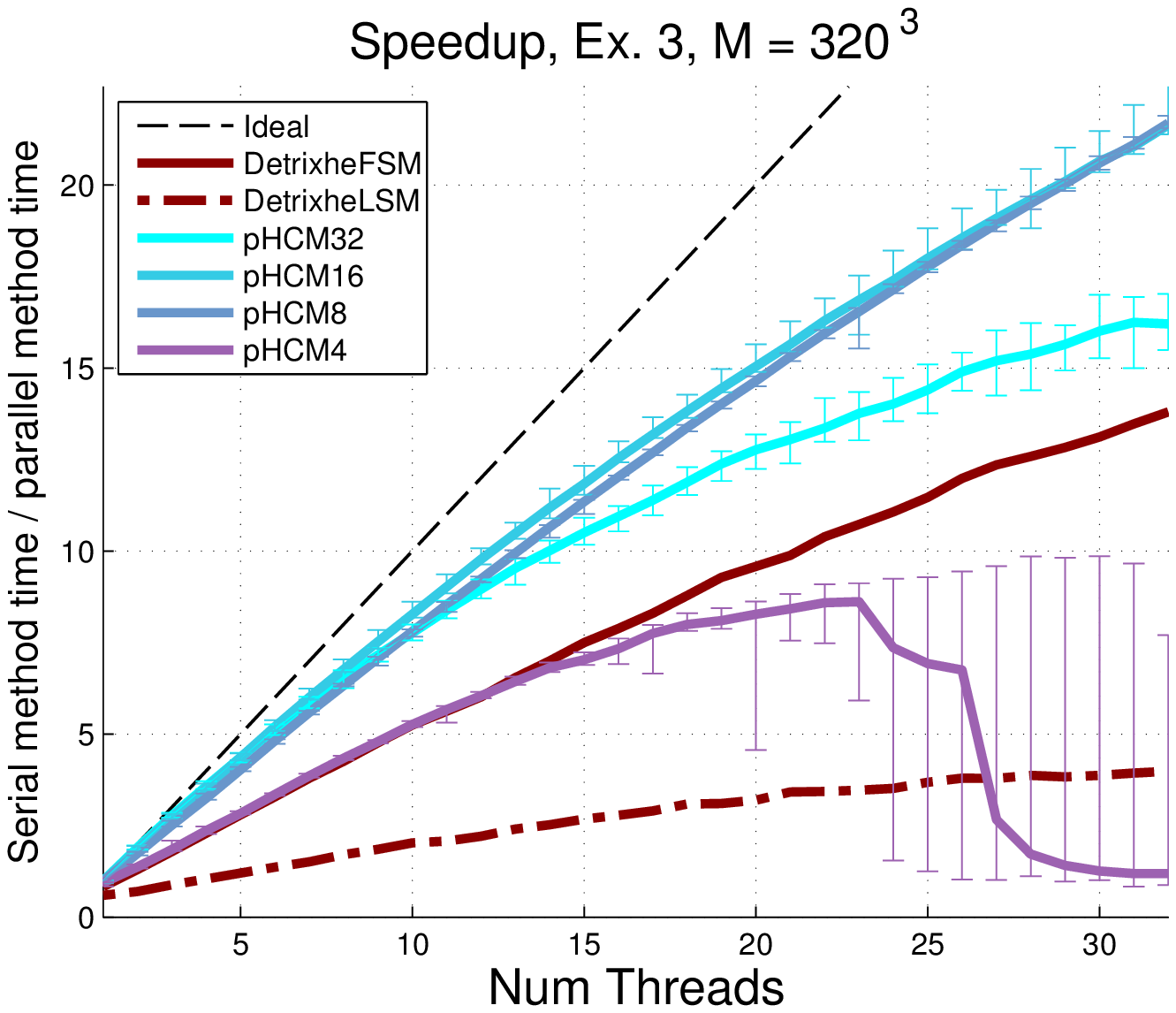}\\

D&E&F\\

\includegraphics[scale = .45]{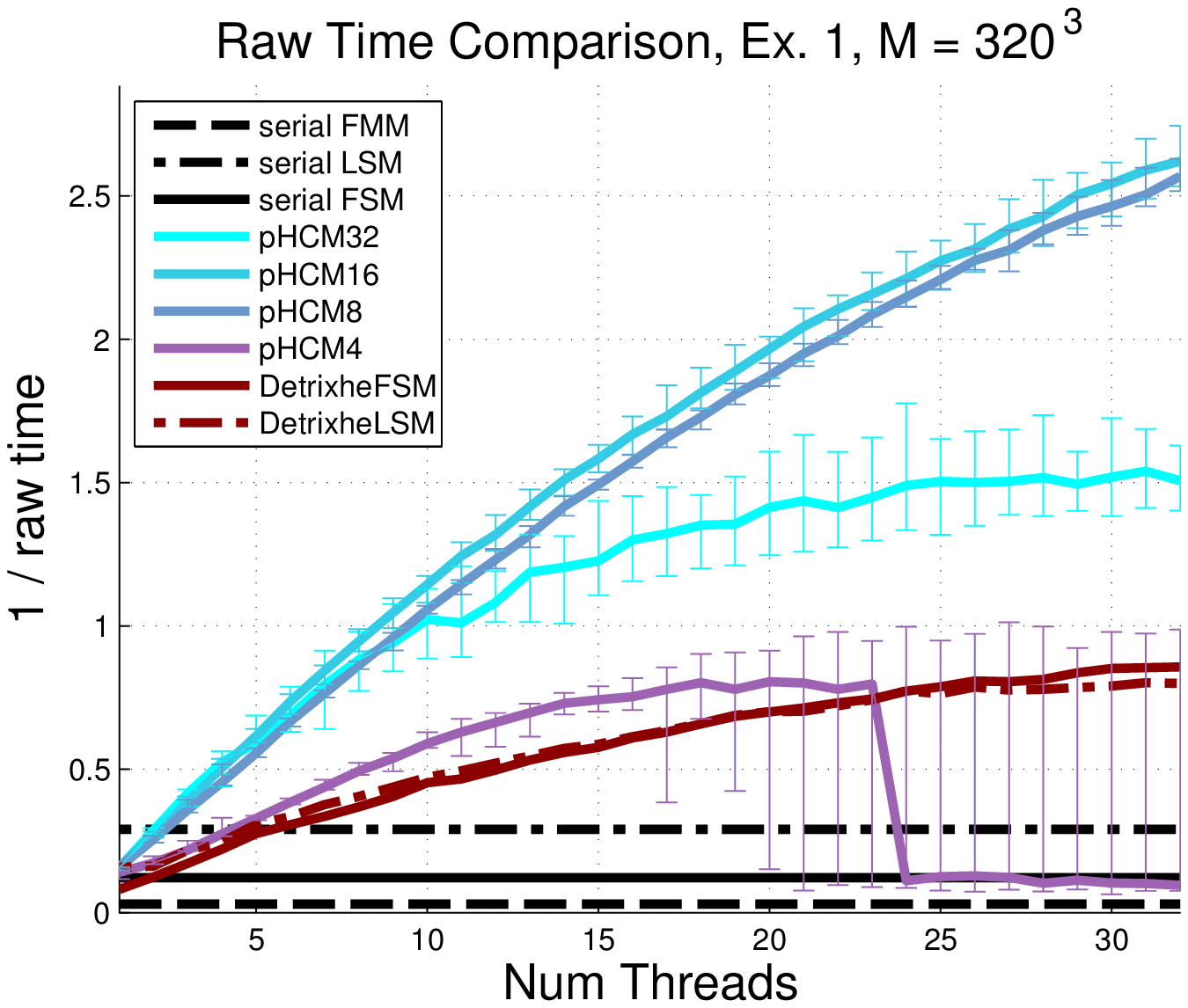}&
\hspace{-.2in}
\includegraphics[scale = .45]{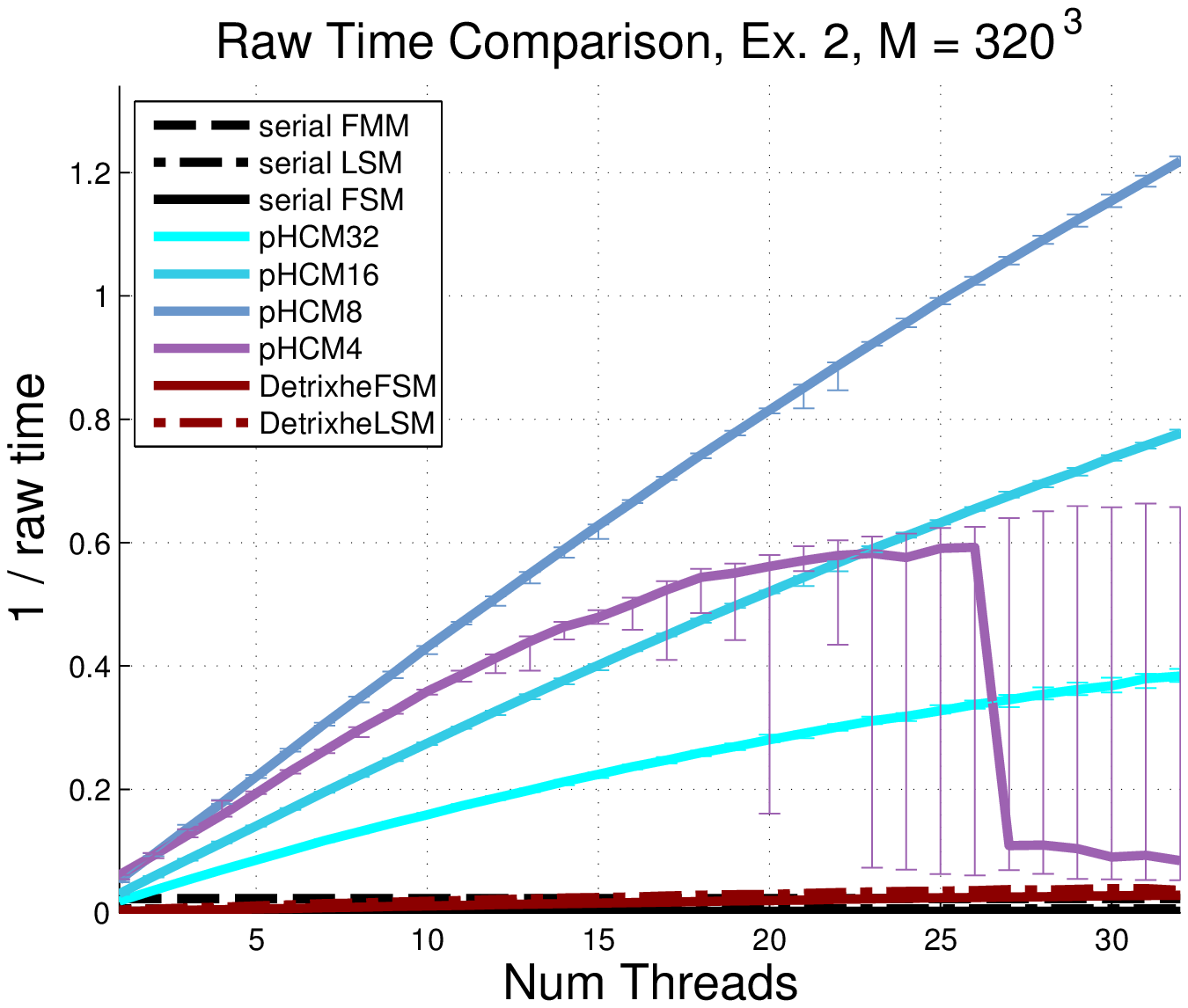}&
\hspace{-.2in}
\includegraphics[scale = .45]{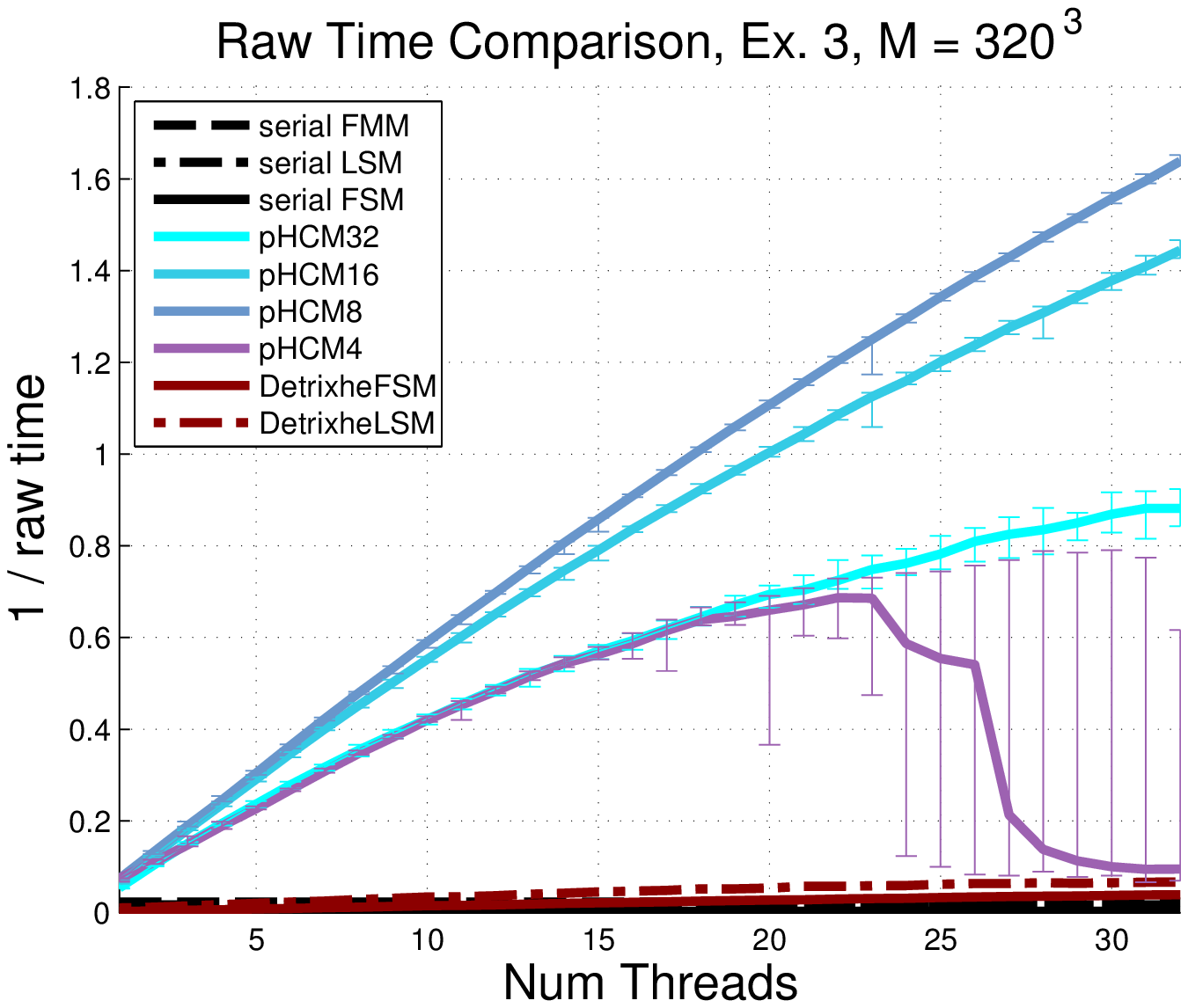}\\

G&H&I

\end{array}
$

\caption{\footnotesize Single precision testing (subsection \ref{ss:single_prec}).
Top row: performance of serial methods for different $M$; compare with Figure \ref{fig:serial_charts}.
Two bottom rows: scaling and performance for pHCM at $M = 320^3$; compare with Figure \ref{fig:parallel_charts}.
}
\label{fig:float}
\end{figure}

%COMMENT REMOVED
%COMMENT REMOVED
%COMMENT REMOVED
\begin{table}[H] \footnotesize
\parbox{.45\linewidth}{
\centering
\caption{Single precision}

\begin{tabular}{|c|c|c|}
\hline
  \textbf{sweep \#} & \textbf{max change} & \textbf{\% grid changing} \\
\hline
1 & 1e+09 & 15\\
\hline 2 & 1e+09 & 23.6\\
\hline 3 & 1e+09 & 43.6\\
\hline 4 & 1e+09 & 42.8\\
\hline 5 & 1e+09 & 75.4\\
\hline 6 & 0.258 & 74\\
\hline 7 & 0.242 & 72.1\\
\hline 8 & 0.156 & 73.4\\
\hline 9 & 0.00248 & 69.7\\
\hline 10 & 0.00155 & 69\\
\hline 11 & 0.00213 & 67.6\\
\hline 12 & 0.00151 & 67.9\\
\hline 13 & 0.00151 & 63.7\\
\hline 14 & 0.00147 & 60.8\\
\hline 15 & 0.00111 & 57.5\\
\hline 16 & 0.000641 & 55.3\\
\hline 17 & 0.000216 & 52.7\\
\hline 18 & 0.000104 & 48.7\\
\hline 19 & 0.000165 & 44.2\\
\hline 20 & 9.66e-05 & 40.5\\
\hline 21 & 0.0001 & 37.8\\
\hline 22 & 8.12e-05 & 31.4\\
\hline 23 & 5.51e-05 & 27.4\\
\hline 24 & 2.31e-05 & 22.8\\
\hline 25 & 6.74e-06 & 20.8\\
\hline 26 & 3.58e-06 & 16.9\\
\hline 27 & 4.71e-06 & 15.1\\
\hline 28 & 2.86e-06 & 12.6\\
\hline 29 & 3.28e-06 & 11.9\\
\hline 30 & 2.86e-06 & 8.41\\
\hline 31 & 2.8e-06 & 7.51\\
\hline 32 & 2.74e-06 & 6.24\\
\hline 33 & 2.44e-06 & 4.48\\
\hline 34 & 2.26e-06 & 3.54\\
\hline 35 & 1.85e-06 & 3.15\\
\hline 36 & 2.15e-06 & 2.4\\
\hline 37 & 2.38e-06 & 2.07\\
\hline 38 & 1.67e-06 & 1.14\\
\hline 39 & 1.67e-06 & 1.04\\
\hline 40 & 1.61e-06 & 0.668\\
\hline 41 & 1.79e-06 & 0.552\\
\hline 42 & 1.19e-06 & 0.335\\
\hline 43 & 1.73e-06 & 0.367\\
\hline 44 & 1.43e-06 & 0.146\\
\hline 45 & 1.13e-06 & 0.053\\
\hline 46 & 9.54e-07 & 0.0202\\
\hline 47 & 1.13e-06 & 0.0153\\
\hline 48 & 8.34e-07 & 0.013\\
\hline 49 & 8.94e-07 & 0.00687\\
\hline 50 & 4.17e-07 & 0.00305\\
\hline 51 & 9.54e-07 & 0.00305\\
\hline 52 & 2.98e-07 & 0.00114\\
\hline 53 & 0 & 0\\
\hline & & \\
\hline & & \\
\hline & & \\
\hline
\end{tabular}
\label{tab:sweepConvergence1}
}
%COMMENT REMOVED
%COMMENT REMOVED
\hfill
\parbox{.45\linewidth}{
\centering
\caption{Double precision}

\begin{tabular}{|c|c|c|}
\hline
\hline \textbf{sweep \#} & \textbf{max change} & \textbf{\% grid changing} \\
\hline
1 & 1e+09 & 15\\
\hline 2 & 1e+09 & 28.7\\
\hline 3 & 1e+09 & 54.5\\
\hline 4 & 1e+09 & 56.3\\
\hline 5 & 1e+09 & 87.1\\
\hline 6 & 0.258 & 88.3\\
\hline 7 & 0.242 & 94.3\\
\hline 8 & 0.156 & 98\\
\hline 9 & 0.00248 & 86.9\\
\hline 10 & 0.00155 & 88.2\\
\hline 11 & 0.00213 & 93.7\\
\hline 12 & 0.00151 & 96.7\\
\hline 13 & 0.00151 & 85.3\\
\hline 14 & 0.00147 & 86.3\\
\hline 15 & 0.00111 & 90.8\\
\hline 16 & 0.00064 & 94.7\\
\hline 17 & 0.000217 & 84\\
\hline 18 & 0.000105 & 84.4\\
\hline 19 & 0.000165 & 86.9\\
\hline 20 & 9.68e-05 & 89.6\\
\hline 21 & 0.0001 & 79.3\\
\hline 22 & 8.11e-05 & 79.6\\
\hline 23 & 5.57e-05 & 80.1\\
\hline 24 & 2.25e-05 & 82.4\\
\hline 25 & 7.04e-06 & 73.9\\
\hline 26 & 2.74e-06 & 72.7\\
\hline 27 & 4.11e-06 & 70.7\\
\hline 28 & 1.52e-06 & 70.3\\
\hline 29 & 1.67e-06 & 62.4\\
\hline 30 & 8.78e-07 & 59.5\\
\hline 31 & 4.82e-07 & 55.8\\
\hline 32 & 1.82e-07 & 52.9\\
\hline 33 & 5.14e-08 & 47.8\\
\hline 34 & 1.7e-08 & 44.2\\
\hline 35 & 1.7e-08 & 40.6\\
\hline 36 & 5.07e-09 & 36.4\\
\hline 37 & 6.38e-09 & 32\\
\hline 38 & 1.19e-09 & 26.9\\
\hline 39 & 6.36e-10 & 23.4\\
\hline 40 & 3.96e-10 & 19.8\\
\hline 41 & 8.61e-11 & 17.1\\
\hline 42 & 2.3e-11 & 13.7\\
\hline 43 & 1.24e-11 & 12.3\\
\hline 44 & 7.12e-12 & 10.5\\
\hline 45 & 5.87e-12 & 8.93\\
\hline 46 & 6.39e-13 & 5.8\\
\hline 47 & 2.77e-13 & 5.2\\
\hline 48 & 2.26e-13 & 4.34\\
\hline 49 & 4.12e-14 & 3.01\\
\hline 50 & 1.11e-14 & 2.07\\
\hline 51 & 7.22e-15 & 1.83\\
\hline 52 & 3.77e-15 & 1.51\\
\hline 53 & 3.66e-15 & 1.43\\
\hline \ldots & \ldots & \ldots \\
\hline 68 & 5.55e-16 & 0.000381\\
\hline 69 & 0 & 0\\

\hline
\end{tabular}
\label{tab:sweepConvergence2}
}
\end{table}

\pagebreak

\subsection{Performance on a different computer architecture}
%COMMENT REMOVED
\label{ss:octopus}

The performance/scaling of parallel methods is often strongly affected by hardware features of a particular shared memory
implementation. All parallel methods considered here scale better when the ratio of memory bandwidth to CPU speed is higher.
In addition, the scaling is affected by the network topology of the cores.  Stampede has ``dual, eight-core sockets," so communication between processors is necessarily slower when $P > 16$.

To explore the influence of these features, we repeat our main three examples on a different platform (``Octopus''): a computer with 8 Dual Core AMD Opteron 880 microprocessors running at 2.4 GHz, with 128 GB total RAM under the Scientific Linux v5.1 operating system.  We have implemented all methods in C++ and compiled with the
{\tt -O2} level of optimization using the g++ compiler v4.2.1.  %COMMENT REMOVED
The scaling was tested on up to 16 threads.
%COMMENT REMOVED
All other experimental settings are exactly the same as described for ``Stampede'' at the beginning of section \ref{s:experiments}.
The results are reported in Figure \ref{fig:octopus}.

While the main conclusions are the same as in subsections \ref{ss:serial_performance}-\ref{ss:parallel_performance}, this change in hardware architecture yields noticeably different relative performance even for serial methods.  We observe that FMM seems to benefit more from larger cache sizes than FSM and LSM do; thus, on Octopus the sweeping methods appear more competitive on large grids than in the previous tests on Stampede.  The HCM2, whose algorithmic behavior is similar to FMM, is also less advantageous on Octopus, while HCM16 and HCM32 (whose computational cost is dominated by cell-sweeping) appear to be more advantageous here for large grids.

As for scaling (Figures \ref{fig:octopus}D - \ref{fig:octopus}F), all parallel methods seem to do much better on Octopus than on Stampede, even when only the first 16 threads are accounted for on Stampede.  For example, on Octopus the pHCM8 median scaling curve has approximate slopes of .6, .92, and .83 on the three examples, while on Stampede the slopes up to $P = 16$ are approximately .5, .8, and .73.  For pHCM4 on Octopus, the slopes are approximately .33, .73, and .67 (making pHCM4 very competitive on Octopus), while on Stampede the slopes up to $P = 16$ are only .27, .43, and .43.
The scaling for DFSM not only improves on Octopus, but
the slope of the scaling curve appears to be higher when the number of threads exceeds 8.
%COMMENT REMOVED
%COMMENT REMOVED

\noindent
\begin{figure}[H]
\hspace{-.8 in}
$
\begin{array}{ccc}

\includegraphics[scale = .45]{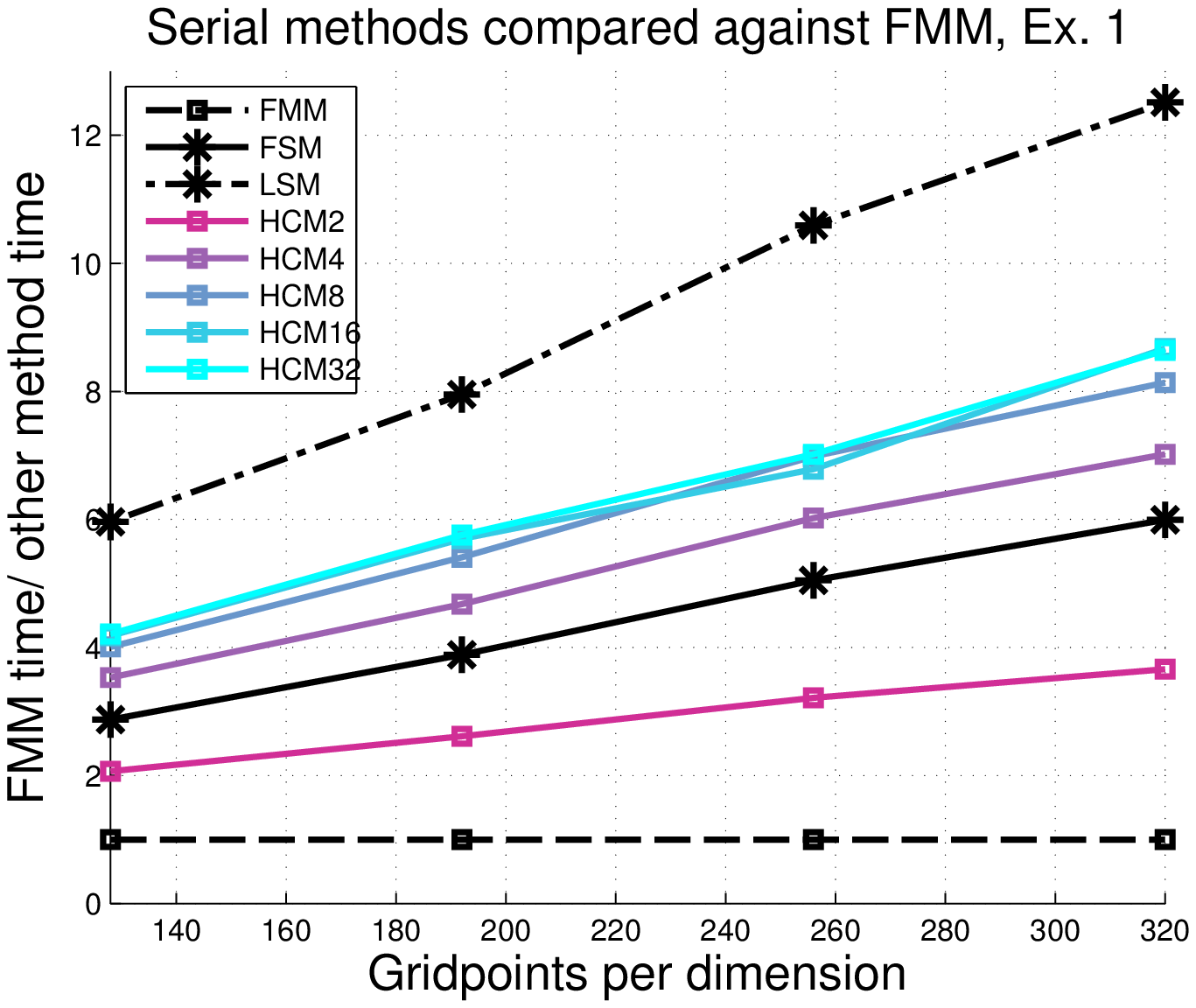}&
\hspace{-.2in}
\includegraphics[scale = .45]{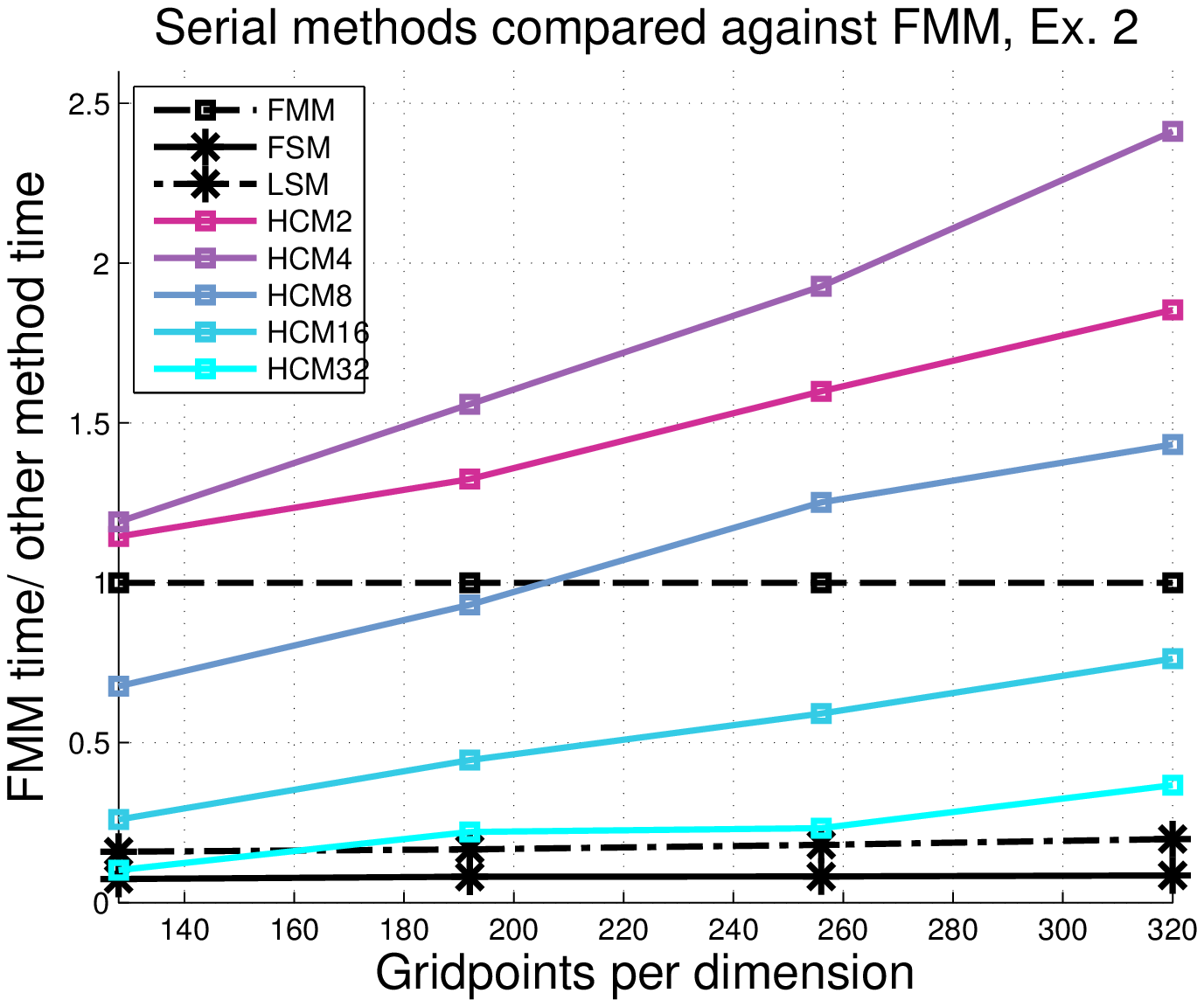}&
\hspace{-.2in}
\includegraphics[scale = .45]{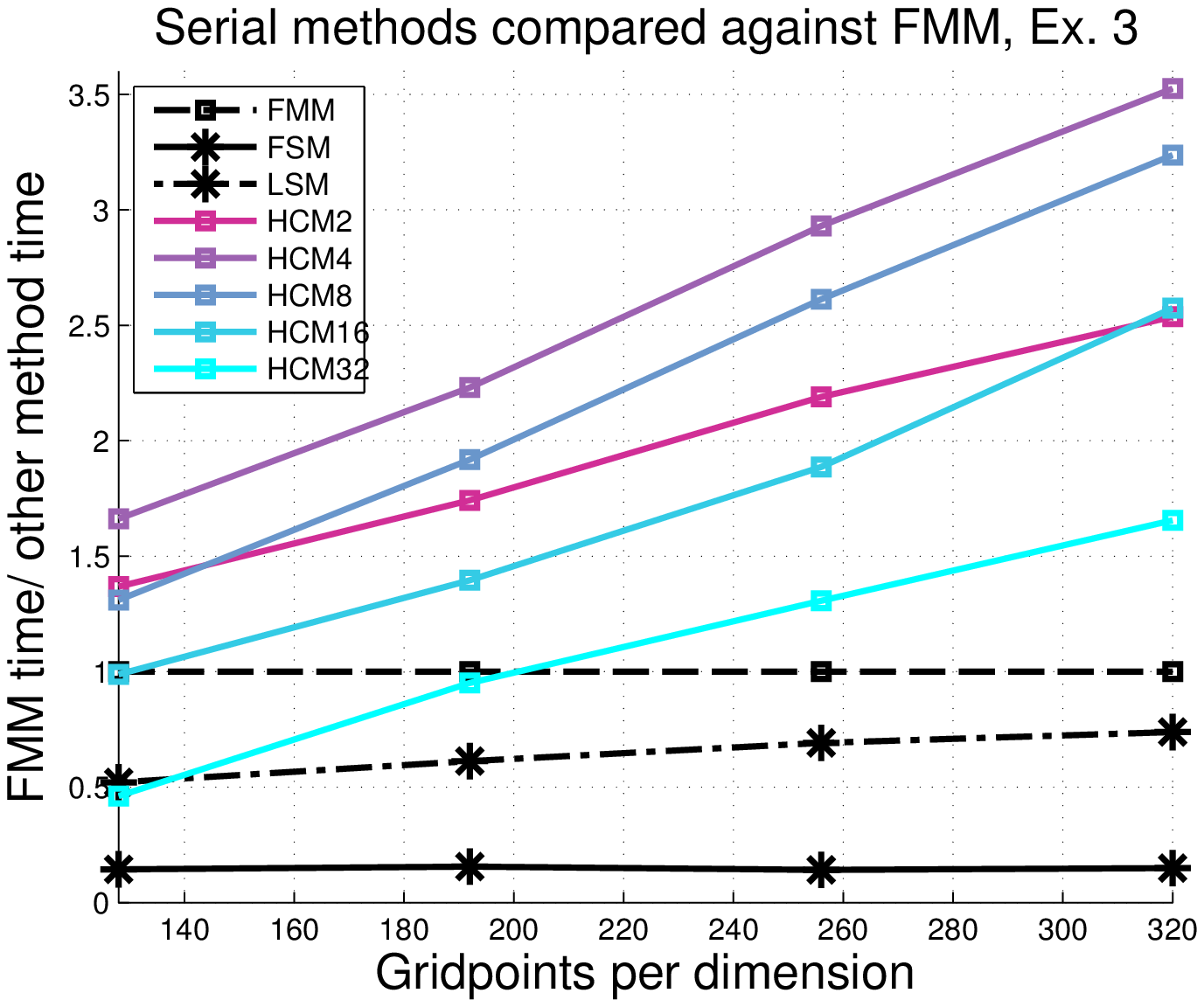}\\

A&B&C\\

\includegraphics[scale = .45]{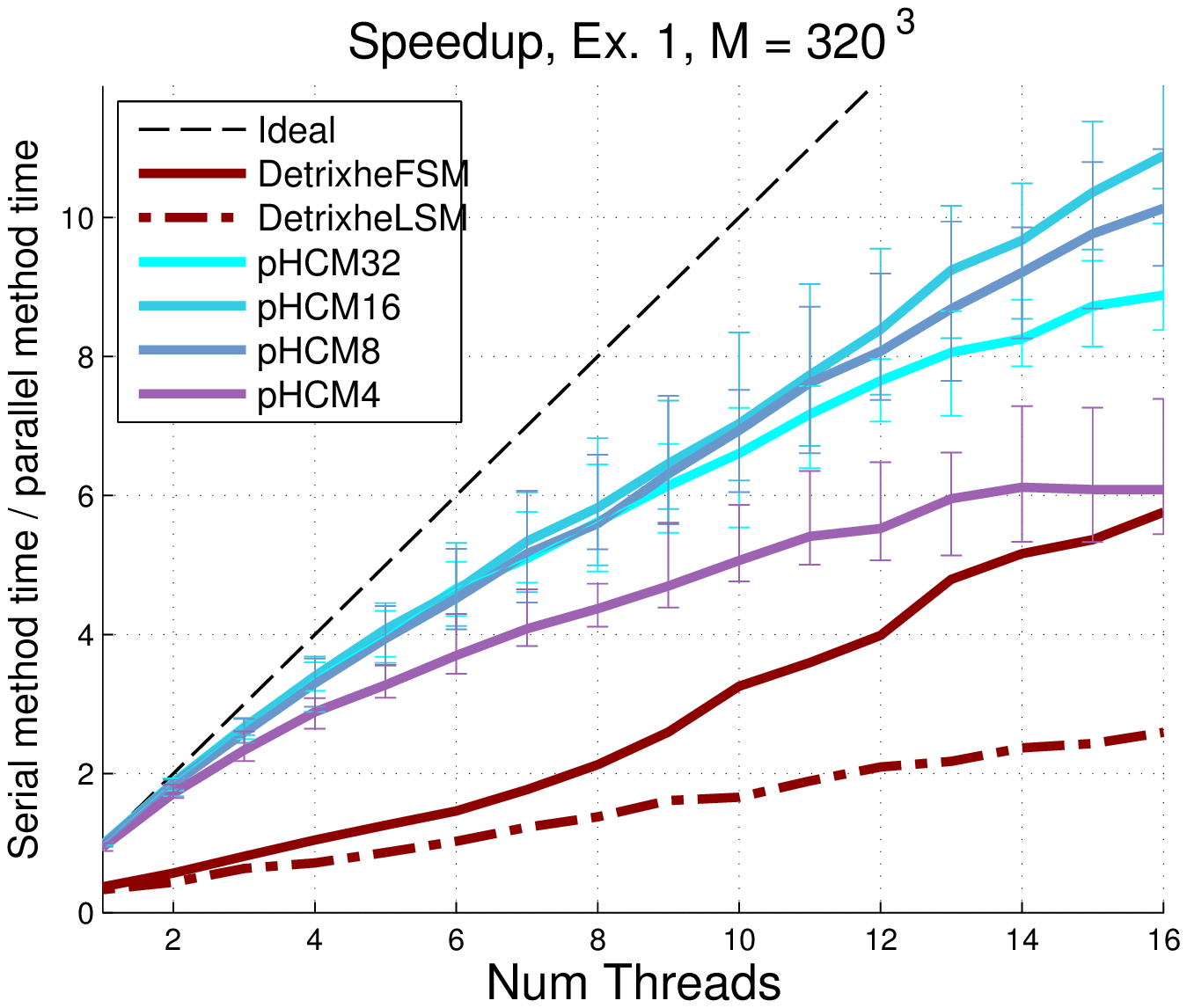}&
\hspace{-.2in}
\includegraphics[scale = .45]{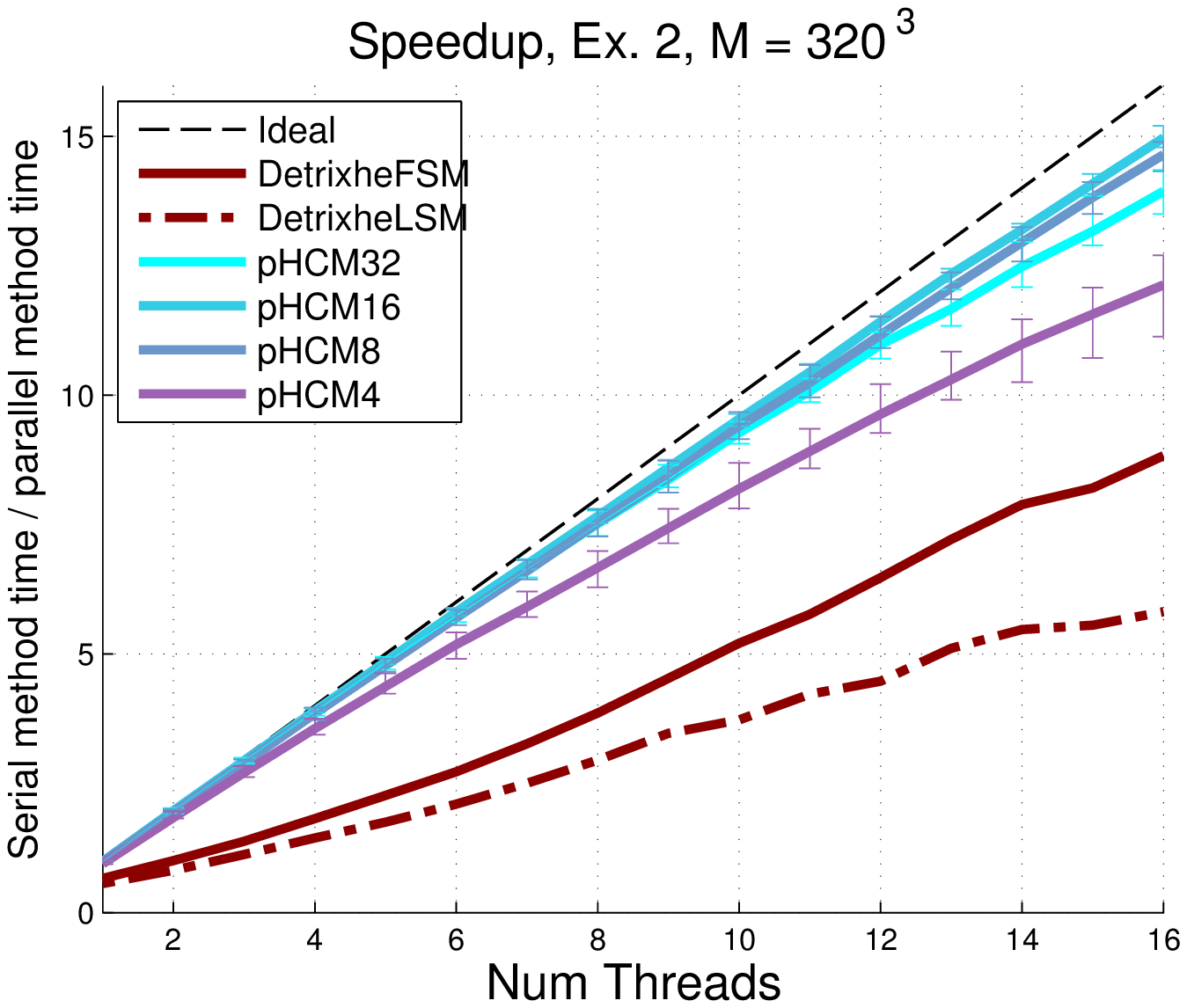}&
\hspace{-.2in}
\includegraphics[scale = .45]{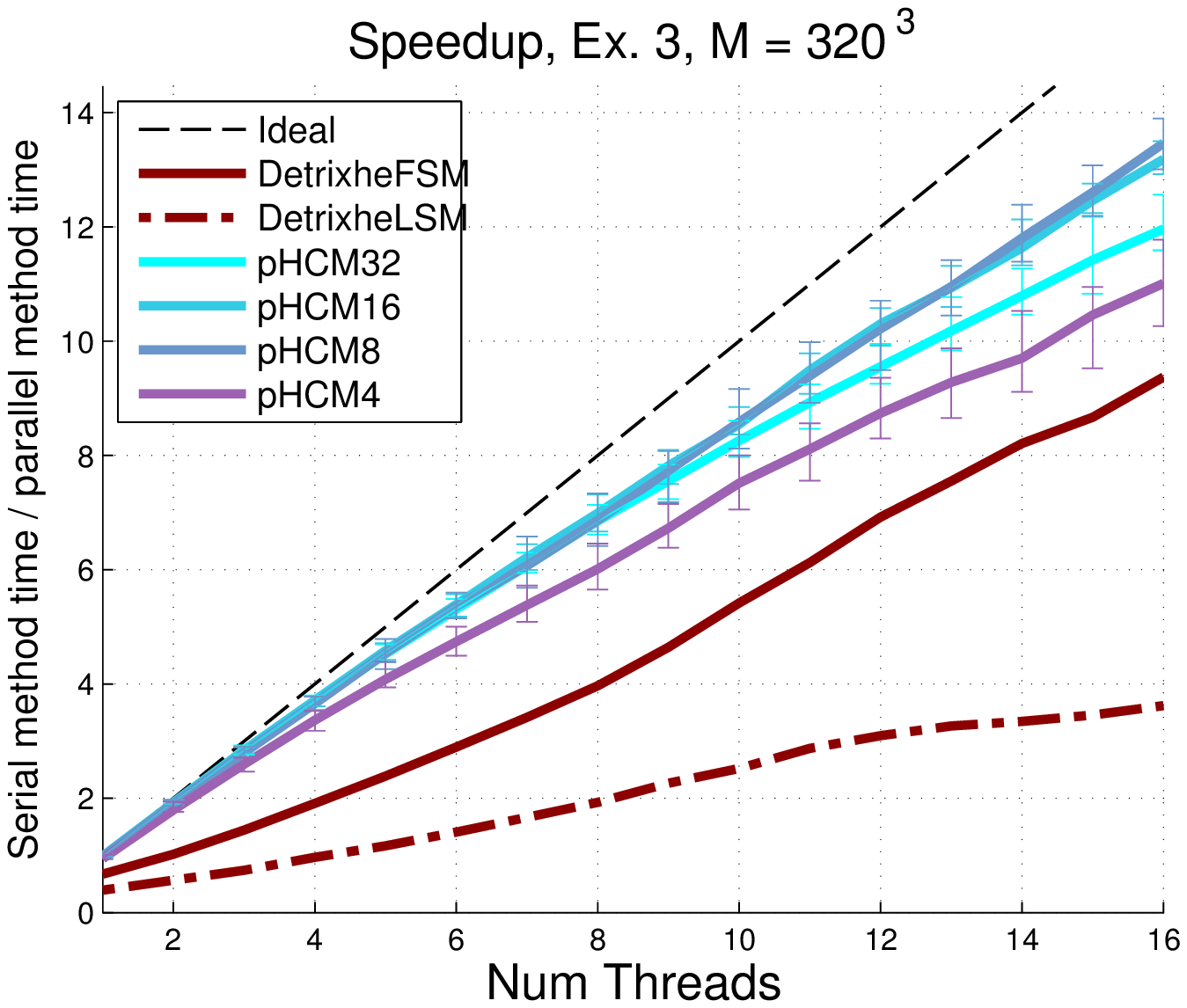}\\

D&E&F\\

\includegraphics[scale = .45]{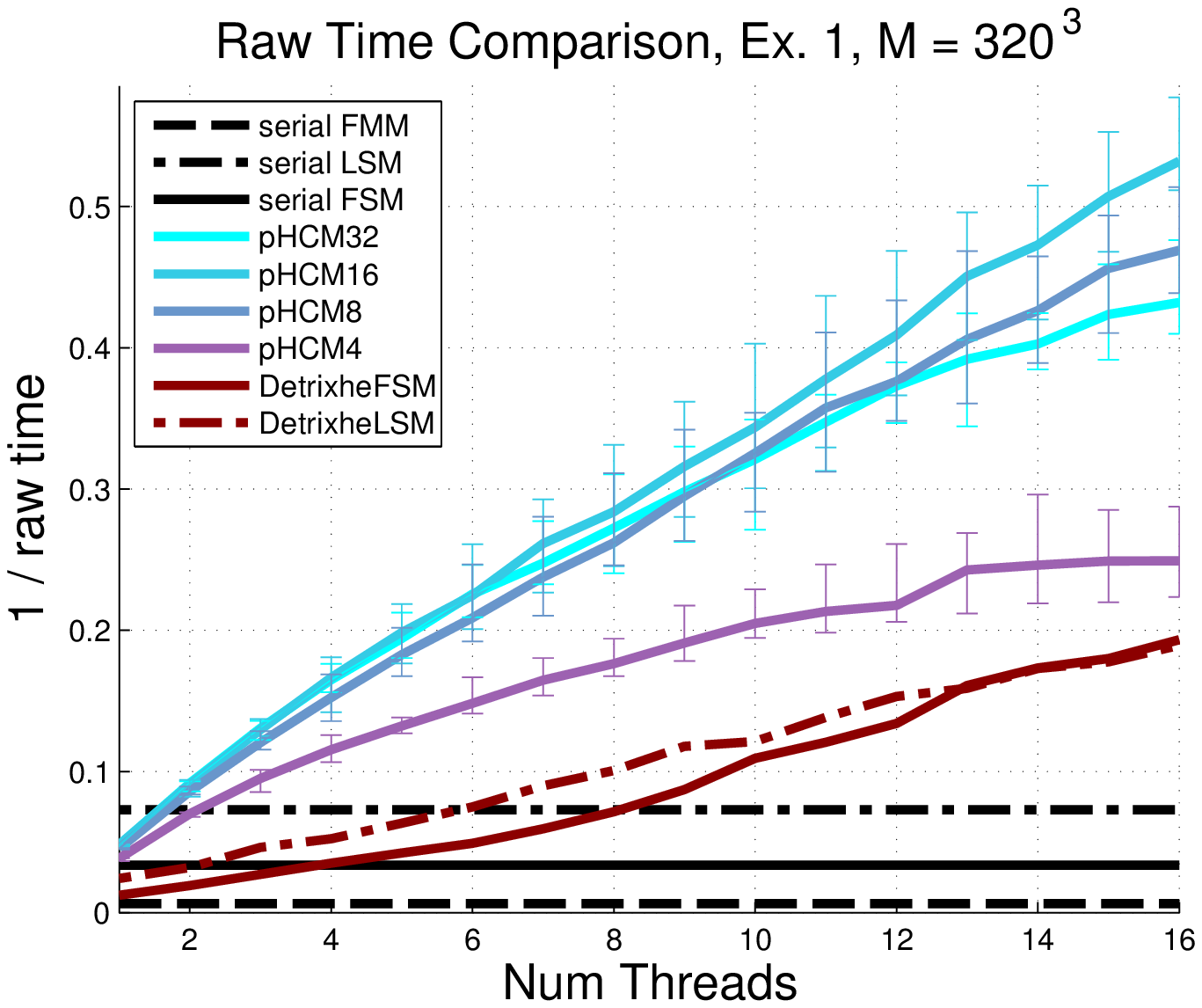}&
\hspace{-.2in}
\includegraphics[scale = .45]{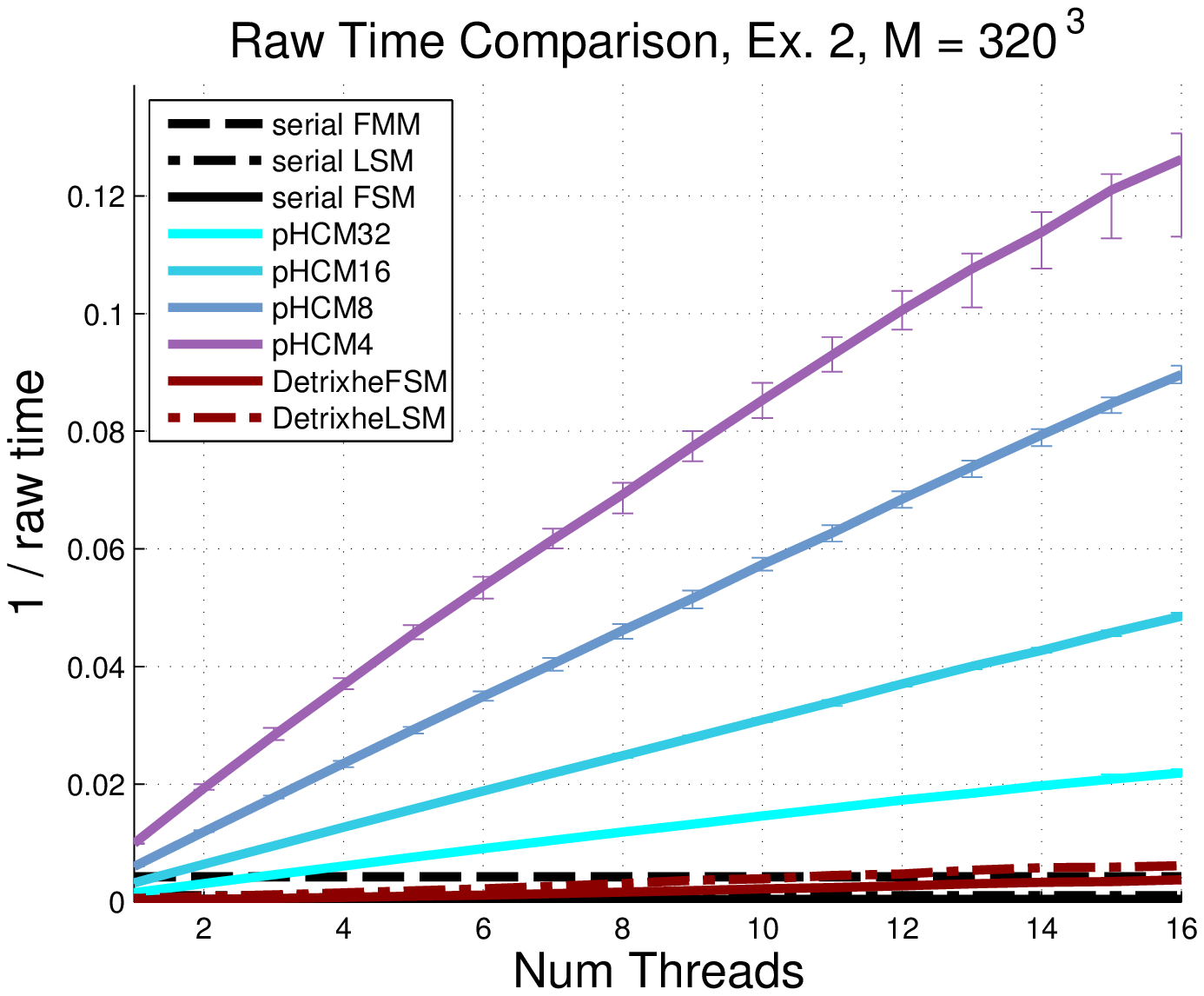}&
\hspace{-.2in}
\includegraphics[scale = .45]{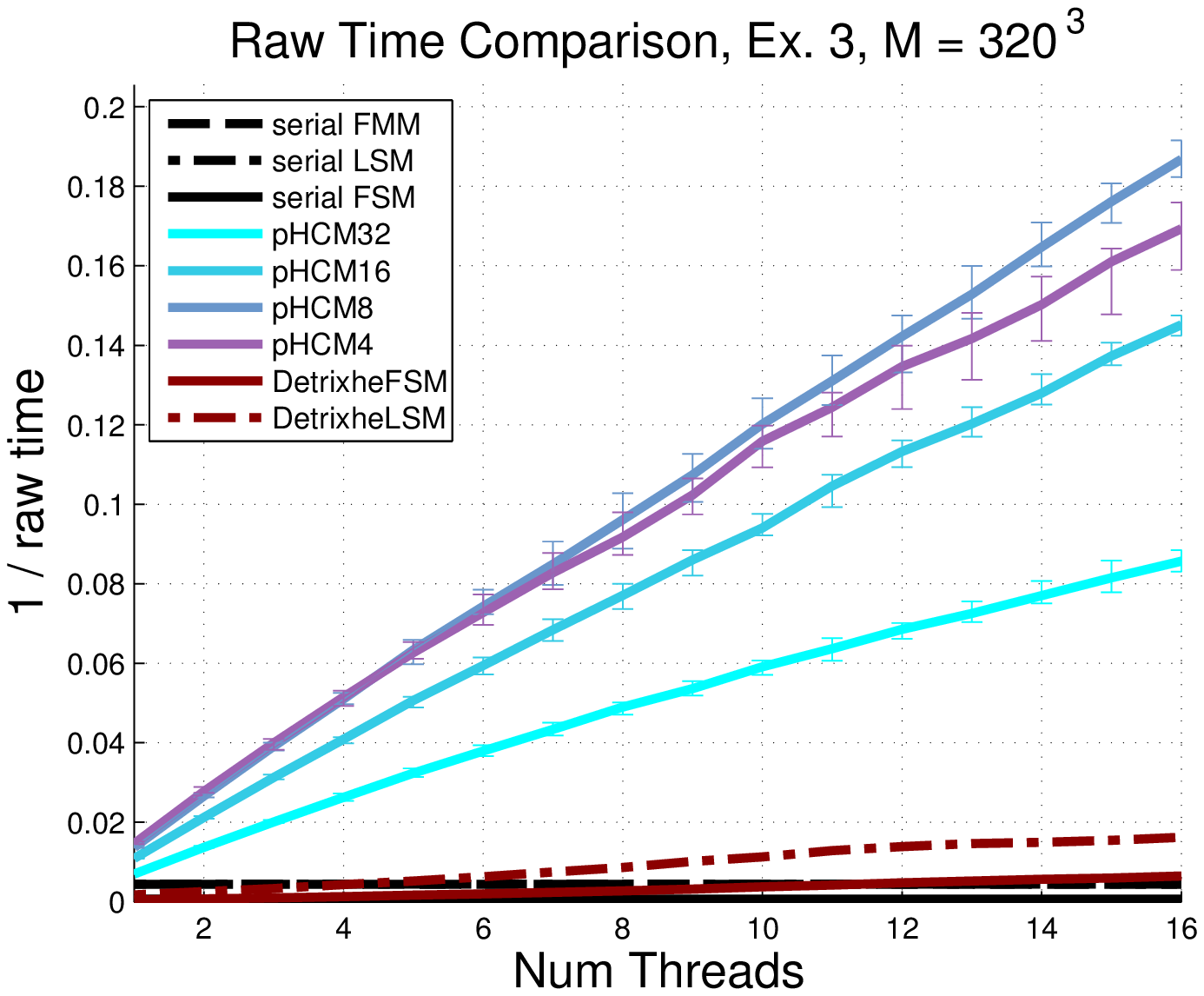}\\

G&H&I

\end{array}
$

\caption{\footnotesize ``Octopus" testing (subsection \ref{ss:octopus}).
Top row: performance of serial methods for different $M$; compare with Figure \ref{fig:serial_charts}.
Two bottom rows: scaling and performance for pHCM at $M = 320^3$; compare with Figure \ref{fig:parallel_charts}.
}
\label{fig:octopus}
\end{figure}

\pagebreak

\subsection{Additional examples: checkerboard speed functions}
\label{ss:checkerboard}

We consider two additional examples with periodic piecewise constant speed functions, which generalize the 2D checkerboard test problems of \cite{ChacVlad, ChacVlad_expanded}.
These examples arise in the numerical computation of effective Hamiltonians in highly oscillatory problems; see also \cite{ObTaVlad}.

\begin{figure}[H]
\hspace{-.6 in}
$
\begin{array}{ccc}
\includegraphics[scale = .45]{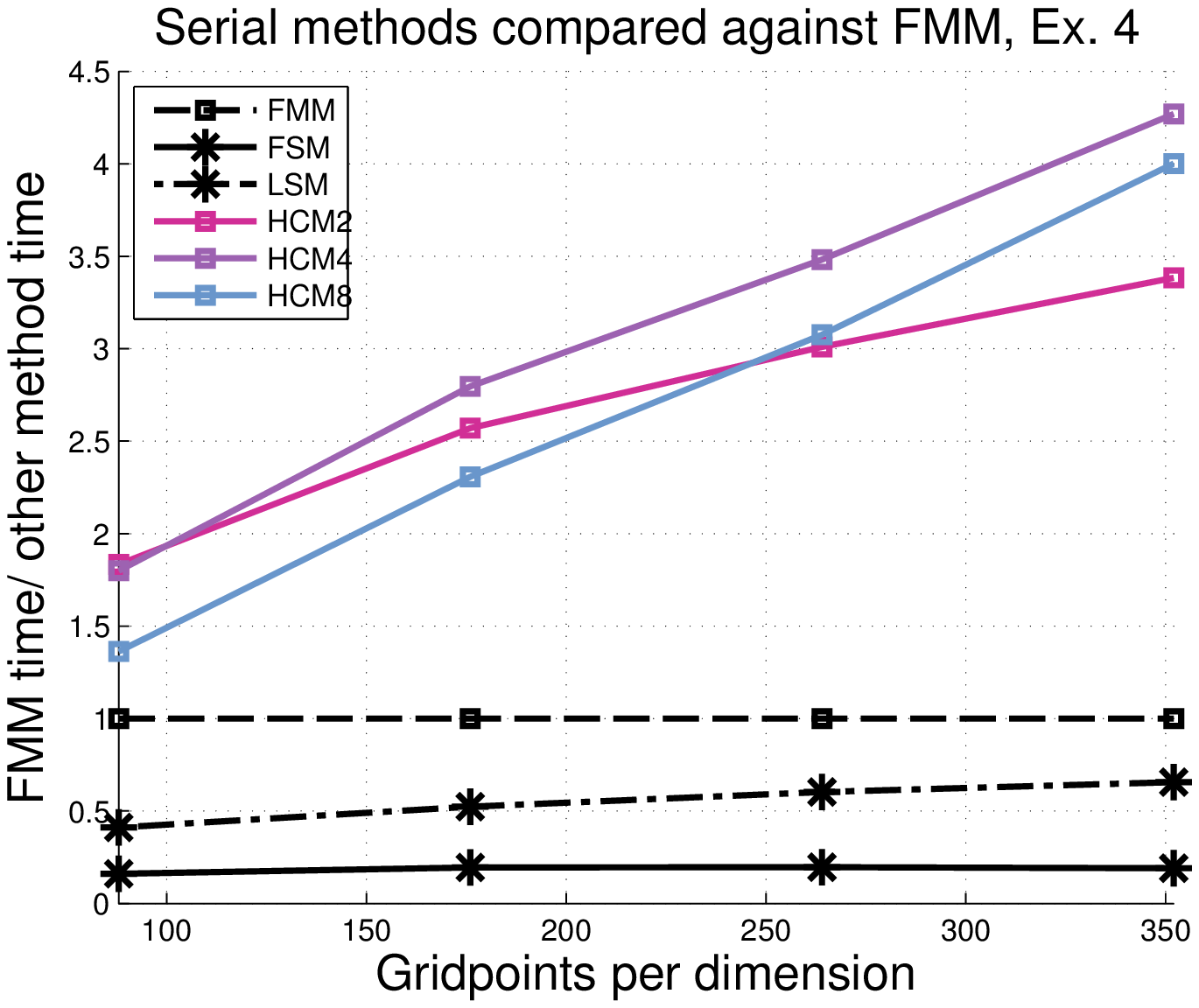}&
\hspace{-.2in}
\includegraphics[scale = .45]{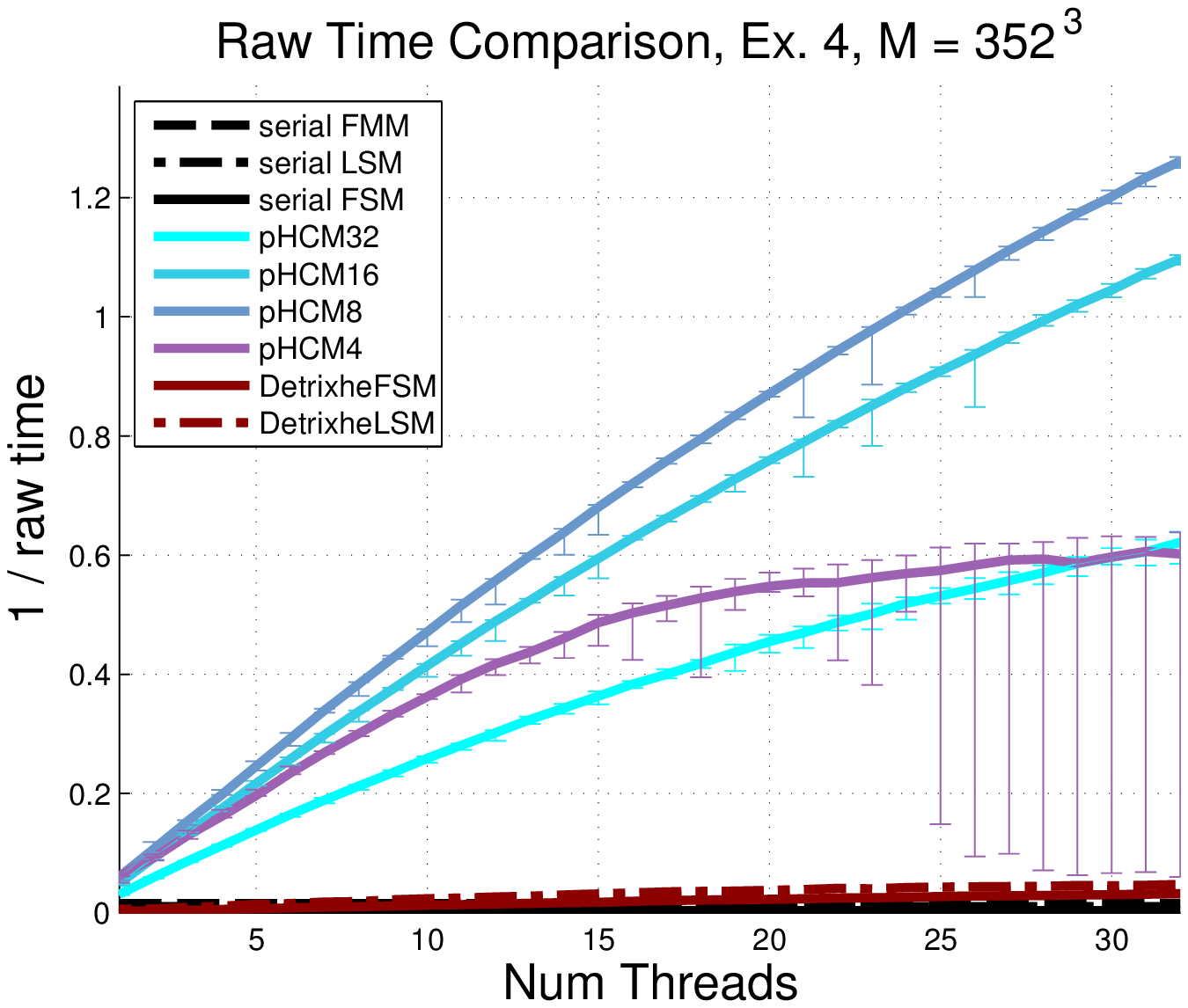}&
\hspace{-.2in}
\includegraphics[scale = .45]{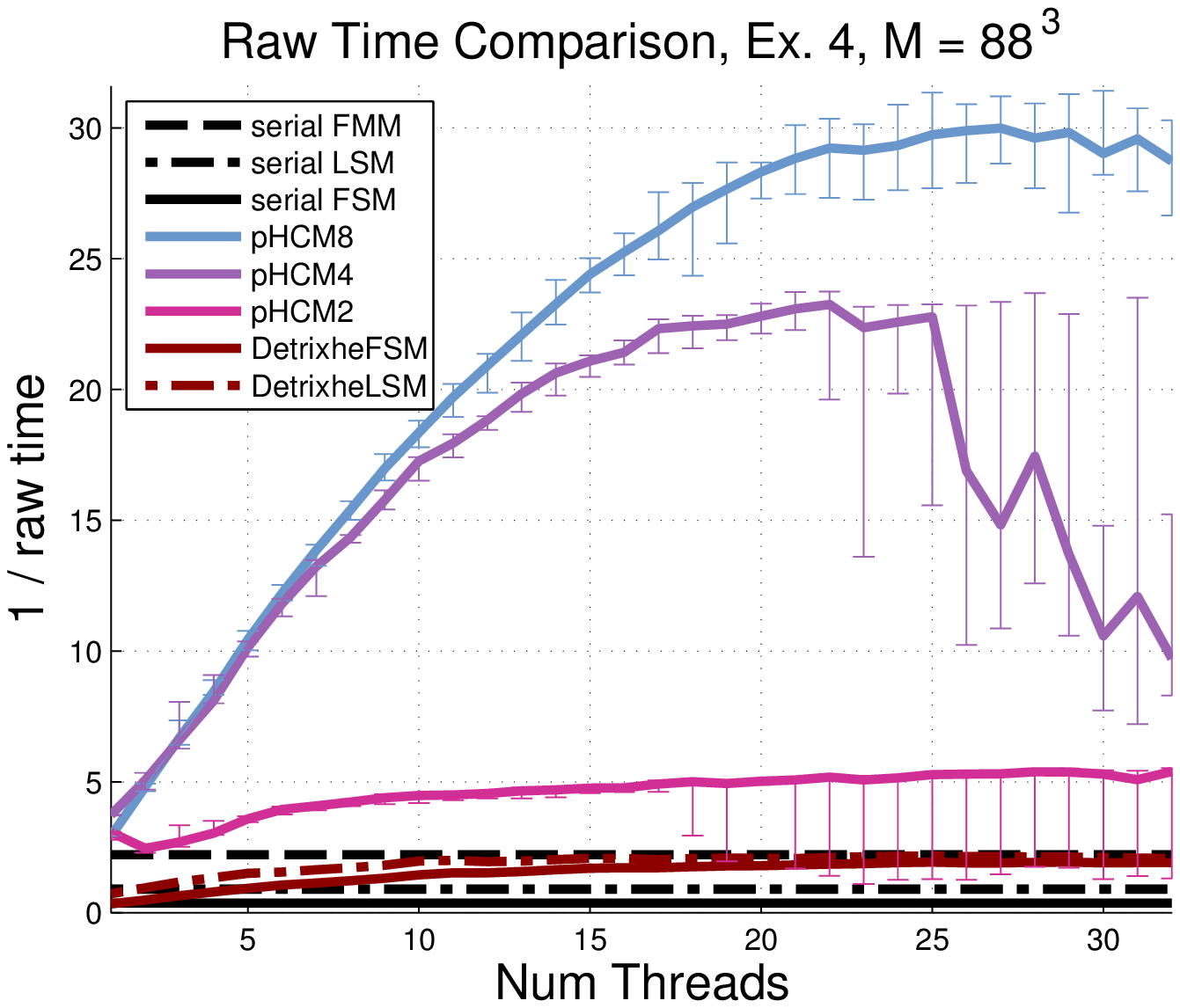}\\

A&B&C\\

\includegraphics[scale = .45]{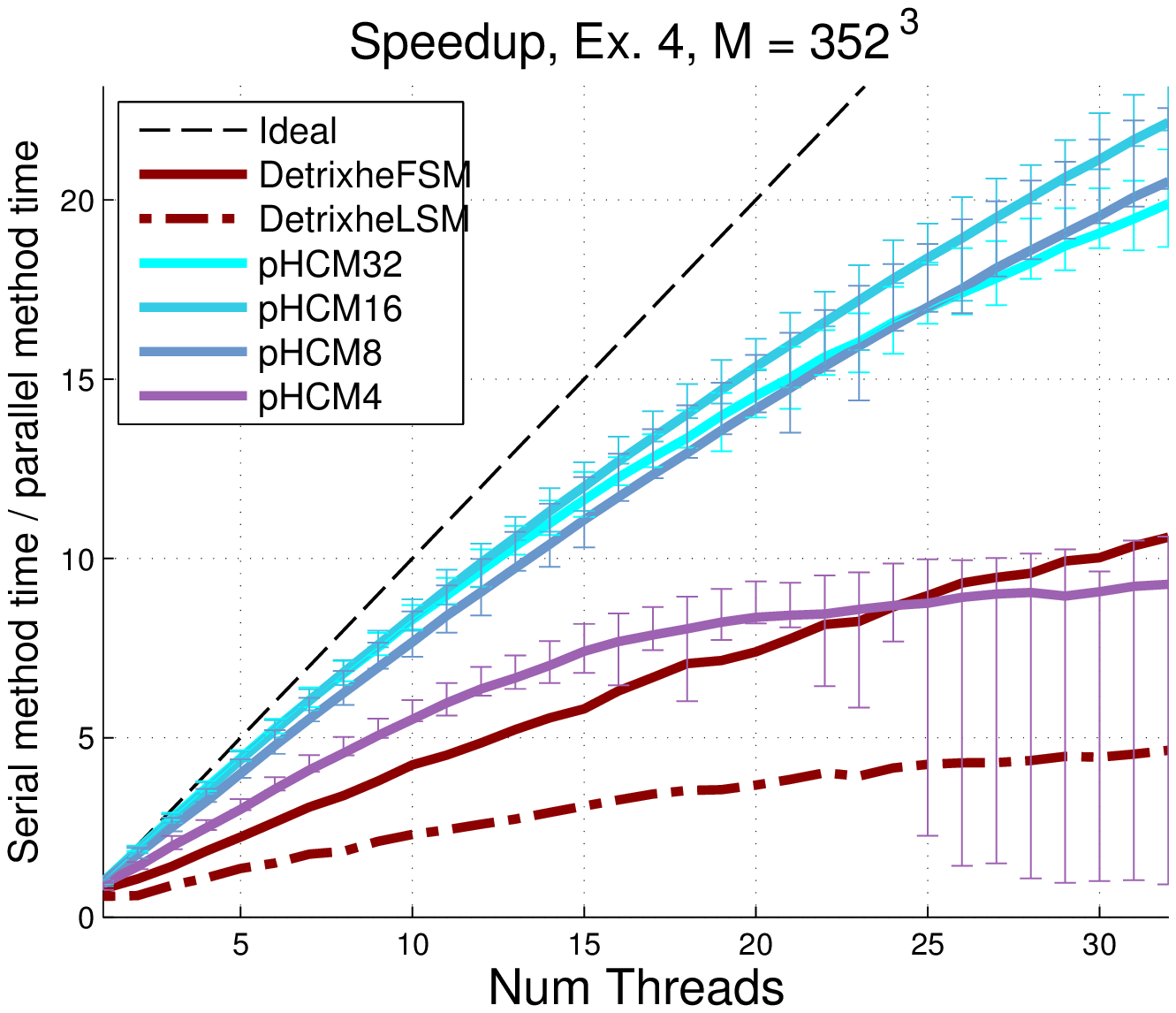}&
\hspace{-.2in}
\includegraphics[scale = .45]{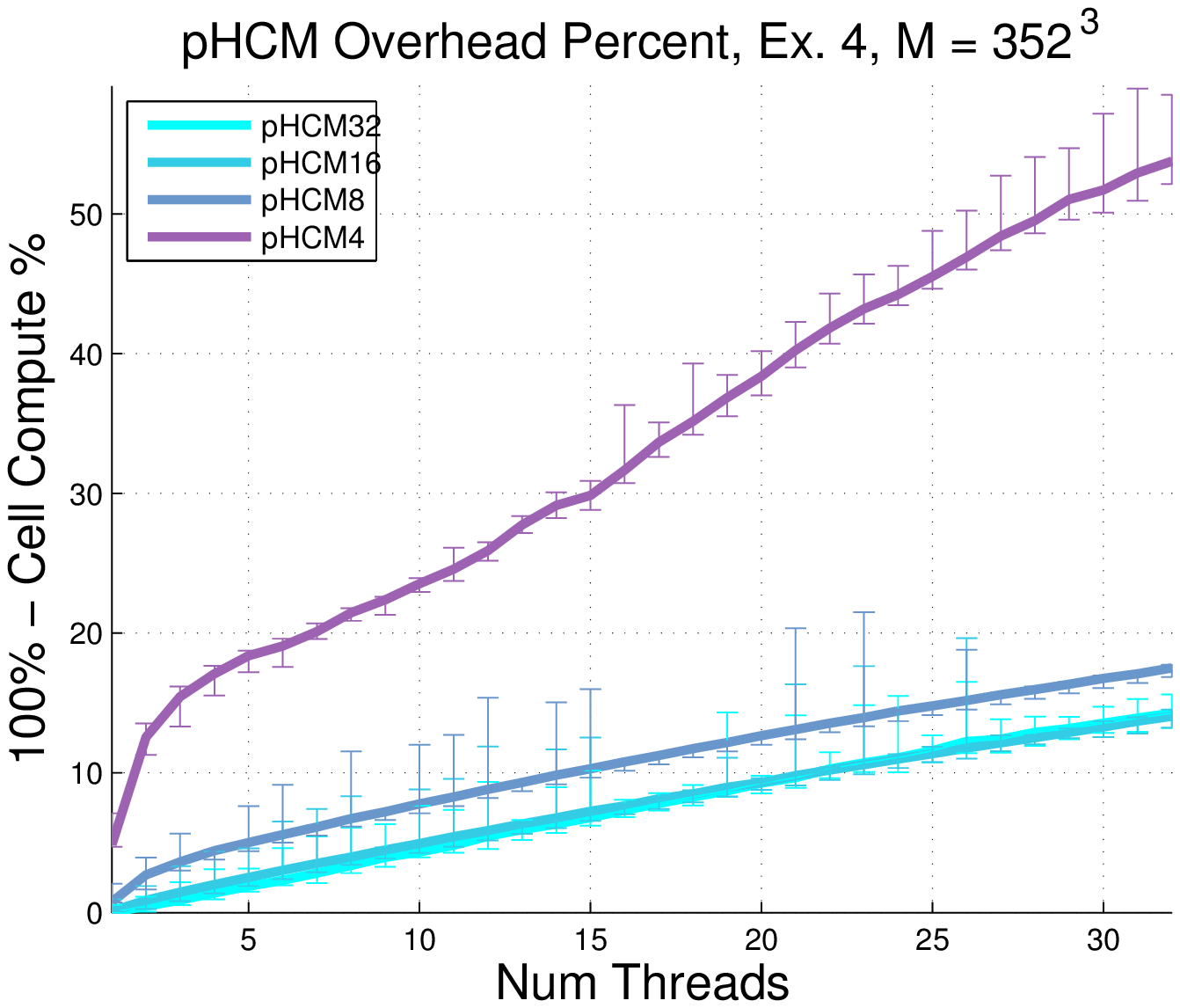}&
\hspace{-.2in}
\includegraphics[scale = .45]{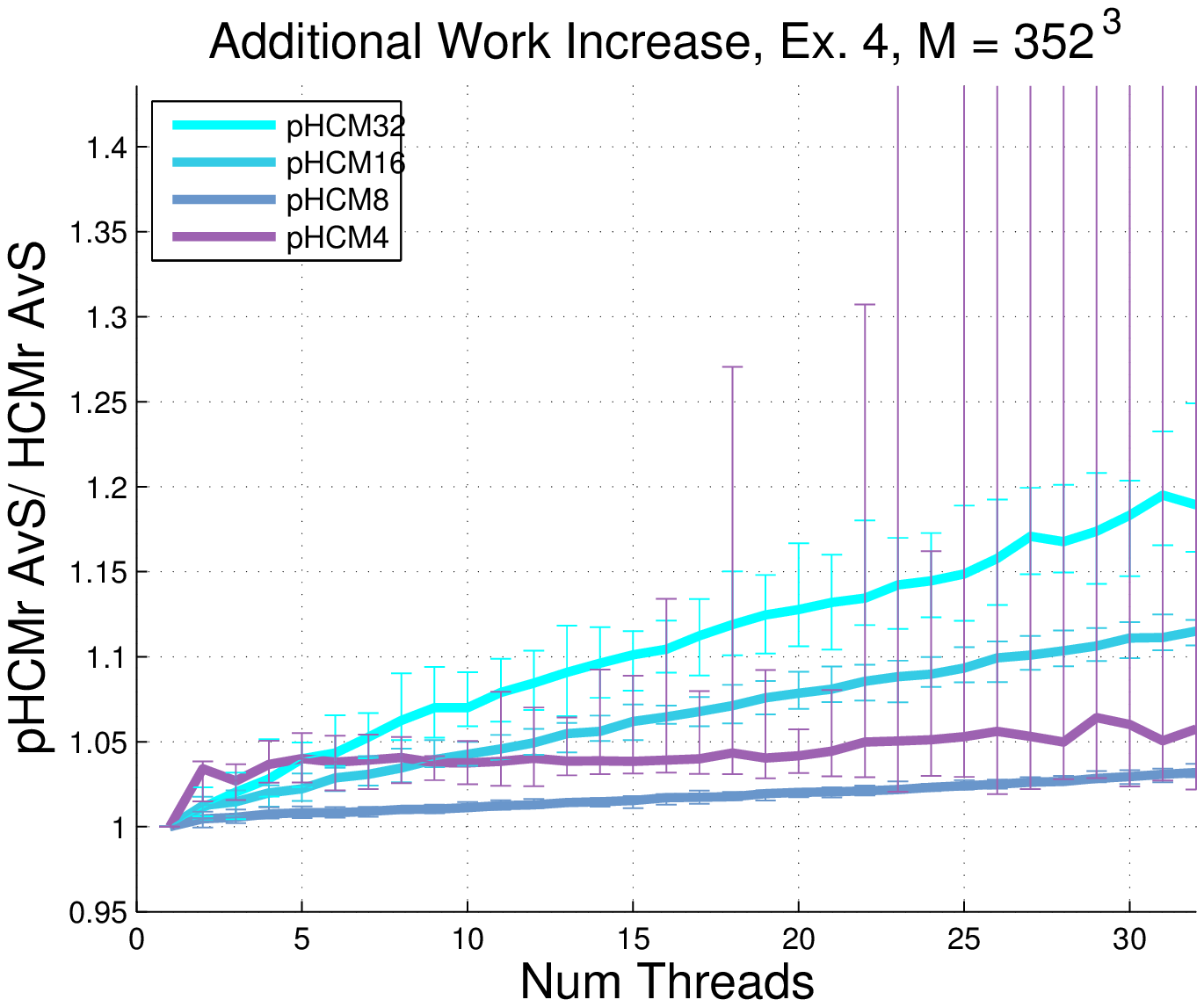}\\

D&E&F\\

\includegraphics[scale = .45]{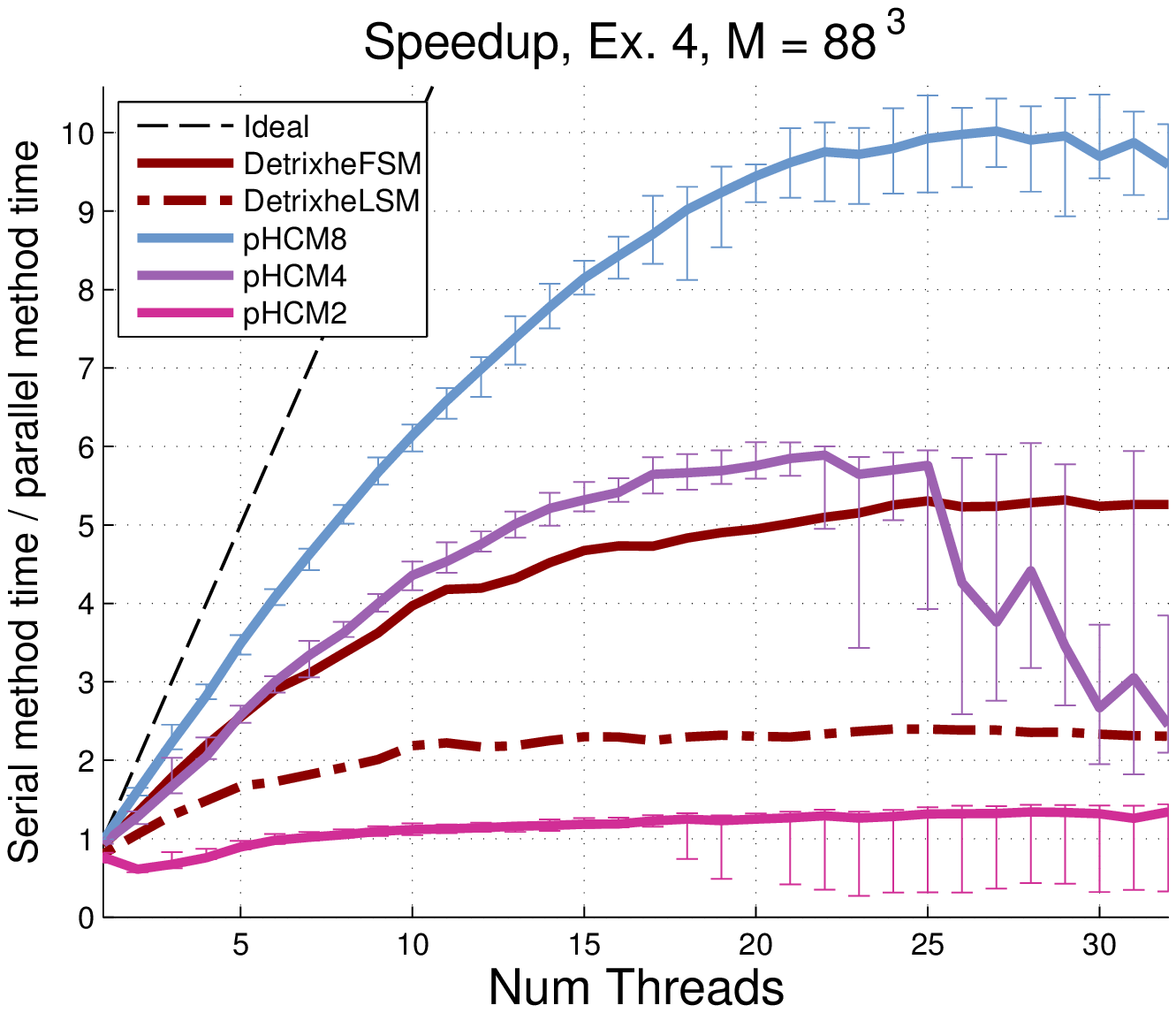}&
\hspace{-.2in}
\includegraphics[scale = .45]{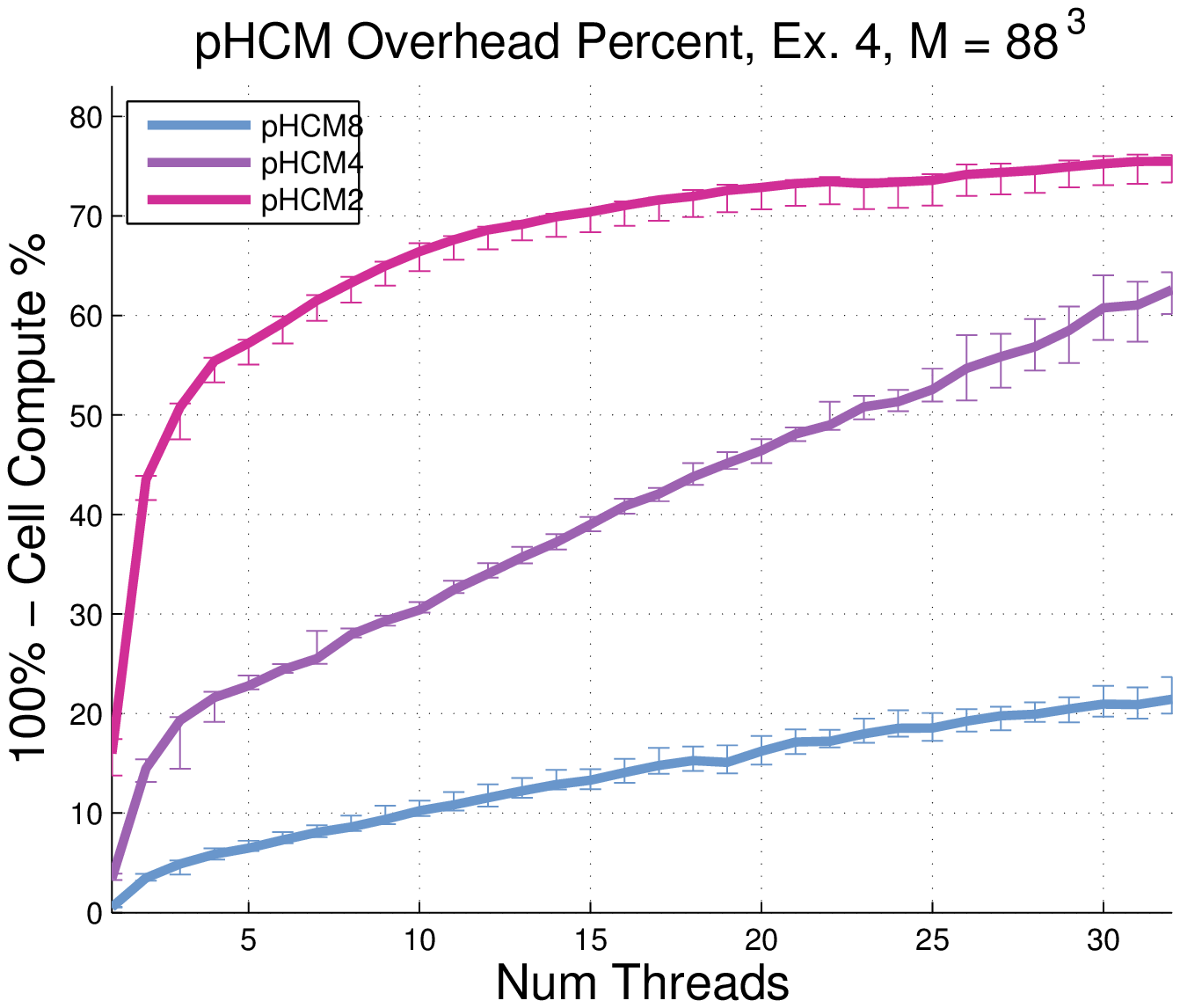}&
\hspace{-.2in}
\includegraphics[scale = .45]{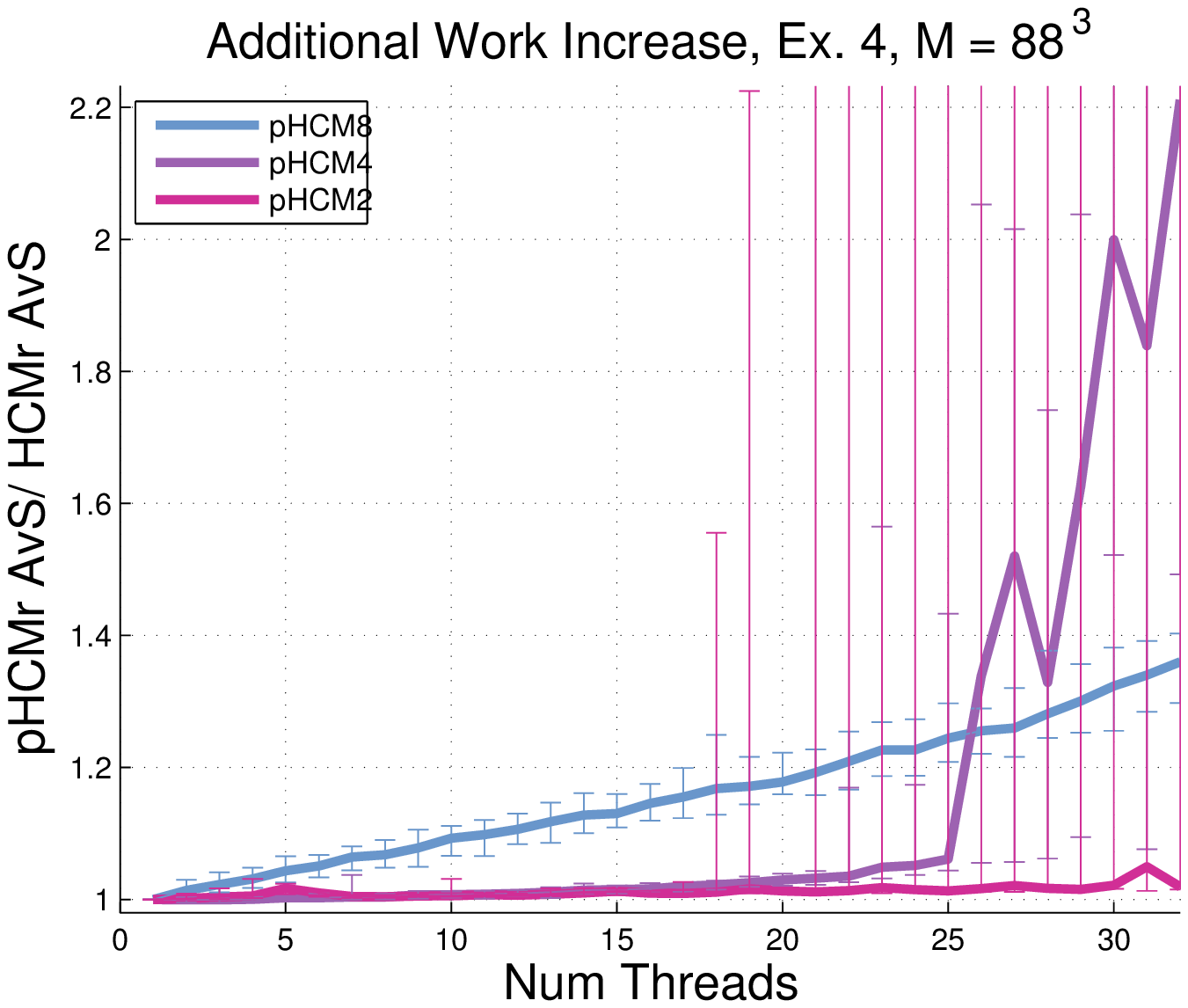}\\

G&H&I\\
\end{array}
$
\caption{\footnotesize 3D Checkerboard example with $K=11$ (subsection \ref{ss:checkerboard}).  Chart $A$ is a comparison of serial methods for different $M$; compare with Figure \ref{fig:serial_charts}.  Scaling/performance for parallel methods with $M=352^3$ is shown in charts $B$ and $D$; compare with Figure \ref{fig:parallel_charts}.  Parallel overhead and additional work with $M=352^3$ are shown in charts $E$ and $F$; compare with Figure \ref{fig:overheads_320}.  The same information for a coarser grid with $M=88^3$ is shown in charts $C$ and $G-I$.
}
\label{fig:checkerboard_11_cubed}
\end{figure}

\pagebreak
Suppose that the unit cube $\cdomain$ is split into $K^3$ smaller cubes (or ``3D checkers")
of edge length $1/K$.  Suppose these smaller cubes are divided into two types (``black" and ``white'') so that no two cubes of the same type have a face in common.  The speed function $F$ is defined to be 2 on black cubes and 1 on white cubes\footnote{We can also take $F=2$ on the boundary of the cubes.  Computationally, the issue does not arise since our gridsizes are selected to ensure that each gridpoint is in the interior of either black or white cube.}.  The exit set $Q$ again consists of a single point in the center of $\cdomain$ and, given the even number of gridpoints, the set $Q'$ consists of 8 gridpoints.

We conduct experiments on 2 different 3D checkerboards: with $K=11$ and $K=41$; the respective performance/scaling results are summarized in Figures \ref{fig:checkerboard_11_cubed} and \ref{fig:checkerboard_41_cubed}.  As observed in \cite{ChacVlad}, HCM performs very well on problems where the discontinuities of the speed function align with cell boundaries.
The scaling trends for $K = 11$ are most similar to those observed in Example 2, where the speed function is also highly oscillatory.  For $K = 41$, the speedup for pHCM4 is surprisingly large and stable.

\begin{figure}[H]
\hspace{-.6 in}
$
\begin{array}{ccc}
\includegraphics[scale = .45]{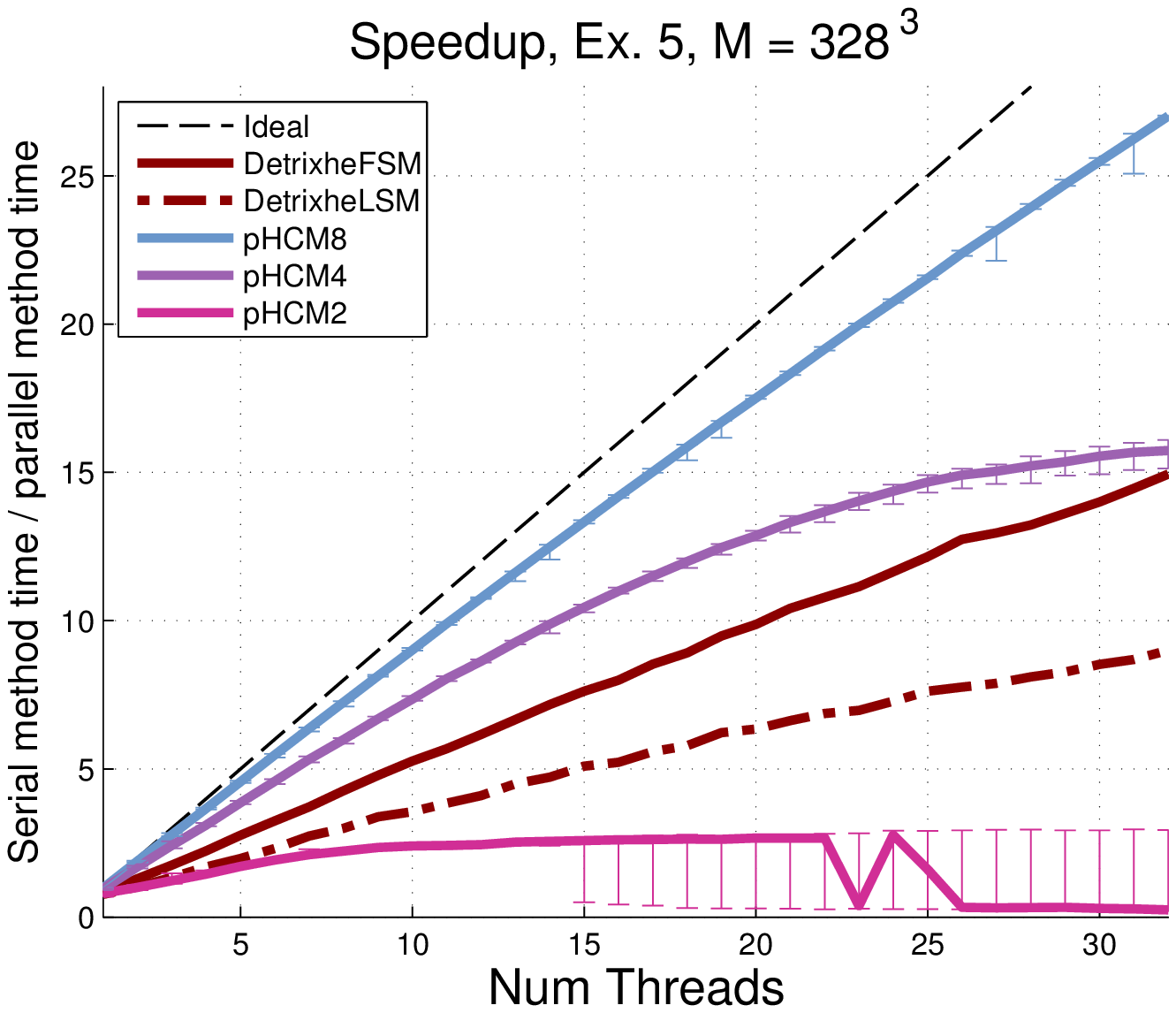}&
\hspace{-.2in}
\includegraphics[scale = .45]{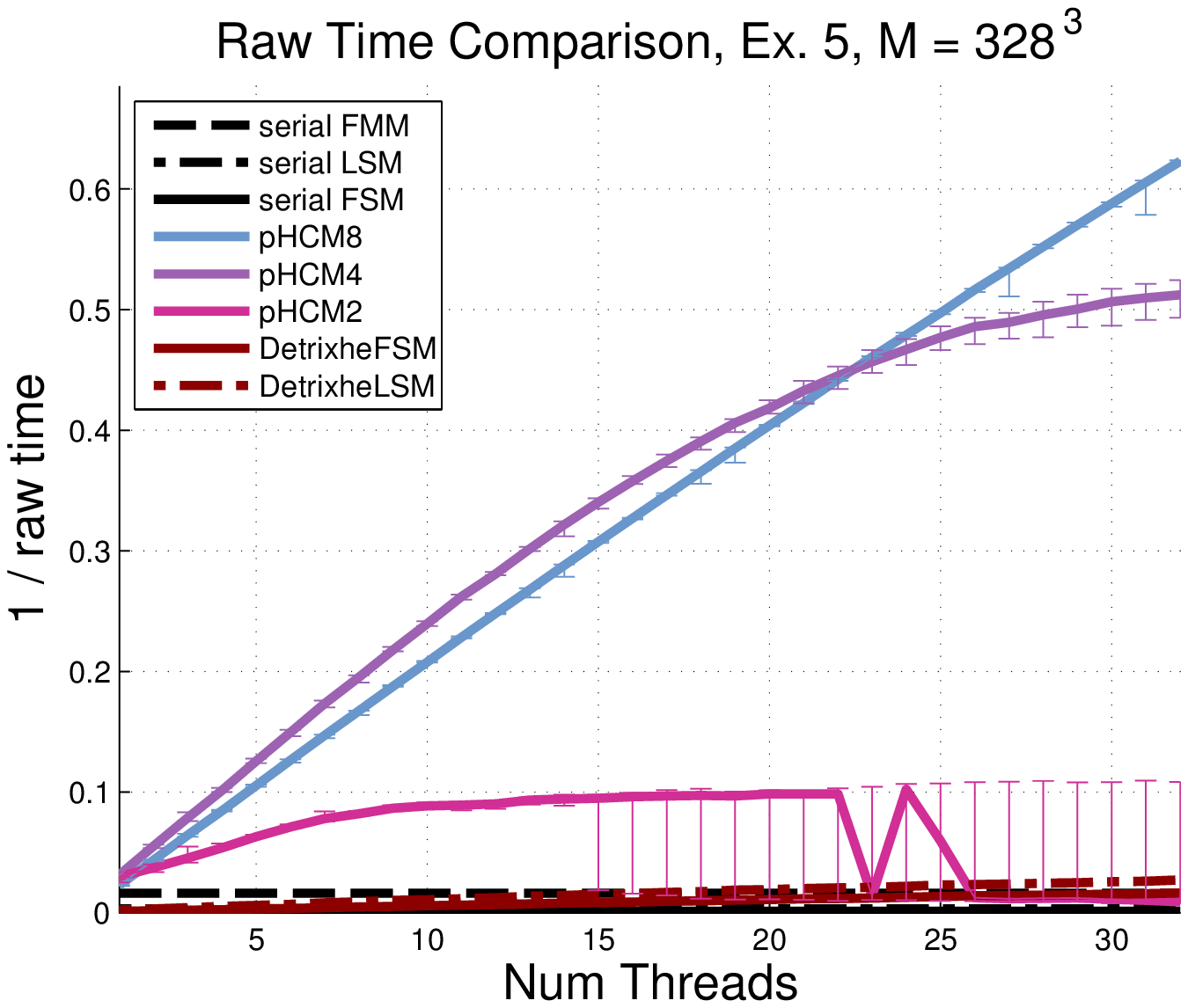}&
\hspace{-.2in}
\includegraphics[scale = .45]{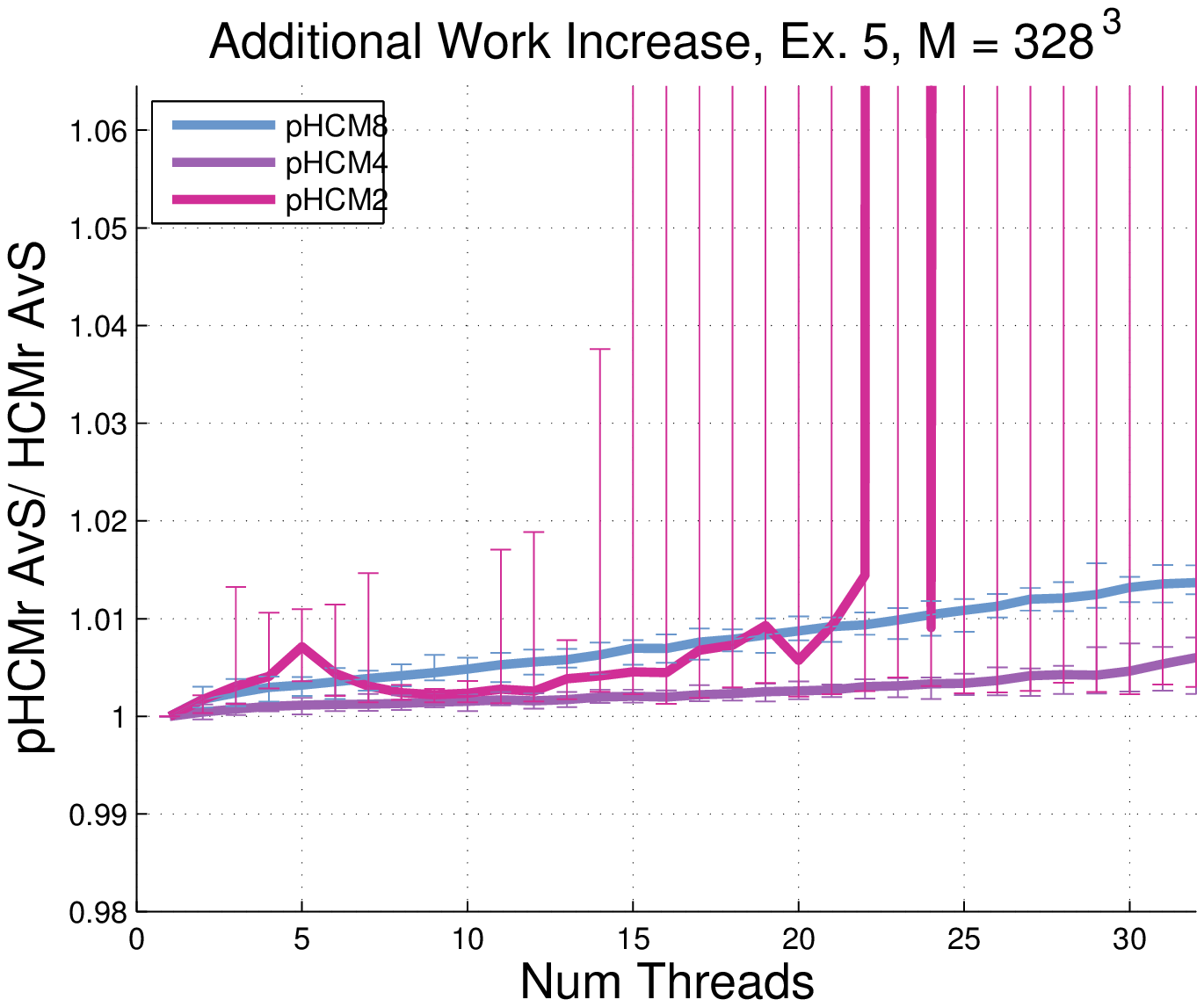}\\

A&B&C\\
\end{array}
$
\caption{\footnotesize 3D Checkerboard example with $K=41$ (subsection \ref{ss:checkerboard}).
}
\label{fig:checkerboard_41_cubed}
\end{figure}

\pagebreak

\subsection{Additional examples: maze speed functions}
\label{ss:shellMaze}

Suppose the domain contains four concentric spherical ``barriers" of thickness $t$ that have openings on alternating sides.  Specifically, $\cdomain = [-1,1]^3$, $Q = \{0,0,0\}$, and $F = 1$ outside the set of (slowly permeable) barriers and .001 inside, with the barriers described as follows:

\begin{align*}
A_1 &= \{\x | .3 < |\x| < .3 + t \} \backslash \left( \{x^2 + y^2 < w\} \cap \{z < 0\}\right) \\
A_2 &= \{\x | .5 < |\x| < .5 + t \} \backslash \left( \{x^2 + y^2 < w\} \cap \{z > 0\}\right) \\
A_3 &= \{\x | .7 < |\x| < .7 + t \} \backslash \left( \{x^2 + y^2 < w\} \cap \{z < 0\}\right) \\
A_4 &= \{\x | .9 < |\x| < .9 + t \} \backslash \left( \{x^2 + y^2 < w\} \cap \{z > 0\}\right) \\
\end{align*}

\noindent where $t = 1/12$ and $w = 1/10$.
This is a modified version of an example from \cite{Detrixhe}, where the barriers considered were impermeable (i.e., with $F=0$).
Unlike the checkerboard examples, here the discontinuities of the speed function {\bf do not} align with the cell boundaries in any special way.  In that sense, this problem is also analogous to the second ``comb maze" example from our previous work; see Section 4.1 in \cite{ChacVlad}.

First, Figure \ref{fig:shellMaze}$A$ shows HCM$r$ is very effective for each $r$.
One of the drawbacks of the original version of HCM \cite{ChacVlad} was precisely the slow convergence on problems of this type.
%COMMENT REMOVED
The greatly improved performance shown here is due to the use of the new cell value heuristic (equation \eqref{eq:cell_value}).

The pHCM's speedup (Fig. \ref{fig:shellMaze} $B$), on the other hand, is significantly lower here %COMMENT REMOVED
%COMMENT REMOVED
(while for DFSM the speedup here is still typical).  We believe this is due to certain level sets of the value function getting ``pinched" at the locations where there is a hole in one of the barriers.  %COMMENT REMOVED
If the ordering of non-barrier cells is strictly causal, this means that, at several stages of the algorithm, there is only one cell upon which all still-to-be-computed cells depend.
(For example, since $w = .1$, in pHCM16 at most one cell will fit through the hole in each barrier.)
Furthermore, as mentioned in section \ref{ss:parallel_performance}, pHCM sees an increase in work over HCM for problems with a strictly causal cell ordering.  However, due to the large-enough advantage that HCM holds over other serial methods, the performance of pHCM is still significantly better than that of DFSM/DLSM; see Fig. \ref{fig:shellMaze} $C$.

%COMMENT REMOVED

\begin{figure}[H]
\hspace{-.6 in}
$
\begin{array}{ccc}
\includegraphics[scale = .45]{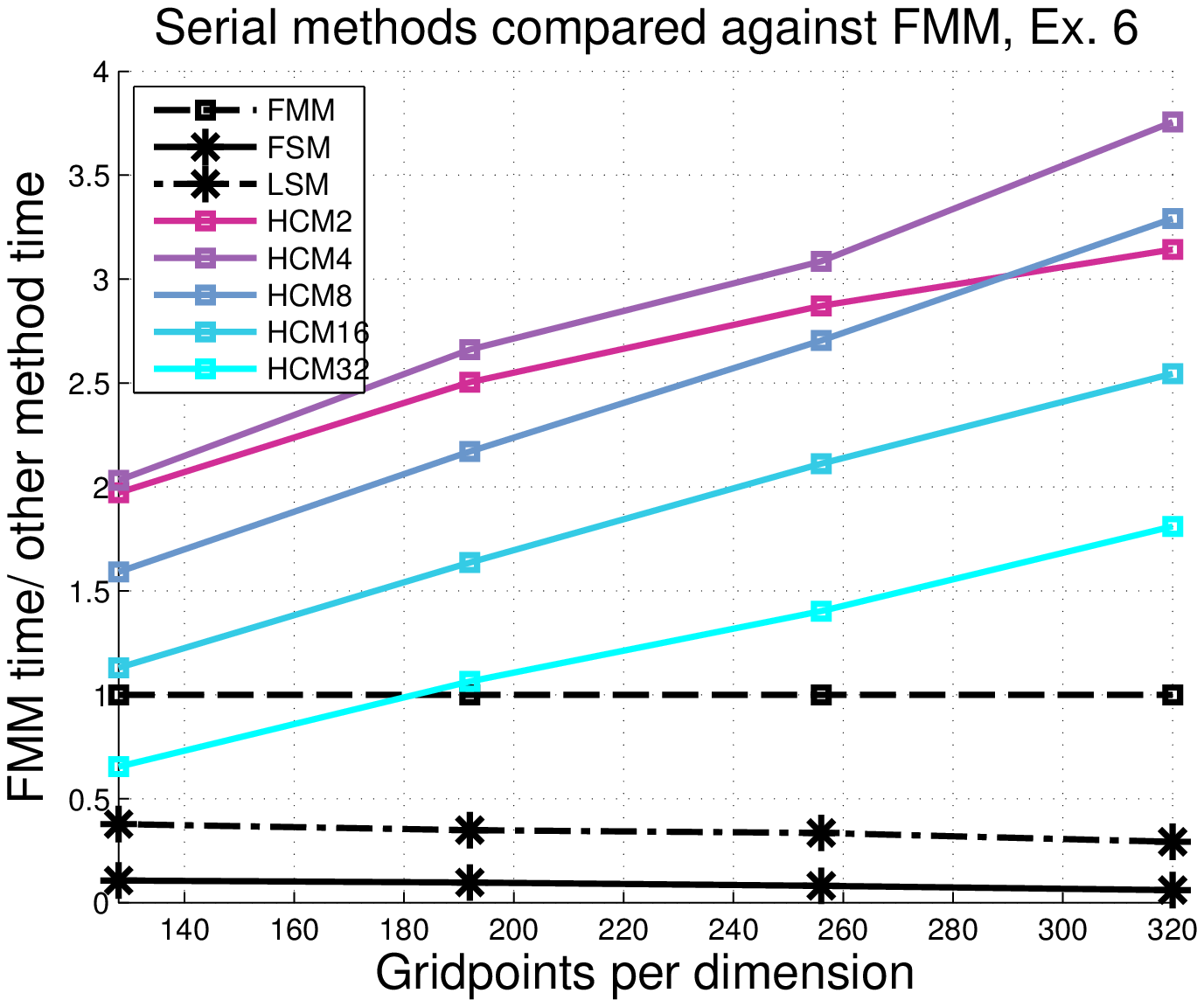}&
\hspace{-.2in}
\includegraphics[scale = .45]{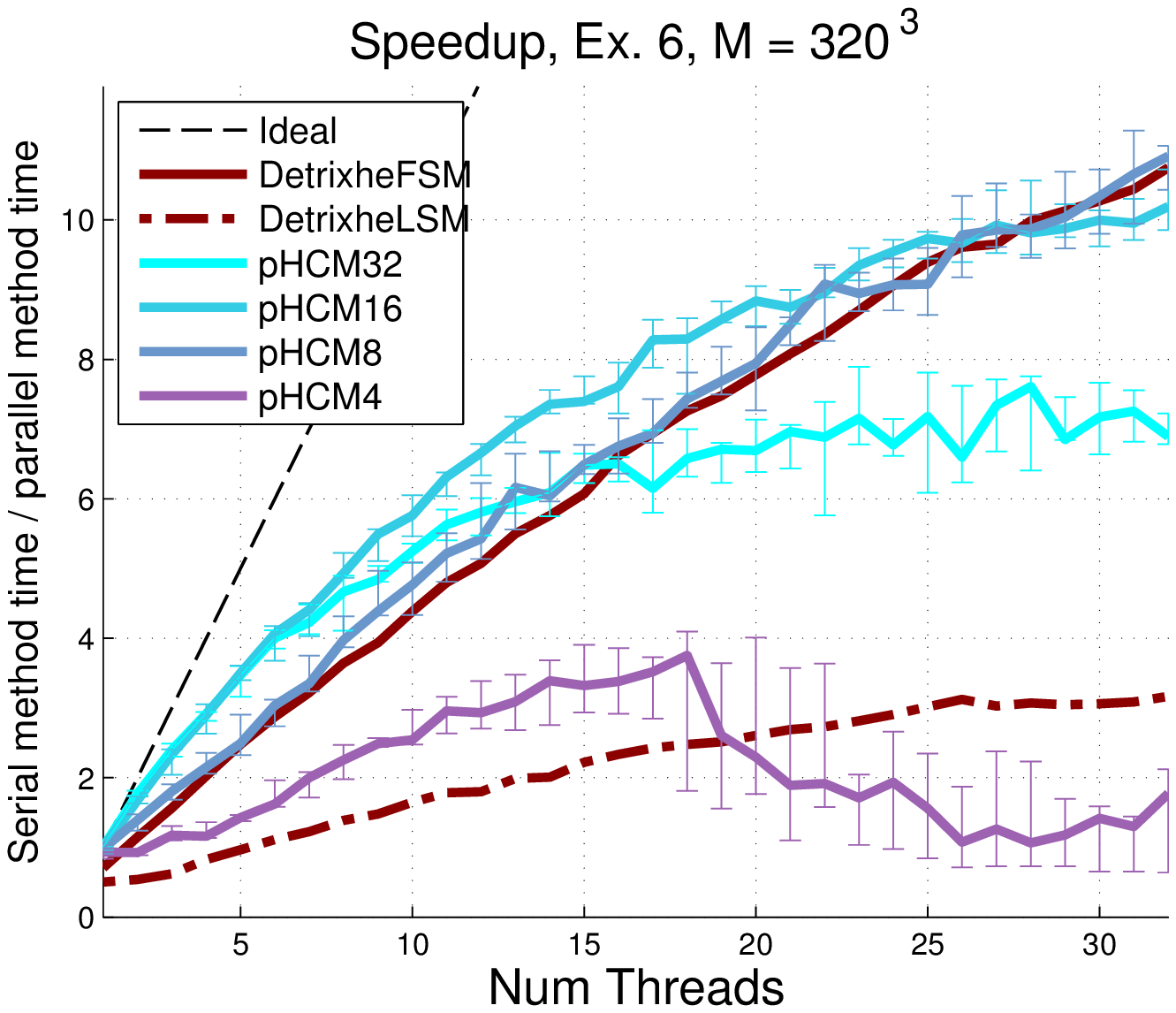}&
\hspace{-.2in}
\includegraphics[scale = .45]{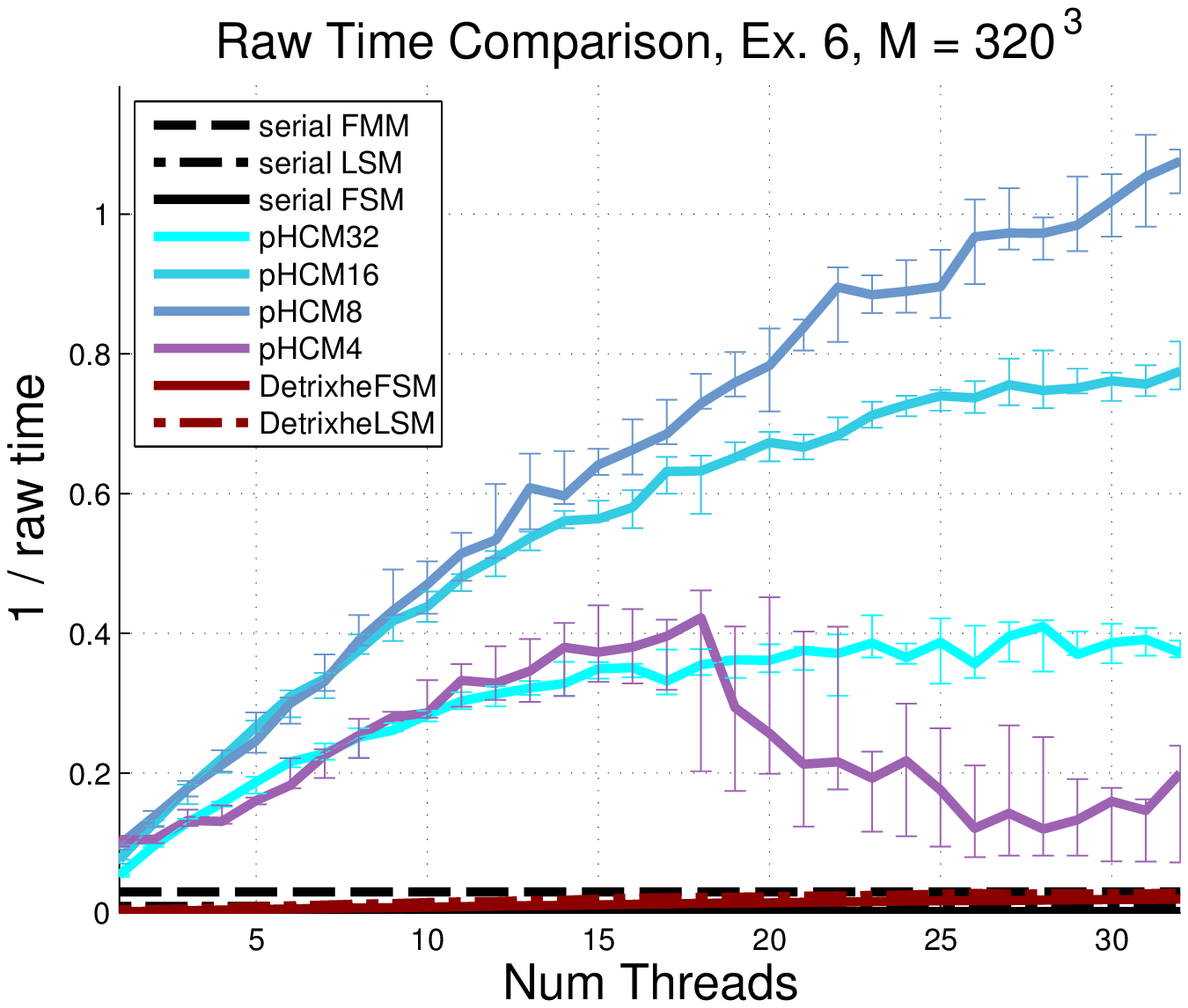}\\
A&B&C\\

\end{array}
$
\caption{\footnotesize Permeable Shell Maze example: serial $M$-scaling comparison ($A$),   parallel scaling at $M = 320^3$ ($B$), and comparison of all methods at $M = 320^3$ ($C$).
}
\label{fig:shellMaze}
\end{figure}

\pagebreak

\subsection{Other cell values}
\label{ss:otherCellVals}

Interestingly, pHCM seems less influenced by the particular choice of cell value heuristic than the serial HCM.  As noted in section \ref{ss:parallel_performance}, if the cell value is a very good predictor of information flow, pHCM will usually see an increase in the total amount of work by not being able to process cells exactly in their causal ordering.  However, pHCM can also partially mitigate the effect of poor cell values; instead of the cell with the lowest value always being processed, we can think of pHCM as simultaneously processing cells in the lowest range of values.  If it is always the case that the true ``most upwind" cell has a value in that range, then pHCM will need fewer heap removals than HCM. Furthermore, neighboring cells that are simultaneously processed may be able to resolve their interdependencies, which would also reduce the total number of heap removals and the number of sweeps per cell (see Figure \ref{fig:work_oldCellVal}A).

We have tested both HCM and pHCM with several other cell value heuristics, including the one from our previous work \cite{ChacVlad}.  We describe it here in Figure \ref{fig:cell_values_old} and equation \eqref{eq:cell_value_old}, supposing $A$ and $B$ are two adjacent cells, with $A$ currently processed.  As before, we define $A_{new} \subset N(B) \bigcap A $ as the set of newly updated inflowing gridpoints of $A$  along the relevant cell border (colored in blue in Figure \ref{fig:cell_values_old}).

\vspace{.4cm}

%COMMENT REMOVED
\begin{figure}[H]
\center{
\begin{tikzpicture}
[scale = 1.8]
\draw[thick] (-2,1) rectangle +(4,-2);
\draw[thick]  (0,1)--(0,-1);

\draw(-1.7,.7) node{$A$};
\draw(1.7,.7) node{$B$};

%COMMENT REMOVED
%COMMENT REMOVED

\draw[red,very thick,->](-.1, .5)--(1,.5);

\draw[dotted](1,1)--(1,-1);

\draw(-.3,.5) node{$V_j^{*}$};

\foreach \x in {-.7, -.5, .1, .3, .7, .9}
	\filldraw(-.1,\x) circle(.03);
\foreach \x in {-.9, -.3, -.1, .5}
	\filldraw[blue](-.1,\x) circle(.03);

\draw(1.2,0) node{$\y_B$};

\end{tikzpicture}
}
\caption{
{\footnotesize
When cell $A$ tags $B$ as downwind, the value computed for $B$ is an approximation to the value of a point along a center axis of $B$; see equation \eqref{eq:cell_value_old}. }}
\label{fig:cell_values_old}
\end{figure}
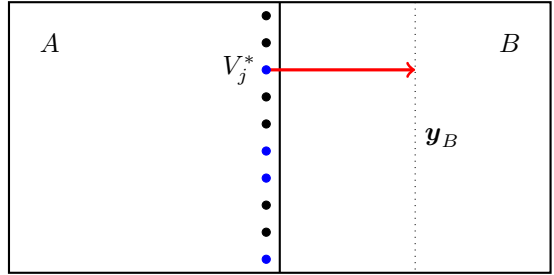

\begin{equation}
\begin{split}
V_{max} \leftarrow \max_{i \in A_{new}} V(\x_i)  \qquad
D \leftarrow \frac{h^c + h}{2} \\\\
\widetilde{V}^c(B) \leftarrow V_{max} + \frac{D}{F(\y)}\\\\
 V^c(B) \leftarrow \min(V^c(B), \widetilde{V}^c(B))
\end{split}
\label{eq:cell_value_old}
\end{equation}
See Figure \ref{fig:cell_values_old} for a geometric interpretation.  For consistency with \cite{ChacVlad}, we tested this heuristic \emph{without} resetting cell values to $+\infty$ each time a cell is processed (see line 5 of Algorithm \ref{alg:HCM_main} and line 11 in Algorithm \ref{alg:pHCM}).  We observed that
\begin{itemize}

\item For serial methods, formula \eqref{eq:cell_value} results in better performance than formula \eqref{eq:cell_value_old} if $r$ is large.

\item For smaller $r$ the median raw time and scaling are better when using \eqref{eq:cell_value_old}.

\item For parallel methods, \eqref{eq:cell_value_old} leads to improved scaling for larger cells.  E.g., Figure \ref{fig:work_oldCellVal}A illustrates how pHCM32 performs noticeably \emph{less} work (measured in terms of AvS) than HCM32, though the raw time actually increases compared to heuristic \eqref{eq:cell_value}.

\end{itemize}

However, the main motivation for using the new cell heuristic \eqref{eq:cell_value} is that formula \eqref{eq:cell_value_old} leads to very bad performance on problems where discontinuities in the speed function are not aligned with cell boundaries.  E.g., for the example of subsection \ref{ss:shellMaze} with $M = 64^3$, HCM8 yields 20.4 average sweeps per cell with formula \eqref{eq:cell_value} compared to 8366 average sweeps per cell with formula \eqref{eq:cell_value_old}.

\begin{figure}[H]
\hspace{-.6 in}
$
\begin{array}{ccc}

\hspace{-.2in}
\includegraphics[scale = .45]{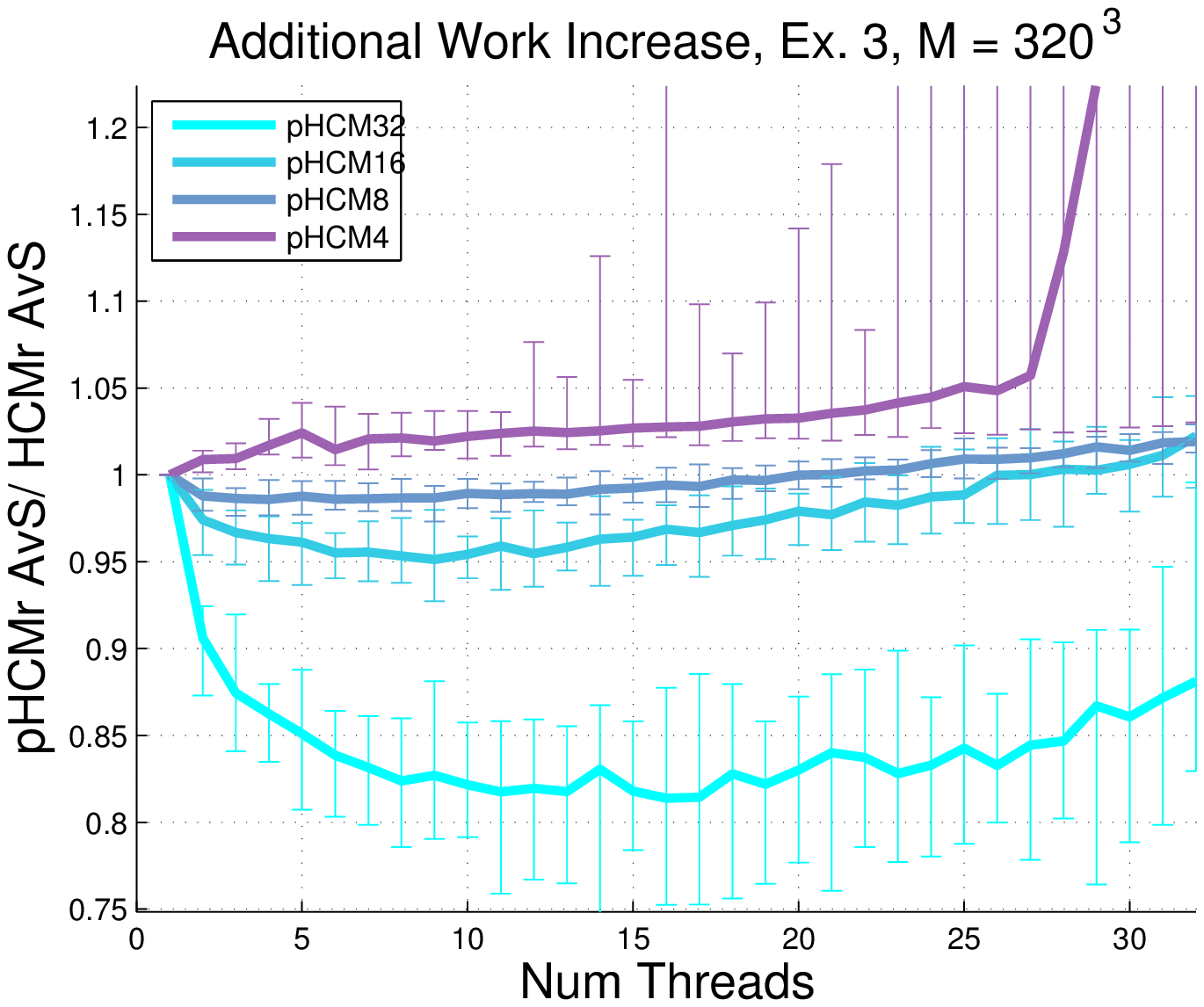}&
\includegraphics[scale = .45]{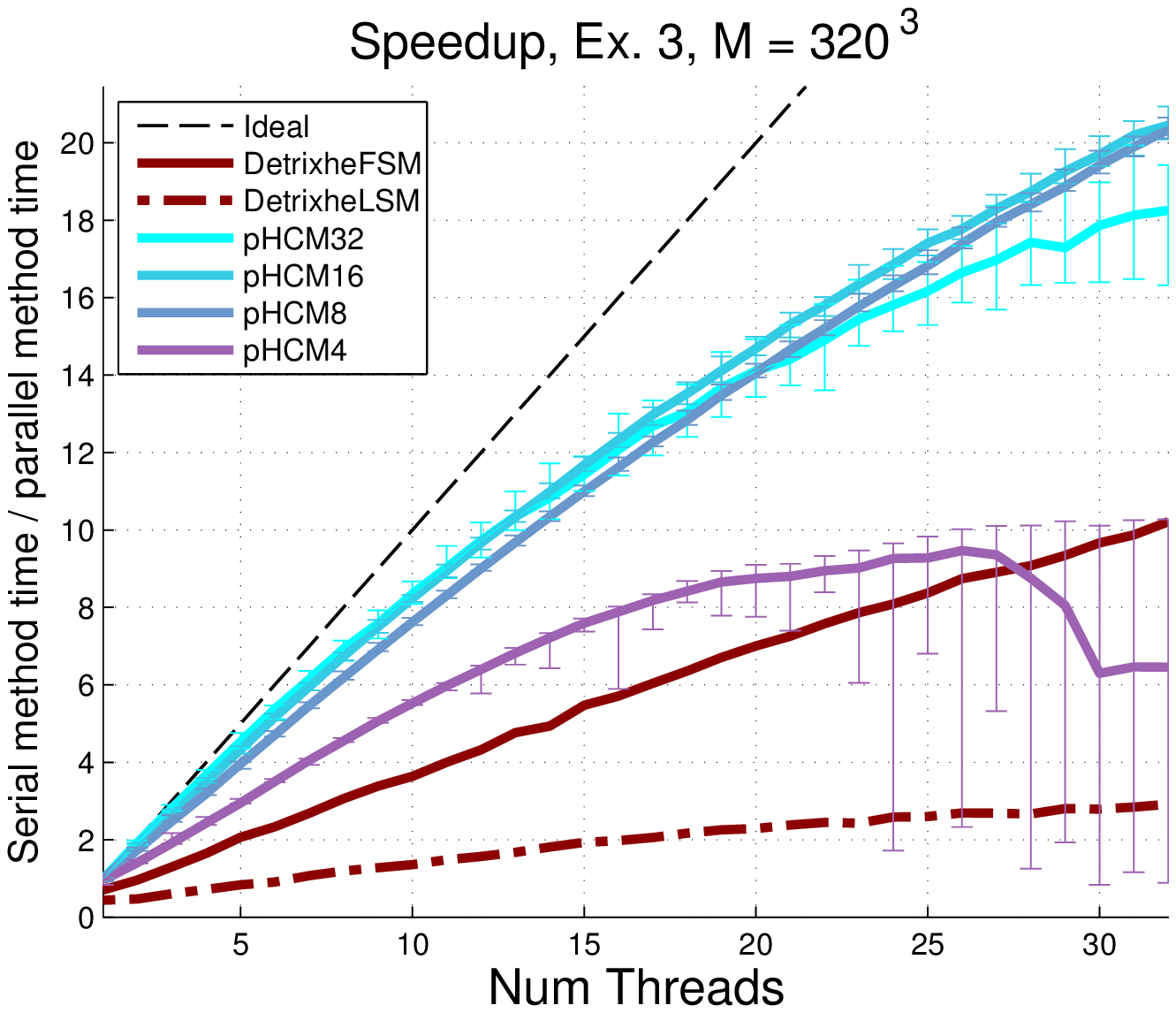}&
\includegraphics[scale = .45]{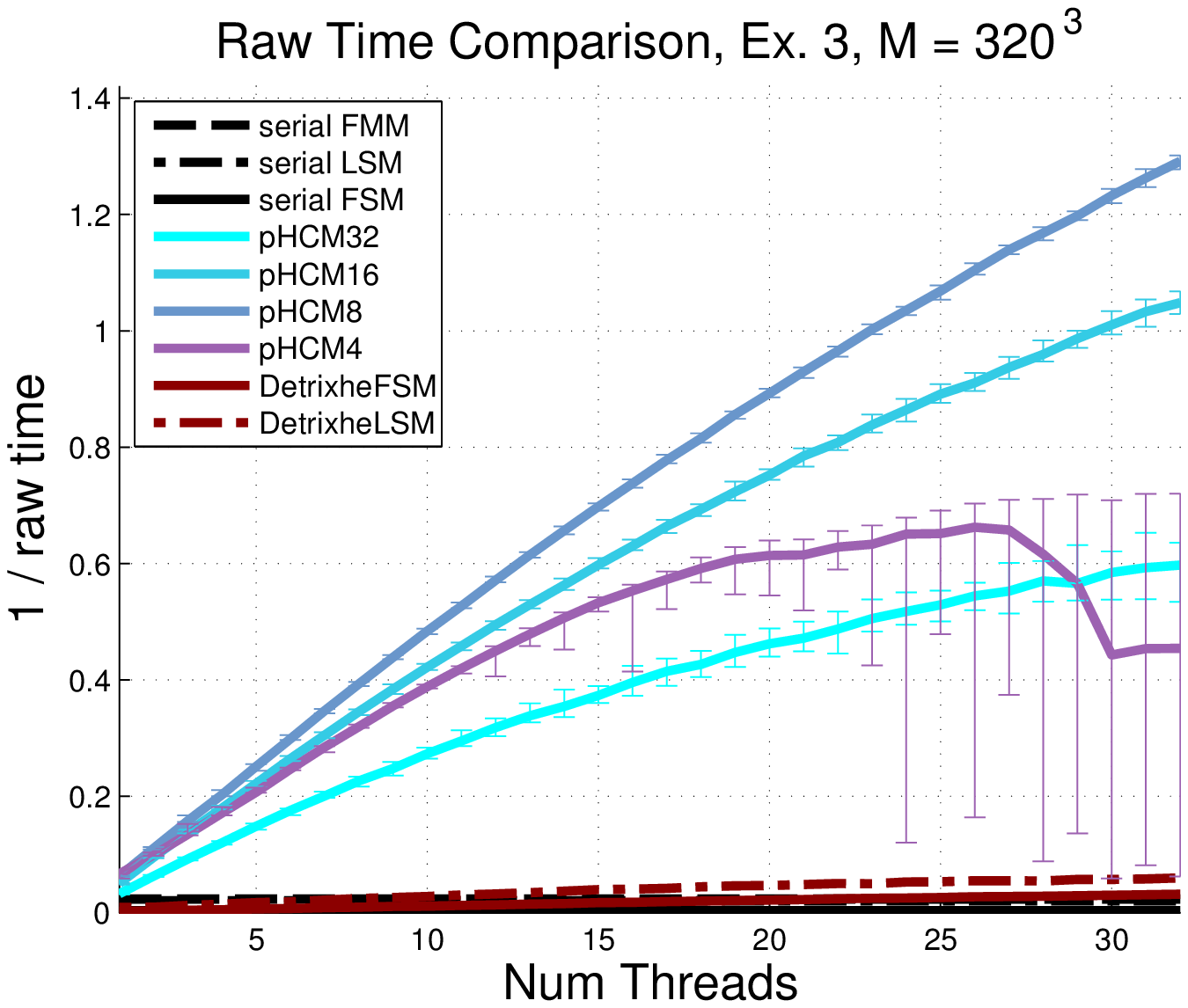}\\
A&B&C\\
\end{array}
$
\caption{\footnotesize An example of pHCM performing less work than HCM for the cell value given by equation \eqref{eq:cell_value_old} on example 3.  Compare with Figures \ref{fig:parallel_charts}C, \ref{fig:parallel_charts}F and \ref{fig:overheads_320}F, and note the difference in scaling in pHCM32.
}
\label{fig:work_oldCellVal}
\end{figure}

%\pagebreak 
\fi

\section{Conclusions}
\label{s:conclusions}
We introduced a new parallel algorithm for the Eikonal equation based on HCM,  a fast two-scale serial solver.  The numerical experiments demonstrated that pHCM achieves its best speedup on problems where the amount of work per cell is high; this occurred when cells were sufficiently large or when the sweeping within cells required more than a few iterations.  As for performance, the combination of HCM's speed and pHCM's good scalability results in a considerable advantage over some of the best serial methods and the parallelization of FSM/LSM.  A comprehensive performance/scaling comparison with other existing parallel Eikonal solvers remains to be performed in the future.

The benchmarking and design of our algorithm was influenced by a particular shared memory architecture, e.g., each thread currently handles the cell-level sweeping serially.  An efficient hybrid GPU/multicore implementation could parallelize %COMMENT REMOVED
the individual cell processing
on a GPU (e.g., as in \cite{WeberDevirBronsBronsKimmel}) while each CPU core would still maintain its own heap.  A possible bottleneck of this approach is the smaller number of GPUs compared to the number of CPU cores in most current systems.
Extensions to a distributed memory architecture appear more problematic since communication times would likely dominate the
%COMMENT REMOVED
cell-processing, at least for the first-order upwind discretization of the Eikonal considered in this paper.
%COMMENT REMOVED

As in HCM, the performance of pHCM for each problem is dependent on a particular cell-decomposition.
E.g., given fixed $P$ and $M$, what value of $J$ will result in the optimal performance?  In this paper we only suggest an answer based on our numerical experiments, but rigorously addressing it will be clearly useful for practitioners.
Ideally, we would like to base (possibly adaptive) cell-decompositions on a posteriori error estimates.  Another interesting direction is the use of non-cubic cells to improve the causal properties of decompositions. %COMMENT REMOVED

The performance analysis in section \ref{s:experiments} suggests a number of possible pHCM improvements.  A smarter memory allocation strategy can be used to increase the spatial and temporal locality of data (particularly in higher dimensional problems).
Rigorous criteria for early sweeping termination would bring additional performance gains to HCM/pHCM (as well as FSM/LSM).  The methods of \cite{Detrixhe} can be substituted in place of LSM within cells, especially for problems with large cell sizes.
In the longer term, we intend to investigate the applicability of our approach to other PDEs and/or discretizations.  Causal problems with a higher amount of work per gridpoint (e.g., discretizations of anisotropic Hamilton-Jacobi) are likely to result in even better pHCM scalability.  We expect this to be also the case for extensions of other parallel Eikonal solvers (e.g., DFSM/DLSM).

Finally, we %COMMENT REMOVED
hope that practitioners will find pHCM useful for applications requiring its efficiency.

\section{Acknowledgements}
\label{s:acknowledgements}
We %COMMENT REMOVED
thank David Bindel, for guidance with parallel computing, and Jeffrey Donatelli for useful correspondence regarding FMM and memory access costs.
We
%COMMENT REMOVED
are also grateful to
Miles Detrixhe for his help in implementing parallel fast sweeping methods and analyzing their performance.  Finally, we
%COMMENT REMOVED
thank the
XSEDE for the computing time allocation and the %COMMENT REMOVED
Texas Advanced Computing Center for the use of their ``Stampede" supercomputer.
%We would also like to thank Miles Detrixhe for his advice in implementing the Detrixhe Sweeping Method and for insight into how different shared memory architectures can lead to vastly different parallel performance. 

\end{document}